\documentclass[11pt,a4paper,twoside]{amsart}

\RequirePackage[colorlinks,citecolor=blue,urlcolor=blue]{hyperref}

\usepackage{enumitem}
\usepackage{graphics}
\usepackage[mathcal]{euscript}
\usepackage{tikz}
\usetikzlibrary{decorations.pathreplacing,decorations.pathmorphing,math}
\usepackage{xpatch}
\usepackage{thmtools}
\usepackage{thm-restate}
\usepackage{xassoccnt}
\usepackage{ytableau}
\usepackage[font=small,labelfont=bf]{caption}
\usepackage{subfig}
\tikzset{snake it/.style={decorate, decoration={random steps,segment length=8pt,amplitude=2.0pt}}}
\usepackage[capitalise,nameinlink]{cleveref}

\makeatletter
\newcommand*{\centerfloat}{%
  \parindent \z@
  \leftskip \z@ \@plus 1fil \@minus \textwidth
  \rightskip\leftskip
  \parfillskip \z@skip}
\makeatother

\makeatletter
\newcounter{globallistcounter}
\newcounter{thmlistcounter}

\let\realItem\item
\NewDocumentCommand\myItem{o}{%
   \IfNoValueTF{#1}{%
      \realItem%
      \edef\countername{\LastRefSteppedCounter}%
      \def\thegloballistcounter{\csname the\countername\endcsname}%
      \refstepcounter{globallistcounter}%
   }{%
      \realItem[#1]%
      \def\thegloballistcounter{#1}%
      \refstepcounter{globallistcounter}%
   }%
}

\AtBeginEnvironment{enumerate}{%
    \let\item\myItem%
    \stepcounter{globallistcounter}%
}

\def\localItemStringLiteral{LOCALITEM@}
\NewDocumentCommand\myLabel{m}{%
      \realLabel{\localItemStringLiteral#1}%
      \setcounter{thmlistcounter}{\value{globallistcounter}}%
      \def\thethmlistcounter{\theparentcounter\thegloballistcounter}%
      \refstepcounter{thmlistcounter}%
      \realLabel[\parentcountername]{#1}
}
\NewDocumentCommand\reflocal{m}{%
    \ref{\localItemStringLiteral#1}%
}

\newlist{thmlist}{enumerate}{1}
\setlist[thmlist]{
    before=\let\item\myItem%
           \let\realLabel\label%
           \let\label\myLabel%
           \stepcounter{globallistcounter}%
           \edef\parentcountername{\LastRefSteppedCounter}%
           \edef\theparentcounter{\csname the\parentcountername\endcsname},
    label=\upshape{(\alph{thmlisti})},
    noitemsep}

\makeatother

\AtEndEnvironment{proof}{\setcounter{claim}{0}}

\declaretheorem[name=Theorem,
                Refname={Theorem,Theorems},
                numberwithin=section,style=plain]{theo}
\declaretheorem[name=Proposition,
               Refname={Proposition,Propositions},
               numberlike=theo,style=plain]{prop}
\declaretheorem[name=Lemma,
               Refname={Lemma,Lemmas},
               numberlike=theo,style=plain]{lemma}
\declaretheorem[name=Conjecture,
               Refname={Conjecture,Conjectures},
               numberlike=theo,style=plain]{conj}
\declaretheorem[name=Definition,
               Refname={Definition,Definitions},
               numberlike=theo,style=remark]{defi}

\declaretheorem[name=Question,
               Refname={Question,Questions},
               numberlike=theo,style=remark]{quest}

\declaretheorem[name=Claim,
               Refname={Claim,Claims},
               style=plain]{claim}

\newcommand{\bbN}{{\mathbb N}}
\newcommand{\bbR}{{\mathbb R}}

\newcommand{\calA}{{\mathcal A}}

\newcommand{\calC}{{\mathcal C}}
\newcommand{\functional}{F}
\newcommand{\calL}{{\mathcal L}}

\newcommand{\calP}{{\mathcal P}}
\newcommand{\calU}{{\mathcal U}}
\newcommand{\calV}{{\mathcal V}}
\newcommand{\bfx}{{\mathbf{x}}}
\newcommand{\bfm}{{\mathbf{m}}}
\newcommand{\bfn}{{\mathbf{n}}}
\newcommand{\bfk}{{\mathbf{k}}}
\newcommand{\bfzero}{{\mathbf{0}}}
\newcommand{\Fmax}{\functional_{\rm max}}

\newcommand{\aas}{a.a.s.}
\newcommand{\Aas}{A.a.s.}

\newcommand{\eps}{\varepsilon}
\newcommand{\decreasing}{D}

\newcommand{\kors}{\#}

\DeclareMathOperator{\supp}{supp}

\DeclareMathOperator{\diam}{diam}
\DeclareMathOperator{\interior}{int}
\DeclareMathOperator{\Diff}{Diff}
\DeclareMathOperator{\Expect}{\mathbb{E}}
\DeclareMathOperator{\dom}{dom}

\providecommand{\abs}[1]{\lvert#1\rvert}
\providecommand{\norm}[1]{\lVert#1\rVert}

\newcommand{\intuition}[1]{\par\noindent{(Intuition: #1)}}
\newcommand{\upsidedown}[1]{\tikz[baseline=(a.north)]\node[yscale=-1,inner sep=0,outer sep=0](a){#1};}
\newcommand{\lagomsqrt}[1]{\sqrt{\vphantom{u_x}\smash[b]{#1}}}

\begin{document}

\title[Monotone subsequences in random permutations]{Monotone subsequences in locally uniform random permutations}
\author{Jonas Sj{\"o}strand}
\address{School of Education, Culture and Communication, Division of Mathematics and Physics, M\"alardalen University, Box~883, 721~23 V\"aster\aa s, Sweden}
\email{jonas.sjostrand@mdu.se}
\date{March, 2023}

\begin{abstract}
A locally uniform random permutation is generated by sampling $n$ points
independently from some absolutely continuous distribution $\rho$
on the plane and interpreting them as a permutation
by the rule that $i$ maps to $j$ if the $i$th point from the left is the $j$th point from below.
As $n$ tends to infinity, decreasing subsequences in the permutation will appear as curves in the plane,
and by interpreting these as level curves, a union of decreasing subsequences give rise to a surface. We show that, under the correct scaling, for any $r\ge0$, the largest union of
$\lfloor r\sqrt{n}\rfloor$ decreasing subsequences approaches a limit surface as $n$ tends to infinity, and the limit surface is a solution to a specific variational problem.
As a corollary, we prove the existence of a limit shape for the Young diagram
associated to the random permutation under the Robinson--Schensted correspondence.
In the special case where $\rho$ is the uniform distribution
on the diamond $\abs{x}+\abs{y}<1$ we conjecture that the limit shape is triangular, and assuming the conjecture is true we find an explicit formula for the limit surfaces of a uniformly random permutation and recover the famous limit shape of Vershik, Kerov and Logan, Shepp.
\end{abstract}

\maketitle

\section{Introduction}\label{sec:introduction}
It has been known since the 1970s that
the longest decreasing (or increasing) subsequence of a random
permutation of $\{1,2,\dotsc,n\}$ has length approximately
$2\sqrt{n}$ for large $n$. More generally, the (scaled) limit
of the cardinality of the largest union of
$\lfloor r\sqrt{n}\rfloor$ disjoint decreasing subsequences is known
for any $r\ge0$, where $\lfloor\cdot\rfloor$ denotes the integral part. 
But what does this union typically look like in the
permutation diagram? And what if the permutation is not sampled
from the uniform distribution? The aim of this paper is to answer
these questions, at least for a family of distributions called \emph{locally uniform}.

Let $\sigma$ be a finite set of points in the plane,
no two of which have the same $x$- or $y$-coordinate.
We can interpret any such $\sigma$ as a permutation
by letting
$\sigma(i)=j$ if the $i$th point from the left is the $j$th point from
below. If $\sigma$ consists of $n$ points that are sampled independently
from some given absolutely continuous distribution $\rho$ on the plane,
$\sigma$ is said to be \emph{locally uniform (with density $\rho$)}.
In particular, if $\rho$ is the uniform distribution on the unit square $(0,1)^2$, then, as a permutation,
$\sigma$ is uniformly distributed among
all permutations of order $n$.

\begin{figure}
\begin{tikzpicture}
\node[xscale=1,yscale=-1,inner sep=0,outer sep=0](0,0){\includegraphics[scale=0.23]{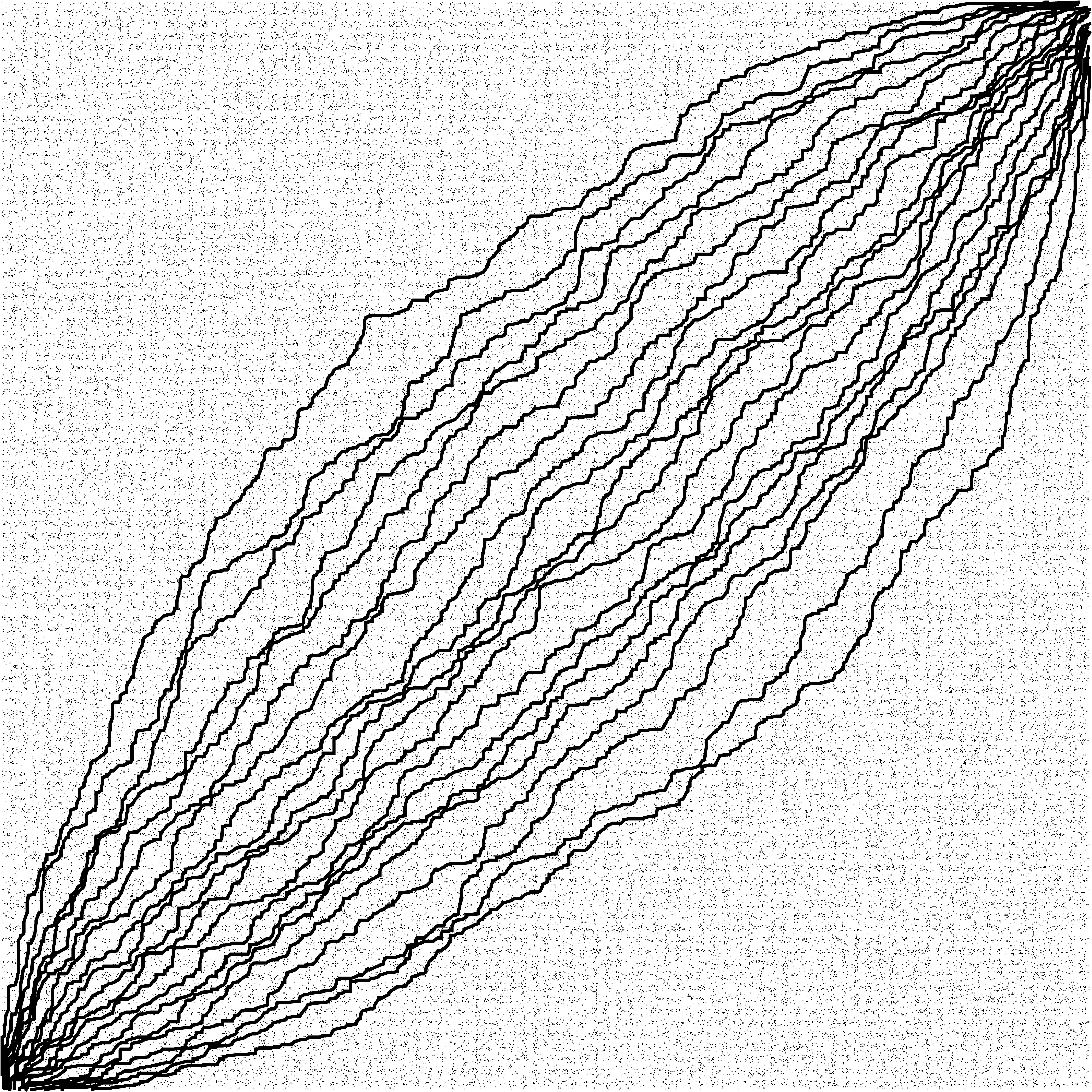}};
\draw[very thick] (2.2,-0.2) -- (3.2,-1.5) -- (2.8,-2.05) -- (1.8,-0.75) -- cycle;
\end{tikzpicture}
\caption{The location of the largest union of 20 decreasing subsequences in
a random permutation of order 100{,}000 drawn from the uniform distribution
on the square $(0,1)^2$.
The small dots are the points in the permutation, and adjacent points in
the decreasing subsequences are connected by line segments. We have also sketched a local parallelogram.}\label{fig:onion}
\end{figure}

In this geometric setting, decreasing subsequences of $\sigma$
appear as ``decreasing subsets'' of the permutation points in the plane,
and we may talk about the \emph{location} of a union of decreasing
subsequences. \Cref{fig:onion} shows an example.
One could imagine a two-dimensional surface whose level curves follow the
decreasing subsets, and as $n$ tends to infinity, under some rescaling
one might hope to obtain a \emph{limit surface} for a maximal
union of $k$ decreasing subsets, where $k$ depends on $n$.
(It is not hard to see that we must require that $k$ grows as $\sqrt{n}$.) This appears to be a new question already for uniform random permutations, and we will motivate below why we think it is a both natural and powerful one.

\section{Background and significance}\label{sec:background}
We will give a brief introduction to the history and current situation
of the research area of monotone subsequences in random permutations.
For a comprehensive review, we refer to Romik~\cite{RomikBook}.

Let $\sigma$ be a \emph{permutation} of order $n$.
A \emph{subsequence} of $\sigma$ is an ordered sequence
$(\sigma(i_1),\sigma(i_2),\dotsc,\sigma(i_k))$ where
$i_1<i_2<\dotsb<i_k$. It is \emph{increasing} if $\sigma(i_1)<\sigma(i_2)<\dotsb<\sigma(i_k)$
and \emph{decreasing} if $\sigma(i_1)>\sigma(i_2)>\dotsb>\sigma(i_k)$.
Let $\mathsf{L}(\sigma)$ denote the length of (any of) the longest increasing
subsequence(s) of $\sigma$, that is
\[
\mathsf{L}(\sigma):=\max\{k\,:\,
\text{there is an increasing subsequence of length $k$}\}.
\]

In 1961 Ulam \cite{Ulam} asked the following question, sometimes known as
{\bf Ulam's problem}:
If $\sigma_n$ is chosen randomly from the uniform distribution of
all permutations of order $n$, what is the expected value of $\mathsf{L}(\sigma_n)$?

About ten years later, Hammersley \cite{Hammersley}
was able to prove that there exists
a limit not only for the expected value $\mathbb{E}\mathsf{L}(\sigma_n)$ but also
for $\mathsf{L}(\sigma_n)$ itself.
\begin{theo}[Hammersley, 1972]\label{th:hammersley}
The limit $\Gamma=\lim_{n\rightarrow\infty}\mathbb{E}\mathsf{L}(\sigma_n)/\sqrt{n}$
exists, and $\mathsf{L}(\sigma_n)/\sqrt{n}\rightarrow\Gamma$ in probability.
\end{theo}
Though simulations
suggested that $\Gamma=2$, this was not proven until 1977
by Vershik and Kerov \cite{VershikKerov}.
They used the fact that $\mathsf{L}(\sigma)$ equals the
length of the first row in the Young diagram corresponding to $\sigma$
under the Robinson--Schensted bijection, and in fact, they were able to describe a
\emph{limit shape} for the Young diagram corresponding to
a random permutation as $n$ grows to infinity. The latter result was also obtained independently by Logan and Shepp \cite{LoganShepp}.

The next break-through happened in 1999, when
Baik, Deift and Johansson \cite{BaikDeiftJohansson}
were able to describe the asymptotic behaviour of $\mathsf{L}(\sigma_n)$ in
detail.
\begin{theo}[Baik--Deift--Johansson 1999]
The random variable
\[
\frac{\mathsf{L}(\sigma_n)-2\sqrt{n}}{n^{1/6}}
\]
converges in distribution to the
Tracy--Widom distribution as $n\rightarrow\infty$.
\end{theo}

\subsection{Where is the longest increasing subsequence?}
Hammersley's approach was to think about a random permutation as a set of
dots randomly and independently positioned in the unit square, and
that setting will be convenient also for us, so let us redefine
our terminology a bit, in accordance with \cref{sec:introduction}.

Let $\sigma$ be a finite set of points in the plane,
no two of which have the same $x$- or $y$-coordinate.
We can interpret any such $\sigma$ as a permutation
by letting
$\sigma(i)=j$ if the $i$th point from the left is the $j$th point from
below.
A subset $I$ of $\sigma$ is \emph{increasing}
if, for any pair of points $(x,y)$ and $(x',y')$ belonging to $I$,
$x<x'$ if and only if $y<y'$. It is \emph{decreasing} if
$x<x'$ if and only if $y>y'$. This corresponds exactly
to the increasing and decreasing subsequences that we defined earlier.

In this framework a natural question arises:
\begin{quest}\label{qu:wherelongest}
For a random permutation $\sigma$ generated by sampling $n$ points uniformly in the unit square, where in the plane does
the longest increasing subsequence typically reside?
\end{quest}
It follows quite easily from Hammersley's work that,
with high probability, all maximal increasing subsets
will be contained in a small region around the diagonal
of the unit square as $n$ tends to infinity. A much stronger
result on the limit distribution of maximal increasing subsets was obtained in a recent paper by Dauvergne and Vir\'ag~\cite{DauvergneVirag21}.

But the new formulation also calls for a generalization:
What if the points in $\sigma$ are sampled from some \emph{non-uniform}
distribution?
Deuschel and Zeitouni \cite{DeuschelZeitouni95} considered this
generalization and were able to describe a limit curve for the
maximal increasing subset when $\sigma$ is a locally uniform random
permutation.

We will concern ourselves with
the following generalized version of \cref{qu:wherelongest}
for locally uniform random permutations.
\begin{quest}\label{qu:whereunion}
Where in the plane does a maximal union of $r\sqrt{n}$
increasing subsets typically reside?
\end{quest}

\subsection{Novelty}
Except for the above-mentioned result on the limit curve
of the longest increasing subsequence by Deuschel and Zeitouni,
the idea to look at not only the cardinality but also the \emph{location}
of monotone subsequences seems to be completely original. As we will see in
\cref{sec:generalpde} below, this approach is potentially
powerful since it enables us to use analytical tools to study
maximal unions of decreasing subsequences in random permutations sampled
from a non-uniform distribution.

There are some results on the \emph{cardinality} of increasing subsequences
for non-uniform random permutations, but to the best of our knowledge
those are all concerned with $q$-analogues of the uniform distribution,
where a permutation $\pi$ is sampled with probability $q^{f(\pi)+f(\pi^{-1})}$
for some classical permutation statistics $f(\cdot)$ like
the number of inversions (the Mallows distribution) \cite{BhatnagarPeled, MuellerStarr} or
the majorant index \cite{Fulman}. In contrast, the family of locally uniform
distributions is much larger; it is uncountably infinite-dimensional.

\subsection{Relation to limit shapes of Young diagrams}
By a theorem of Greene \cite{Greene} (see \cref{pr:curtisgreene} below),
the cardinality of a maximal union of $k$ decreasing subsequences
is encoded in the Young diagram
corresponding to the permutation under the Robinson--Schensted bijection, so
the asymptotic behavior of such cardinalities corresponds to a \emph{limit shape} of a random Young diagram as $n$ tends to infinity. If the random permutation
is drawn from the uniform distribution, the corresponding Young diagram is
drawn from the Plancherel distribution, and its limit shape is the well-known
result by Vershik, Kerov and Logan, Shepp that we mentioned above.
For non-uniform random permutations, however, the limiting behavior of the
corresponding Young diagram is an open problem.
%In fact, potentially a new proof of the
%limit shape formula of Vershik, Kerov and Logan, Shepp might come out as a
%by-product of the project.

Limit shapes of Young diagrams have gained much interest
over the years, and there are results for specific probability
distributions of Young diagrams, often generated by a stochastic process
\cite{JockuschProppShor, Seppalainen98}.
More recent examples include
\cite{ErikssonJonssonSjostrand12,
ErikssonJonssonSjostrand18Bulgarian,
ErikssonJonssonSjostrand18Markov,
ErikssonJonssonSjostrand20}.

\subsection{Relation to permutation limits}
Locally uniform random permutations, our main objects of study, 
appear naturally as limit objects in the sense of
Hoppen et~al.~\cite{HoppenEtAl}.
Their main result is
a definition of convergence for permutation sequences
and an equivalence between such sequences and
(essentially) locally uniform random permutations, and their paper is the
first step towards a theory for permutations analogous to
the emerging theory of limits of graphs created by
Lov\'asz and many coauthors;
see~\cite{LovaszSzegedy} for
an overview.

\section{Terminology and results}\label{sec:results}

\subsection{Probabilistic setting}
Our probability space will always be a simple point process
in the plane viewed as a
random set of points. We will define complex statements about
such random point sets without worrying about the measurability of the
truth-value of the statement.
Typically, such statements will be parameterized by a real number $\gamma$, and
we say that the statement holds \emph{asymptotically almost surely (\aas)
as $\gamma\rightarrow\infty$}
if it is implied by an event which happens with probability tending to one.
To be precise, we make the following definition.
\begin{defi}
Let $(\Omega, \mathcal{F}, P)$ be a probability space and
let $\{E_\gamma\subseteq\Omega\}_{\gamma>0}$ be a collection of outcome sets
indexed by
a parameter $\gamma$. (Note that we do not require the sets to be elements of the $\sigma$-algebra $\mathcal{F}$.)

We say that $E_\gamma$ happens
\emph{asymptotically almost surely (\aas) as $\gamma\rightarrow\infty$}
if there is a family $\{F_\gamma\in\mathcal{F}\}_{\gamma>0}$ such that
$F_\gamma\subseteq E_\gamma$ for any $\gamma$, and $P(F_\gamma)\rightarrow1$.

Also, if $\{X_\gamma\}_{\gamma>0}$ is a family of functions from $\Omega$ to $\bbR$, we say that $X_\gamma\rightarrow x$ \emph{in probability} if, for any $\eps>0$,
the event $\{\omega\in\Omega\ :\ \abs{X_\gamma(\omega)-x}<\eps\}$ happens
\aas\ as $\gamma\rightarrow\infty$.
\end{defi}

\subsection{Increasing and decreasing sets}
A set $P$ of points in $\bbR^2$ is \emph{increasing}
if, for any pair of points $(x,y)$ and $(x',y')$ belonging to $P$,
$x<x'$ if and only if $y<y'$. It is \emph{decreasing} if
$x<x'$ if and only if $y>y'$. It is \emph{$k$-increasing}
(resp.~$k$-decreasing) if
it is a union of $k$ increasing (resp.~decreasing) sets.

\subsection{The local parallelogram}
Consider a situation like that in \cref{fig:onion}, with a
permutation embedded in the unit square and a chosen maximal
union of $k$ decreasing subsets.
Our main idea is to exploit the local uniformity of the random permutation
by zooming in on a small region.
Let us choose the region to have the shape of a narrow parallelogram as depicted in
\cref{fig:onion}
where the long edges are parallel to the decreasing subsets passing by (the curves in the figure) and the short edges have the same slope but with a positive sign.
Then, we might ask what proportion of the permutation points inside the parallelogram are ``picked up'' by the passing curves. Intuitively,
this value is only dependent on the local density and slope of the curves
together with the local density of permutation points. In fact, since
the property of a set being decreasing is invariant under rescaling of the $x$- and $y$-axes,
the question can be reduced to a problem about a narrow rectangle of
45-degree slope.
Given a ``density'' of lines with slope minus one,
what proportion of the permutation points inside the narrow rectangle
are ``picked up'' by the lines?
The following theorem follows directly from \cref{pr:narrowrectangles},
proven in \cref{sec:phi}.
\begin{theo}\label{th:narrowrectanglesresult}
Let $\Omega$ be the open rectangle (depicted in \cref{fig:slopedbetaonerectangle})
\[
0<(x+y)/\sqrt2<1,\ \ 0<(y-x)/\sqrt2<\beta
\]
for some $\beta>0$,
and let $r\ge0$.
For each $\gamma>0$, let $\sigma_\gamma$ be a Poisson point
process in the plane with homogeneous intensity $\gamma$.
Define the random variable $\Lambda^{(\gamma)}$ as the size
of a maximal $\lfloor r\sqrt{\gamma}\rfloor$-decreasing subset
of $\sigma_\gamma\cap\Omega$. Then, as $\gamma$
and $\beta$ tends to infinity,
$\Lambda^{(\gamma)}/\beta\gamma$ converges in $L^1$ to a constant $\Phi(r)$.
\end{theo}
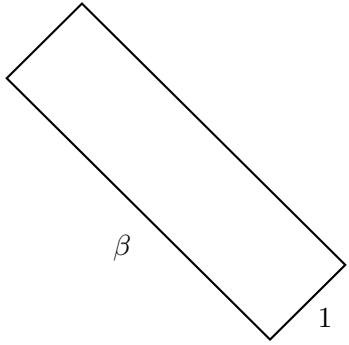
\begin{figure}
\begin{tikzpicture}[scale=0.7]
\draw[thick,rotate around={-45:(0,0)}] (0,0) rectangle (7,2);
\draw [rotate around={-45:(0,0)}]
(7,2) -- (7,0) node [black,midway,right,xshift=0pt,yshift=-6pt]
{1};
\draw [rotate around={-45:(0,0)}]
(7,0) -- (0,0) node [black,midway,below,xshift=-6pt,yshift=-6pt]
{$\beta$};
\end{tikzpicture}
\caption{The rectangle $\Omega$ in \cref{th:narrowrectanglesresult}.}
\label{fig:slopedbetaonerectangle}
\end{figure}

The function $\Phi$ defined by the preceding theorem will play a main role throughout the paper.

\subsection{Limit shape under Robinson--Schensted}
While we have introduced $\Phi$ as a tool for showing our main results below,
it turns out that $\Phi$ is surprisingly interesting in its own right: The derivative of $\Phi$ is a limit shape of
the Young diagram associated with $\sigma$ under the Robinson--Schensted
correspondence, where $\sigma$ is generated by a homogeneous Poisson point process on a diamond square.

\begin{figure}
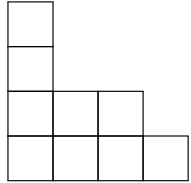

\ydiagram{1,1,3,4}
\caption{A Young diagram $\lambda$ with
column lengths given by $(\lambda_1,\lambda_2,\dotsc)=(4,2,2,1,0,\dotsc)$.}
\label{fig:youngdiagram}
\end{figure}

A \emph{Young diagram} $\lambda$ (in French notation) is a finite collection of
unit cells, arranged in bottom-justified columns whose lenghts are in
non-increasing order from the left. The length of the $i$th column from
the left is denoted by $\lambda_i$, and we let $\lambda_i=0$ if $i$ is
larger than the number of columns.\footnote{Note that this deviates from the standard notation.
In the literature, usually $\lambda_i$ denotes the length of the $i$th longest row, but since we are interested in decreasing rather than increasing subsets,
the columns will be most important to us.}
\Cref{fig:youngdiagram} shows an example.

The \emph{Robinson--Schensted correspondence} is a bijection between permutations and pairs of \emph{standard Young tableaux} of the same shape. We will not define this bijection or even the concept of standard Young tableaux since all we need is contained in the proposition of Greene below. For a comprehensive review we refer to
\cite{StanleyEnum2}.

\begin{defi}
Suppose $\sigma$ is a finite set of points in general position in the
sense that no two points share the same $x$- or $y$-coordinate.
Then we define the \emph{permutation corresponding to $\sigma$} to be the permutation
of $\{1,2,\dotsc,\kors\sigma\}$ defined by letting
$\pi(i)=j$ if the $i$th point from the left in $\sigma$ is
the $j$th point from below. The \emph{Young diagram corresponding to $\sigma$}
is the shape of the standard Young tableaux corresponding to $\pi$ under
Robinson--Schensted.
\end{defi}

In our setting, Greene's \cite{Greene} beautiful connection between the Young diagram and the decreasing (or increasing) subsequences of the permutation
can be formulated as follows.
\begin{prop}[Greene]\label{pr:curtisgreene}
Suppose $\sigma$ is a finite set of points in general position in the
sense that no two points share the same x or y coordinate. For each $k$,
let $\Lambda_k$ be the size of a maximal $k$-decreasing subset of $\sigma$.
Then
\[
\Lambda_k=\sum_{i=1}^k \lambda_i.
\]
\end{prop}

In \cref{pr:phiexists} we will show that $\Phi$ is concave and thus differentiable almost everywhere.
Our next theorem, proven in \cref{sec:mainproof}, states that $\Phi'$ is a
limit shape of the Young diagram corresponding to a homogeneous
Poisson point process on the diamond region $\abs{x}+\abs{y}<1/\sqrt2$.
\begin{restatable*}{theo}{rhombuslimitshapeexists}\label{th:rhombuslimitshapeexists}
Let $\Omega$ be the open diamond square $\abs{x}+\abs{y}<1/\sqrt2$ and,
for each $\gamma>0$, let $\sigma_\gamma$ be a
Poisson point processes on $\Omega$ with
intensity $\gamma$. Then the Young diagram $\lambda^{(\gamma)}$
corresponding to $\sigma_\gamma$ approaches the
limit shape $\Phi'$ in the sense that, for any $r>0$ where $\Phi'(r)$ exists,
\[
\frac1{\sqrt\gamma} \lambda^{(\gamma)}_{\lfloor r\sqrt{\gamma}\rfloor+1}
\rightarrow \Phi'(r)
\]
in probability as $\gamma\rightarrow\infty$.
\end{restatable*}

\begin{figure}
\centerfloat
\begin{tabular}{ccc}
\upsidedown{\includegraphics[height=42mm]{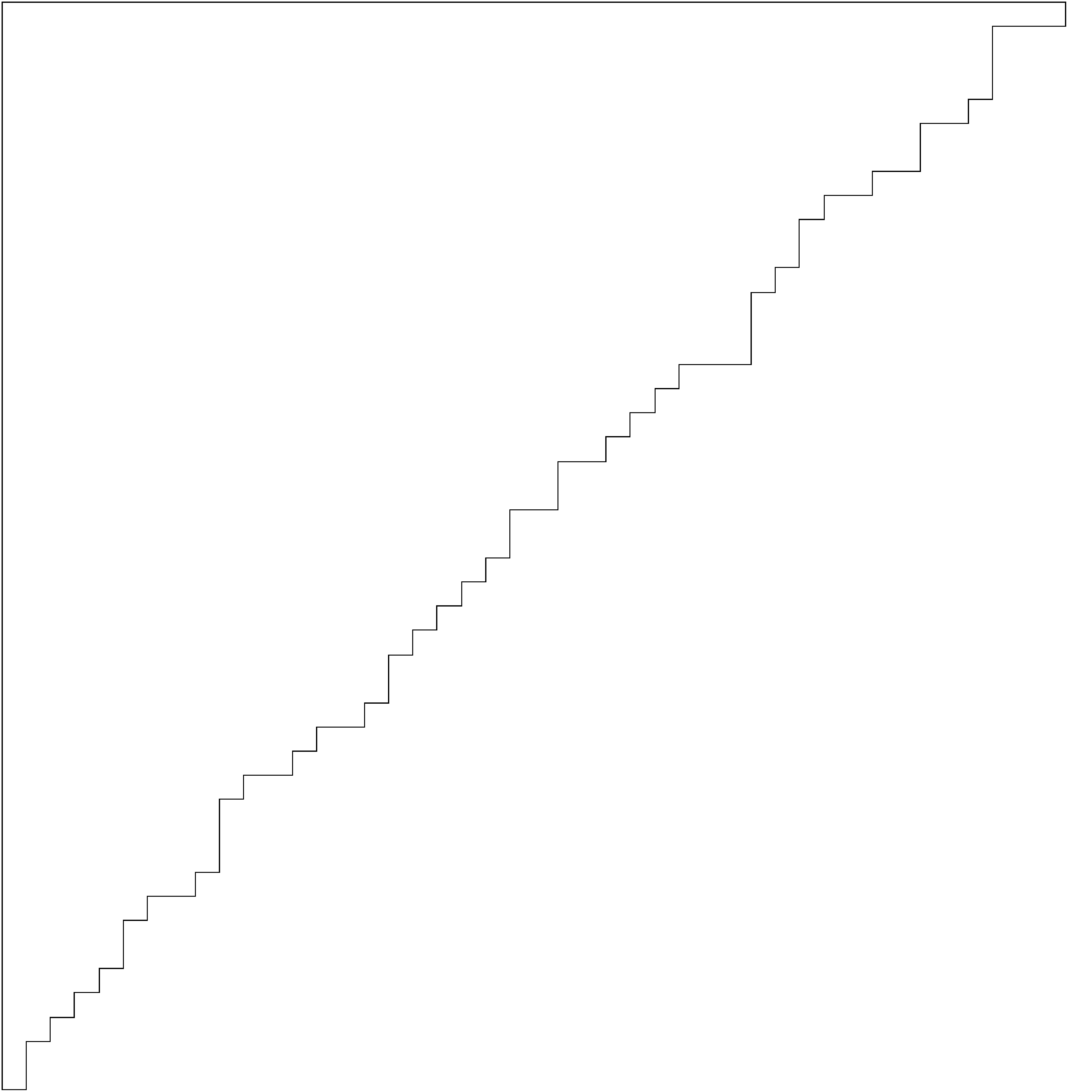}}
&
\upsidedown{\includegraphics[height=42mm]{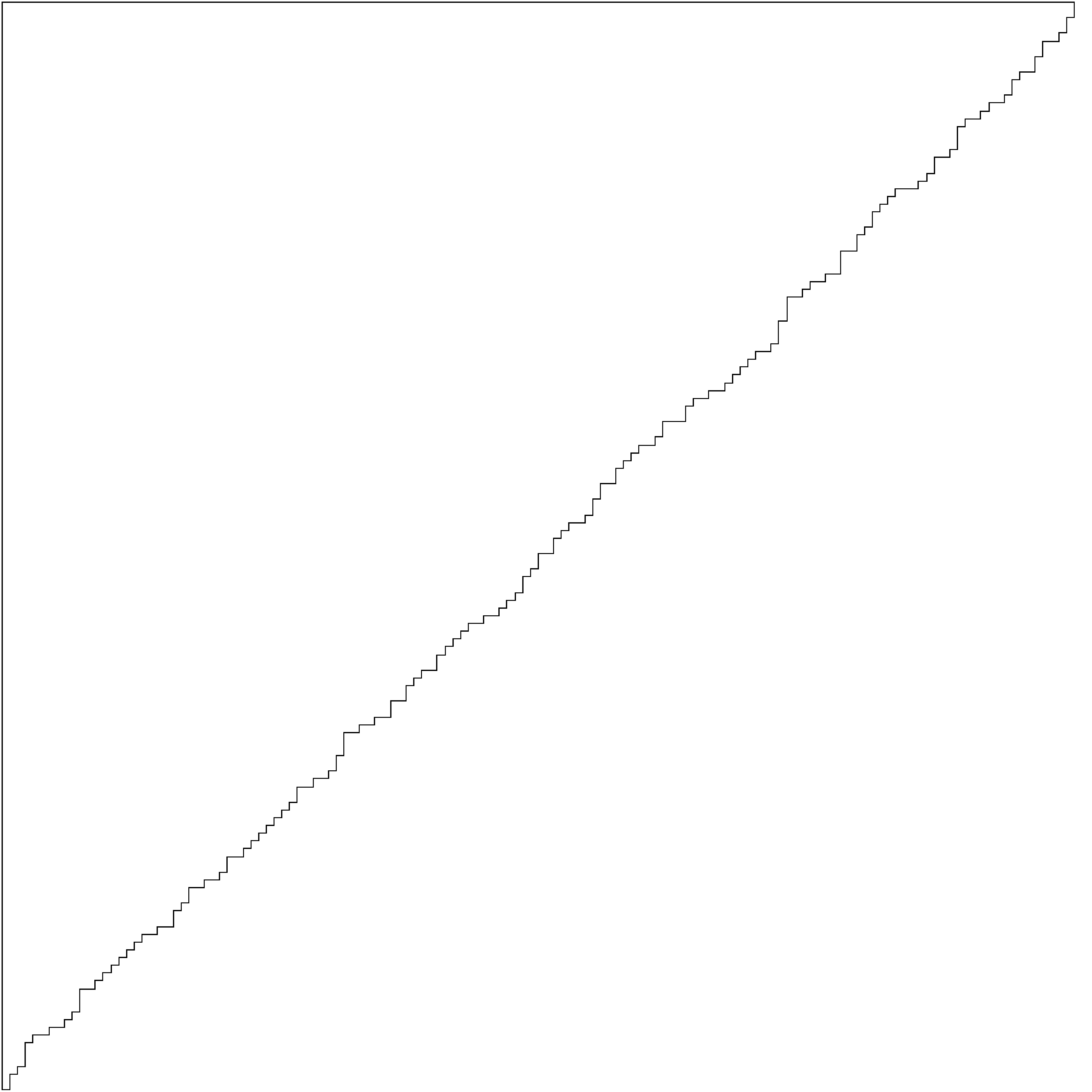}}
&
\upsidedown{\includegraphics[height=42mm]{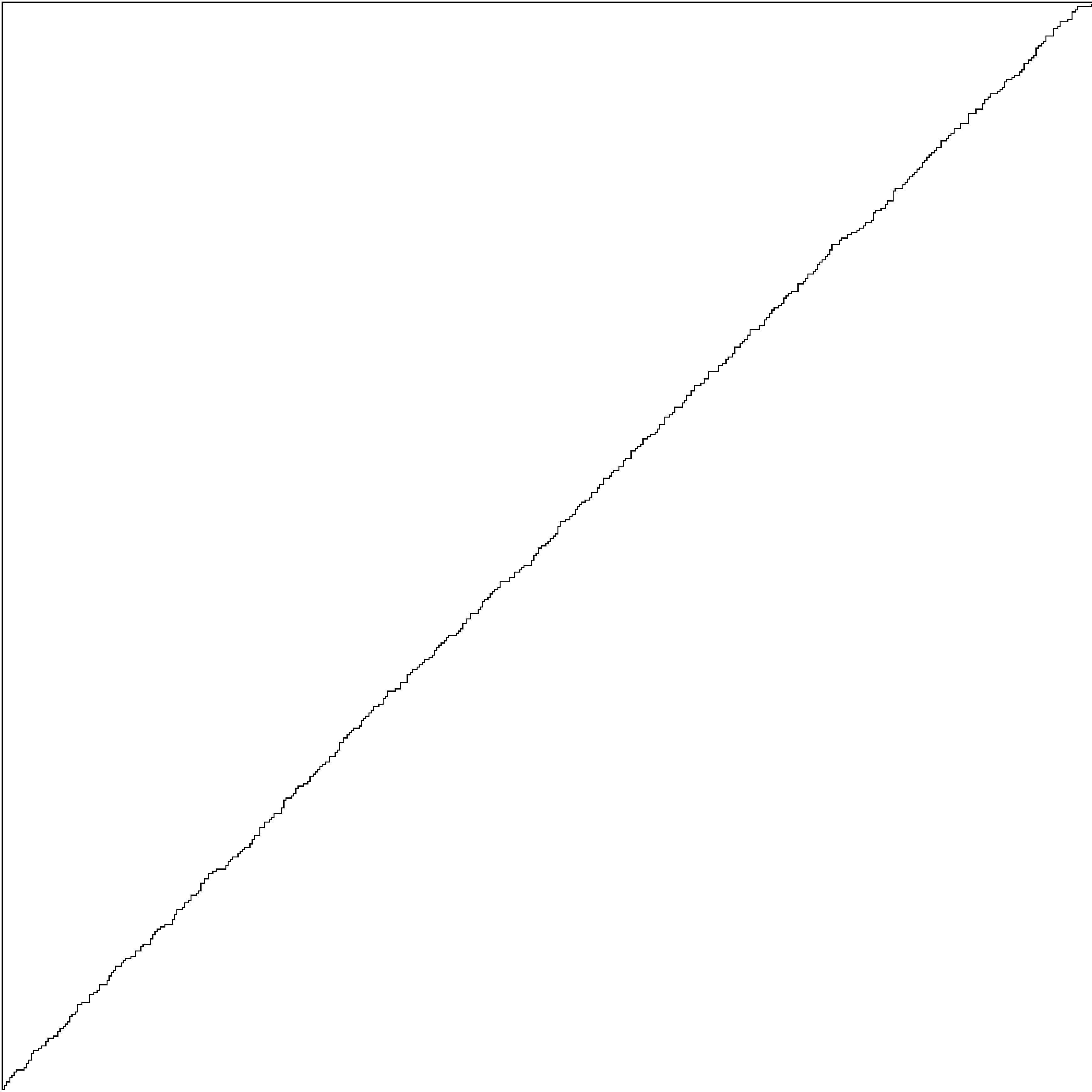}} \\
$n=1000$ & $n=10{,}000$ & $n=100{,}000$
\end{tabular}
\caption{The contours of Young diagrams corresponding
to $n$ points sampled from the uniform distribution on the diamond shape $\abs{x}+\abs{y}<1$.}\label{fig:triangularlimitshape}
\end{figure}

In fact, as is evident in \cref{fig:triangularlimitshape},
computer simulations
strongly suggest that the limit shape is an isosceles triangle!
We make the
following conjecture.
\begin{conj}\label{con:triangularlimitshape}
\[
\Phi'(r)=
\begin{cases}
\sqrt2-r & \text{if $0\le r \le \sqrt2$,} \\
0 & \text{if $r>\sqrt2$.}
\end{cases}
\]
\end{conj}

\subsection{Doubly increasing functions}
As mentioned in the introduction, we want to find some kind of limit object
to a bundle of decreasing subsets like those depicted by the curves in \cref{fig:onion}, and a natural idea is to view the curves as level curves
of a two-dimensional surface. Such surfaces can be described by functions
of $x$ and $y$ that are increasing in both variables.

Define a partial order $\le$
on $\bbR^2$ by letting $(x_1,y_1)\le(x_2,y_2)$ if
$x_1\le x_2$ and $y_1\le y_2$.
For any subset $A$ of $\bbR^2$,
a function $u:A\rightarrow\bbR$ is \emph{doubly increasing} if
$u(x_1,y_1)\le u(x_2,y_2)$ whenever $(x_1,y_1)\le(x_2,y_2)$.
For $r\ge0$, we let $\calU_r(A)$ denote the set of
doubly increasing functions
$u$ on $A$ with $\diam u(A)\le r$, and we let $\calU(A):=\bigcup_{r\ge0}\calU_r(A)$
denote the set of all bounded
doubly increasing functions on $A$.
Let $\calU_{h,r}(A)$ denote the subset of $\calU_r(A)$ consisting of functions with values in $[h,h+r]$.

\begin{figure}
\begin{tikzpicture}[scale=0.8]
\draw (1,8) -- (1,4) -- (3,4) -- (3,1) -- (8,1);
\draw[fill] (1,4) circle(2pt) (3,1) circle(2pt);
\draw (2,8) -- (2,6) -- (4,6) -- (4,3) -- (5,3) -- (5,2) -- (8,2);
\draw[fill] (2,6) circle(2pt) (4,3) circle(2pt) (5,2) circle(2pt);
\draw (6,8) -- (6,5) -- (8,5);
\draw[fill] (6,5) circle(2pt);
\draw (1,2) node {0} (3,5) node {1} (5,5) node {2} (7,6) node {3};
\end{tikzpicture}
\caption{The $\kappa_P$ function for a set $P$ of six points.}
\label{fig:kappa}
\end{figure}
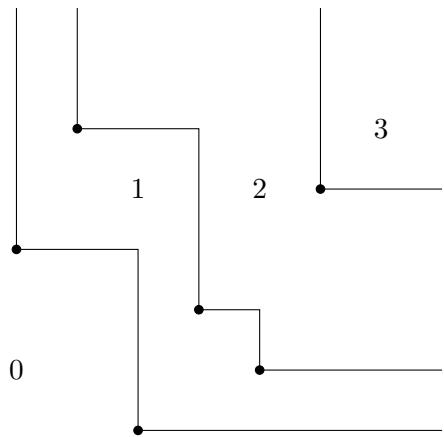

Exactly how should a $k$-decreasing set be interpreted as level curves
and mapped to a doubly increasing function? A simple idea is to convert each
decreasing subset to a curve by joining adjacent points with line segments, but this
is problematic: First, curves from different decreasing subsets might intersect,
and second, there might be multiple ways of partitioning the $k$-decreasing set into
$k$ decreasing subsets.
We can avoid both of these problems by focusing instead on the \emph{increasing} subsets of the $k$-decreasing set, thanks to the
following well-known combinatorial fact (for which we provide a proof
for completeness).
\begin{prop}\label{pr:decreasingincreasingrelation}
Let $P$ be a finite set of points in general position in the
sense that no two points share the same $x$- or $y$-coordinate.
Then, $P$ is $k$-decreasing if and only if it has no increasing subset
of cardinality larger than $k$.
\end{prop}
\begin{proof}
Suppose $P$ is a union of $k$ decreasing sets.
No two elements of an increasing set can belong to the same decreasing set, so, by the pigeonhole principle, there is no increasing subset of $P$ of cardinality larger than $k$.

For the converse, suppose $P$ has no increasing subset with more than $k$
elements. Let $p_1,\dotsc,p_n$ be the points in $P$ sorted from west to east. Construct a sequence of sets $D_1,D_2,\dotsc$ by the following procedure. Initially, let $D_1,D_2,\dotsc$ be
empty sets. Then, iteratively for $i=1,\dotsc n$, add $p_i$ to the first of the sets $D_1,D_2\dotsc$ that currently contains only points to the north-west of $p_i$.
After this procedure, if $p$ is an element in $D_{k+1}$, then $D_k$
must contain an element south-west of $p$; otherwise, $p$
would have been added to $D_k$ instead. Iterating this argument
yields an increasing set of cardinality $k+1$ which contradicts our
assumption. Thus, $D_{k+1}$ is empty and $P$ is a union of the $k$
decreasing sets $D_1,\dotsc,D_k$.
\end{proof}
In accordance, our interpretation of point sets as doubly increasing functions looks as follows.
\begin{defi}\label{def:kappa}
For any finite set $P$ of points in the plane,
define a map $\kappa_P:\bbR^2\rightarrow\bbN$ by letting
$\kappa_P(x,y)$ be the maximal size of an increasing subset
of $P\cap((-\infty,x]\times(-\infty,y])$.
\end{defi}
See \cref{fig:kappa} for an example.

\subsection{A functional}
Now, we will take a global perspective:
Instead of letting the curves
be defined by a maximal union of decreasing subsets, let us think about them just as a bunch of
decreasing curves that we can bend and move freely. The goal is to
position these curves so that together they pick up as many permutation points
as possible. With our parameterization of the bunch of curves by a two-dimensional
surface, the discrete optimization problem can be approximated and formulated as a
continuous variational problem where we want to choose the two-dimensional surface
that maximizes a certain functional.

Let $\mu$ denote the Lebesgue measure on $\bbR^2$.
By a \emph{density domain} we will mean a pair $(\Omega,\rho)$ where
$\Omega$ is an open subset of $\bbR^2$ of positive finite measure
and $\rho$ is a nonnegative function on $\Omega$ such that
$\int_{\Omega}\rho\,d\mu$ is finite. We will often write $\norm{f}_A$
as a shorthand for $\int_A\abs{f}\,d\mu$.

\begin{defi}\label{df:FrhoandL}
For any $\eta,\theta\ge0$, let
\[
L(\eta,\theta):=\begin{cases}
\eta\,\Phi(\sqrt{2\theta/\eta}) & \text{if $\eta>0$,} \\
0 & \text{if $\eta=0$,}
\end{cases}
\]
and, for any density domain $(\Omega,\rho)$, let
$\functional_\rho:\calU(\Omega)\rightarrow\bbR$
be a (nonlinear) functional given by
\[
\functional_\rho(u) := 
\int_\Omega L(\rho, u_xu_y)\,d\mu=\norm{L(\rho, u_xu_y)}_\Omega,
\]
where $u_x$ and $u_y$ denote partial derivatives.
\end{defi}

We will show in \cref{sec:Fwelldefined} that the integrand is integrable so that
$\functional_\rho$ is well defined.

Intuitively, the factor $\Phi(\sqrt{2u_xu_y/\rho})$ of the integrand in the definition of $\functional_\rho$ measures the ``local efficiency'' of
the surface $u$, that is, the proportion of dots in the neighborhood that the curves (encoded by the surface $u$) will
pick up, where we once again refer to \cref{fig:onion}. Note that
it is invariant under rescaling of the $x$- and $y$-axes if the density
$\rho$ is rescaled accordingly.

When the density domain $(\Omega,\rho)$ is implicitly understood,
we let $\Fmax$ be the map from $\bbR_{\ge0}$ to $\bbR_{\ge0}$
defined by letting
\[
\Fmax(r):=\sup_{u\in\calU_r(\Omega)}\functional_{\rho}(u).
\]

\subsection{Main result}
In the introduction we defined a locally uniform random permutation
to consist of a fixed number $n$ of i.i.d.~points with density $\rho$.
However, in many situations it is more natural to consider an (inhomogeneous) Poisson point process with density $n\rho$ where the total number
of points is a (Poisson-distributed) random variable with mean $n$. As $n$ tends to infinity, the difference between these random point
processes becomes negligible; they approach each other in the following sense.
\begin{defi}
Let $\{\sigma_\gamma\}_{\gamma>0}$ and $\{\tau_\gamma\}_{\gamma>0}$
by two families of random point processes parameterized by $\gamma>0$.
Then we say that \emph{$\tau_\gamma$ approaches $\sigma_\gamma$ as $\gamma\rightarrow\infty$} if $\kors(\sigma_\gamma\triangle\tau_\gamma)/\gamma\rightarrow0$ in probability as $\gamma\rightarrow\infty$,
where $\triangle$ denotes symmetric difference.
\end{defi}

Our first main theorem connects the doubly increasing functions and
$\functional_\rho$ with random point processes. We postpone its proof
until \cref{sec:mainproof}.
\begin{restatable*}{theo}{main}\label{th:main}
Let $(\Omega,\rho)$ be a density domain and
let $\{\tau_\gamma\}_{\gamma>0}$ be random point processes on $\Omega$ approaching a Poisson point process with
intensity function $\gamma\rho$ as $\gamma\rightarrow\infty$.
Then the following holds for any $r\ge0$.
\begin{thmlist}
\item\label{th:mainlimitsurface}
For any $\eps>0$, \aas\ as
$\gamma\rightarrow\infty$, for
every maximal $\lfloor r\sqrt{\gamma}\rfloor$-decreasing subset $P$
of $\tau_\gamma$ there is a $u\in\calU_{0,r}(\Omega)$ with $\functional_\rho(u)=\Fmax(r)$ such that
$\norm{\kappa_P/\sqrt{\gamma}-u}_\Omega<\eps$.
\item\label{th:mainlimittoFmax}
Let $\Lambda^{(\gamma)}$ denote the size of a maximal
$\lfloor r\sqrt{\gamma}\rfloor$-decreasing subset of $\tau_\gamma$. Then,
$\Lambda^{(\gamma)}/\gamma\rightarrow \Fmax(r)$ in probability as $\gamma$ tends to infinity.
\end{thmlist}
\end{restatable*}

As a corollary, we obtain a limit shape for
the Young diagram associated with a locally uniform random permutation. 
\begin{restatable*}{coro}{limitshape}\label{cor:limitshape}
With the same setup as in \cref{th:main},
the Young diagram $\lambda^{(\gamma)}$ corresponding to $\tau_\gamma$
approaches the limit shape $\Fmax'$
in the sense that,
for any $r>0$ where the derivative $\Fmax'(r)$ exists,
\[
\frac1{\sqrt\gamma} \lambda^{(\gamma)}_{\lfloor r\sqrt{\gamma}\rfloor+1}
\rightarrow \Fmax'(r)
\]
in probability as $\gamma\rightarrow\infty$.
\end{restatable*}
The proof will appear in \cref{sec:mainproof}.

Our second main theorem is nonprobabilistic in nature and deals with
the functional $\functional_\rho$ exclusively.
We postpone its proof until \cref{sec:mainFproof}.
\begin{restatable*}{theo}{mainF}
\label{th:mainF}
Let $(\Omega,\rho)$ be a density domain. Then the following holds.
\begin{thmlist}
\item\label{th:mainFmaximizer}
For any $r\ge0$, $\functional_\rho$ attains its maximum on $\calU_r(\Omega)$.
\item\label{th:mainFrhoconcave}
$\functional_\rho$ is a concave function on $\calU(\Omega)$.
\item\label{th:mainFmaxcontinuousincreasingconcave}
$\Fmax$ is continuous, increasing and concave, and $\Fmax(r)\rightarrow\norm{\rho}_\Omega$ as $r\rightarrow\infty$.
\end{thmlist}
\end{restatable*}
In \cref{sec:generalpde}, we show how a maximizer of
$\functional_\rho$ can be found in practice by solving a system of partial differential equations.

Note that \cref{th:mainFmaxcontinuousincreasingconcave} means that there is no
loss of mass in the limit shape of \cref{cor:limitshape}, that is, there
is no macroscopic proportion of the boxes in the longest $o(\sqrt{\gamma})$ rows or columns of the Young diagram.

\subsection{Consequences of \texorpdfstring{\cref{con:triangularlimitshape}}{the conjecture}}
Our final result is concerned with the consequences of \cref{con:triangularlimitshape}.
In \cref{sec:uniform}, provided \cref{con:triangularlimitshape} holds true, we obtain a simple
parameterization of the limit surface for a uniformly random permutation,
and we recover the celebrated limit-shape result of Logan, Shepp and Vershik, Kerov mentioned in \cref{sec:background}.

\subsection{Organization of the paper}
The remainder of the paper is organized as follows.
\begin{description}
\item[\cref{sec:phi}]
We show that $\Phi$ exists and that the functional $\functional_\rho$ is
well defined.
\item[\cref{sec:generalpde}]
We reduce the maximization problem to a PDE system.
\item[\cref{sec:localparallelogram}]
We study the probabilistic behavior of the local parallelogram.
\item[\cref{sec:greatlemma}]
We divide a density domain into small parallelograms.
\item[\cref{sec:semicontinuity}]
We show that $\functional_\rho$ is semicontinuous.
\item[\cref{sec:compactness}]
We show that the set $\calU_{h,r}(\Omega)$ of doubly increasing functions
is compact as a subset of $L^1(\Omega)$.
\item[\cref{sec:mainproof}]
We prove our first main theorem, \cref{th:main}.
\item[\cref{sec:mainFproof}]
We prove our second main theorem, \cref{th:mainF}.
\item[\cref{sec:uniquemaximizers}]
We show that, provided $\Phi$ is reasonably well behaved (which it is
if \cref{con:triangularlimitshape} holds true),
the maximizer of $\functional_\rho$ is essentially unique.
\item[\cref{sec:uniform}]
Provided \cref{con:triangularlimitshape} holds true, we find
the limit surfaces for a uniformly random permutation and the
limit shape for the corresponding Young diagram.
\item[\cref{sec:future}]
We discuss some open questions for future research.
\end{description}

\section{Existence of \texorpdfstring{$\Phi$}{Phi}}\label{sec:phi}
\noindent
In this section we will prove \cref{th:narrowrectanglesresult} and thereby establish the existence of the function $\Phi$. We will also
show some basic properties of $\Phi$ and $L$ and see that the functional $\functional_\rho$ defined in \cref{sec:results} is well defined.

Let $\bbN:=\{0,1,2,\dotsc\}$ denote the set of nonnegative integers.
We will use bold letters like $\bfn$ for points in the
integer lattice $\bbN^2$, and coordinates
will be denoted by the corresponding italic letters with indices,
like $\bfn=(n_1,n_2)$. Define the operators plus, minus and the
relations $\le$ and $<$ coordinate-wise on $\bbN^2$. Also, 
let $\bfm\ast\bfn=(m_1n_1,m_2n_2)$ denote coordinate-wise multiplication.
Finally, we let $\bfn\rightarrow\infty$ mean that
$\min\{n_1,n_2\}\rightarrow\infty$.

Hammersley's proof of the existence of a limit for the size of the longest
increasing subsequence (\cref{th:hammersley}) uses
Kingman's subergodic theorem \cite{Kingman}.
We will need a fancier version of that theorem in order
to show the existence of a limiting behavior of the decreasing subsets
in our two-dimensional setting.
The following theorem is a special case of a result by Sch\"urger
\cite[Theorem~2.1]{Schurger88}\footnote{We put $h=0$ in Sch\"urger's theorem
and change the sign of the random variables.}.
\begin{theo}[Sch\"urger]\label{th:schurger}
Let $\{X_{\bfm,\bfn}\}$ be a family of real random variables,
where the indices span over all $\bfm,\bfn\in\bbN^2$ with
$\bfm<\bfn$. Suppose the following holds.
\begin{itemize}
\item \emph{Translation invariance.}
For any $\bfk\in\bbN^2$, the family
$\{X_{\bfm+\bfk,\bfn+\bfk}\}_{\bfm<\bfn}$ has the same finite
joint probability
distributions as $\{X_{\bfm,\bfn}\}_{\bfm<\bfn}$.
\item \emph{Superadditivity.}
For any $\bfzero<\bfm<\bfn$, we have
\begin{align*}
X_{\bfzero,\bfn} &\ge X_{\bfzero,(m_1,n_2)}+X_{(m_1,0),\bfn}, \\
X_{\bfzero,\bfn} &\ge X_{\bfzero,(n_1,m_2)}+X_{(0,m_2),\bfn}.
\end{align*}
\item \emph{Integrability.}
The set $\{\Expect X_{\bfzero,\bfn}/n_1n_2\ :\ \bfn>\bfzero\}$
is bounded.
\end{itemize}
Then, $\lim_{\bfn\rightarrow\infty}X_{\bfzero,\bfn}/n_1n_2$
exists in $L^1$ and equals
\[
\lim_{\bfn\rightarrow\infty}(n_1n_2)^{-1}
\lim_{\bfm\rightarrow\infty}(m_1m_2)^{-1}
\sum_{\bfzero<\bfk\le\bfm}X_{\bfk\ast\bfn-\bfn,\bfk\ast\bfn},
\]
both limits existing in $L^1$.
\end{theo}

If $I$ and $J$ are sets of points, we say that
$I$ is a \emph{$k$-decreasing set compatible with $J$} if
$I\cup J$ is $k$-decreasing.

\begin{prop}\label{pr:X}
Fix $s>0$ and
let $\sigma$ be a Poisson point process on $\bbR^2$ with
homogeneous intensity $s$.

For any pair $(\bfm,\bfn)$ with $\bfm,\bfn\in\bbN^2$ and $\bfm<\bfn$
we define the interval
\[
[\bfm,\bfn):=\{\bfx\in\bbR^2\ :\ \bfm\le\bfx<\bfn\}.
\]
Let $T:\bbR^2\rightarrow\bbR^2$ be the linear transformation 
that maps $(1,0)$ to $(1,1)$ and $(0,1)$ to $(-1,1)$, and
define the random variable $X_{\bfm,\bfn}$ as the maximal size of
an $(n_1-m_1)$-decreasing subset of $\sigma\cap T[\bfm,\bfn)$
compatible with $T([m_1,n_1)\times\{m_2,n_2\})$, where
$[m_1,n_1)$ denotes the set $\{m_1,m_1+1,\dotsc,n_1-1\}$.

Then, $X_{\bfzero,\bfn}/n_1n_2$ converges in $L^1$ to a constant
$c_s$ as $\bfn\rightarrow\infty$.
\end{prop}
\begin{proof}
We claim that the family $\{X_{\bfm,\bfn}\}_{\bfm<\bfn}$
has the three properties defined in \cref{th:schurger}.
Translation invariance follows immediately from the translation
invariance of $\sigma$. Integrability holds since
$\Expect X_{\bfzero,\bfn}/n_1n_2
\le\Expect \#(\sigma\cap T[\bfzero,\bfn))/n_1n_2=2s$.

To prove superadditivity we consider any $\bfzero<\bfm<\bfn$ and
verify the two superadditivity inequalities.

\smallskip
Let $A$ be an $m_1$-decreasing subset of
$\sigma\cap T[\bfzero,(m_1,n_2))$ compatible with
$T([0,m_1)\times\{0,n_2\})$
and let $B$ be an $(n_1-m_1)$-decreasing subset of
$\sigma\cap T[(m_1,0),\bfn)$ compatible with
$T([m_1,n_1)\times\{0,n_2\})$, as depicted in
\cref{fig:superadditivityone}.
\begin{figure}%
\centerfloat%
\subfloat[][]{%
\begin{tikzpicture}[scale=1.1]
\draw[thick,rotate around={45:(0,0)}] (0,0) rectangle (3,3)
(1.2,0) -- (1.2,3);
\draw[fill=black,rotate around={45:(0,0)}] (0,0) circle (2pt) (0.3,0) circle (1pt) (0.6,0) circle (1pt) (0.9,0) circle (1pt) (1.2,0) circle (2pt) (1.5,0) circle (1pt) (1.8,0) circle (1pt) (2.1,0) circle (1pt) (2.4,0) circle (1pt) (2.7,0) circle (1pt) (3,0) circle (2pt);
\draw[fill=black,rotate around={45:(0,0)}] (0,3) circle (2pt) (0.3,3) circle (1pt) (0.6,3) circle (1pt) (0.9,3) circle (1pt) (1.2,3) circle (2pt) (1.5,3) circle (1pt) (1.8,3) circle (1pt) (2.1,3) circle (1pt) (2.4,3) circle (1pt) (2.7,3) circle (1pt) (3,3) circle (2pt);
\pgfmathsetseed{1}
\draw[rotate around={45:(0,0)},snake it]
	(0,0) -- (0.1,0.5) -- (0.1,2.5) -- (0,3)
	(0.3,0) -- (0.38,0.5) -- (0.38,2.5) -- (0.3,3)
	(0.6,0) -- (0.66,0.5) -- (0.66,2.5) -- (0.6,3)
	(0.9,0) -- (0.94,0.5) -- (0.94,2.5) -- (0.9,3)
	(1.2,0) -- (1.3,0.5) -- (1.3,2.5) -- (1.2,3)
	(1.5,0) -- (1.58,0.5) -- (1.58,2.5) -- (1.5,3)
	(1.8,0) -- (1.86,0.5) -- (1.86,2.5) -- (1.8,3)
	(2.1,0) -- (2.14,0.5) -- (2.14,2.5) -- (2.1,3)
	(2.4,0) -- (2.42,0.5) -- (2.42,2.5) -- (2.4,3)
	(2.7,0) -- (2.70,0.5) -- (2.70,2.5) -- (2.7,3);
\draw[white,fill=white,rotate around={45:(0,0)}]
	(0.6,1.5) circle (10pt)
	(2.1,1.5) circle (10pt);
\draw[rotate around={45:(0,0)}]
	(-0.2,-0.2) node {$\scriptstyle T(0,0)$}
	(1.4,-0.6) node {$\scriptstyle T(m_1,0)$}
	(3.5,-0.5) node {$\scriptstyle T(n_1,0)$}
	(-0.5,3.5) node {$\scriptstyle T(0,n_2)$}
	(1.0, 3.6) node {$\scriptstyle T(m_1,n_2)$}
	(3.2,3.2) node {$\scriptstyle T(n_1,n_2)$}
	(0.6,1.5) node {$A$}
	(2.1,1.5) node {$B$}; 
\end{tikzpicture}%
\label{fig:superadditivityone}%
}%
\hfill%
\subfloat[][]{%
\begin{tikzpicture}[scale=1.1]
\draw[thick,rotate around={45:(0,0)}] (0,0) rectangle (3,3)
(0,1.3) -- (3,1.3);
\draw[fill=black,rotate around={45:(0,0)}] (0,0) circle (2pt) (0.3,0) circle (1pt) (0.6,0) circle (1pt) (0.9,0) circle (1pt) (1.2,0) circle (2pt) (1.5,0) circle (1pt) (1.8,0) circle (1pt) (2.1,0) circle (1pt) (2.4,0) circle (1pt) (2.7,0) circle (1pt) (3,0) circle (2pt)
(0,1.3) circle (2pt) (0.3,1.3) circle (1pt) (0.6,1.3) circle (1pt) (0.9,1.3) circle (1pt) (1.2,1.3) circle (1pt) (1.5,1.3) circle (1pt) (1.8,1.3) circle (1pt) (2.1,1.3) circle (1pt) (2.4,1.3) circle (1pt) (2.7,1.3) circle (1pt) (3,1.3) circle (2pt)
(0,3) circle (2pt) (0.3,3) circle (1pt) (0.6,3) circle (1pt) (0.9,3) circle (1pt) (1.2,3) circle (1pt) (1.5,3) circle (1pt) (1.8,3) circle (1pt) (2.1,3) circle (1pt) (2.4,3) circle (1pt) (2.7,3) circle (1pt) (3,3) circle (2pt);
\pgfmathsetseed{1}
\draw[rotate around={45:(0,0)},snake it]
	(0,0) -- (0.07,0.4) -- (0.07,0.9) -- (0,1.3) -- (0.07,1.8) -- (0.07,2.6) -- (0,3)
	(0.3,0) -- (0.3,1.3) -- (0.3,3)
	(0.6,0) -- (0.6,1.3) -- (0.6,3)
	(0.9,0) -- (0.9,1.3) -- (0.9,3)
	(1.2,0) -- (1.2,1.3) -- (1.2,3)
	(1.5,0) -- (1.5,1.3) -- (1.5,3)
	(1.8,0) -- (1.8,1.3) -- (1.8,3)
	(2.1,0) -- (2.1,1.3) -- (2.1,3)
	(2.4,0) -- (2.4,1.3) -- (2.4,3)
	(2.7,0) -- (2.7,1.3) -- (2.7,3);
\draw[white,fill=white,rotate around={45:(0,0)}]
	(1.5,0.65) circle (10pt)
	(1.5,2.15) circle (10pt);
\draw[rotate around={45:(0,0)}]
	(-0.2,-0.2) node {$\scriptstyle T(0,0)$}
	(-0.6,1.5) node {$\scriptstyle T(0,m_2)$}
	(3.5,-0.5) node {$\scriptstyle T(n_1,0)$}
	(-0.5,3.5) node {$\scriptstyle T(0,n_2)$}
	(3.6, 1.1) node {$\scriptstyle T(n_1,m_2)$}
	(3.2,3.2) node {$\scriptstyle T(n_1,n_2)$}
	(1.5,0.65) node {$A$}
	(1.5,2.15) node {$B$}; 
\end{tikzpicture}%
\label{fig:superadditivitytwo}%
}%
\caption{The two superadditivity situations considered in the
proof of \cref{pr:X}.}
\end{figure}
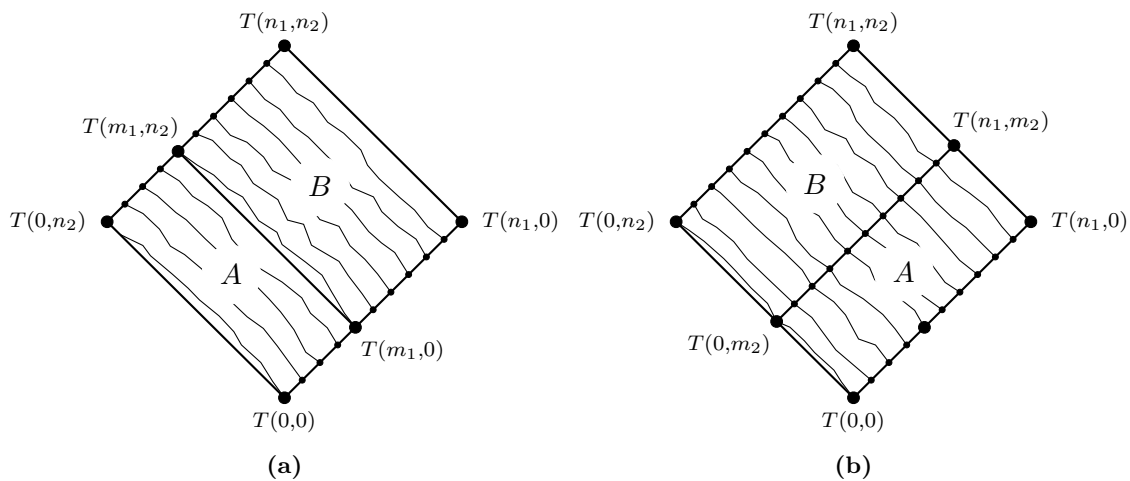
Since the disjoint sets
$A\cup T([0,m_1)\times\{0,n_2\})$ and
$B\cup T([m_1,n_1)\times\{0,n_2\})$ are
$m_1$- and $(n_1-m_1)$-decreasing, respectively,
their union is $n_1$-decreasing. This means that
$A\cup B$ is an $n_1$-decreasing set compatible with
$T([0,n_1)\times\{0,n_2\})$, and it follows that
$X_{\bfzero,\bfn}\ge X_{\bfzero,(m_1,n_2)}+X_{(m_1,0),\bfn}$.

\smallskip
Now, instead let $A$ be an $n_1$-decreasing subset of
$\sigma\cap T[\bfzero,(n_1,m_2))$ compatible with
$T([0,n_1)\times\{0,m_2\})$
and let $B$ be an $n_1$-decreasing subset of
$\sigma\cap T[(0,m_2),\bfn)$ compatible with
$T([0,n_1)\times\{m_2,n_2\})$, as depicted in
\cref{fig:superadditivitytwo}.
The $n_1$-decreasing set
$A\cup T([0,n_1)\times\{0,m_2\})$ is a union of $n_1$ decreasing sets $A_1,A_2,\dotsc,A_{n_1}$. Since no two elements of
$T([0,n_1)\times\{m_2\})$ can belong to the same decreasing set,
we may assume that $T(i-1,m_2)\in A_i$ for $i=1,2,\dotsc,n_1$.
Analogously, $B\cup T([0,n_1)\times\{m_2,n_2\})$ is a union of
$n_1$ decreasing sets $B=B_1\cup B_2\cup\dotsb\cup B_{n_1}$ such that
$T(i-1,m_2)\in B_i$ for $i=1,2,\dotsc,n_1$.
Clearly, $A_i\cup B_i$ is decreasing, and
\[
A\cup B\cup T([0,n_1)\times\{0,m_2,n_2\})
=\bigcup_{i=1}^{n_1} A_i\cup B_i,
\]
so $A\cup B$ is an $n_1$-decreasing subset of $\sigma\cap T[\bfzero,\bfn)$
compatible with $T([0,n_1)\times\{0,m_2,n_2\})$ and hence with
$T([0,n_1)\times\{0,n_2\})$.
It follows that
$X_{\bfzero,\bfn}\ge X_{\bfzero,(n_1,m_2)}+X_{(0,m_2),\bfn}$.

\smallskip
We have showed that the family
$\{X_{\bfm,\bfn}\}_{\bfm<\bfn}$ is translation invariant, superadditive and integrable.
By \cref{th:schurger}, it follows that the limit
$X_\infty:=\lim_{\bfn\rightarrow\infty}X_{\bfzero,\bfn}/n_1n_2$
exists in $L^1$ and equals
\begin{equation}\label{eq:doublelimit}
\lim_{\bfn\rightarrow\infty}(n_1n_2)^{-1}
\lim_{\bfm\rightarrow\infty}(m_1m_2)^{-1}
\sum_{\bfzero<\bfk\le\bfm}X_{\bfk\ast\bfn-\bfn,\bfk\ast\bfn}.
\end{equation}
For any fixed $\bfm,\bfn\in\bbN^2$,
\[
S_{\bfm,\bfn} := (m_1m_2)^{-1}
\sum_{\bfzero<\bfk\le\bfm}X_{\bfk\ast\bfn-\bfn,\bfk\ast\bfn}
\]
is an average of $m_1m_2$ i.i.d.\ random variables of
finite expectation.
By the law of large numbers, as $\bfm\rightarrow\infty$, 
$S_{\bfm,\bfn}$ converges almost surely to
$\Expect X_{\bfzero,\bfn}$. It follows that the limit in
\cref{eq:doublelimit} is concentrated to a constant $c_s$.
\end{proof}

The following proposition defines the function
$\Phi:\bbR_{\ge0}\rightarrow\bbR_{\ge0}$.
\begin{prop}\label{pr:phiexists}
Let $\Omega$ be the open rectangle
\[
0<(x+y)/\sqrt2<\alpha,\ \ 0<(y-x)/\sqrt2<\beta
\]
and let $r\ge0$.
For each $\gamma>0$, let $\sigma_\gamma$ be a Poisson point
process in the plane with homogeneous intensity $\gamma$.
Define the union of lines
\[
L_\gamma:=\bigcup_{i=0}^{\lfloor \alpha r\sqrt{\gamma}\rfloor-1}
\{(x,y)\in\bbR^2\,:\,(x+y)/\sqrt2=i\big/r\sqrt{\gamma}\}
\]
and let $P_\gamma:=\{(x,y)\in L_\gamma\,:\,y-x=0\}$ and
$Q_\gamma:=\{(x,y)\in L_\gamma\,:\,(y-x)/\sqrt2=\beta\}$. 
See \cref{fig:rectanglewithsquares} for an illustration.
\begin{figure}
\begin{tikzpicture}[scale=0.9]
\draw[fill=lightgray, rotate around={-45:(0,0)}] (2,0) rectangle (5,2);
\draw[thick,rotate around={-45:(0,0)}] (0,0) rectangle (7,2);
\draw [decorate,decoration={brace,amplitude=5pt,raise=4pt},rotate around={-45:(0,0)}]
(7,2) -- (7,0) node [black,midway,right,xshift=6pt,yshift=-12pt]
{$\alpha$};
\draw [decorate,decoration={brace,amplitude=5pt,raise=4pt},rotate around={-45:(0,0)}]
(7,0) -- (0,0) node [black,midway,below,xshift=-12pt,yshift=-6pt]
{$\beta$};
\draw [decorate,decoration={brace,amplitude=5pt,raise=4pt},rotate around={-45:(0,0)}]
(0,2) -- (2,2) node [black,midway,above,xshift=12pt,yshift=6pt]
{$\alpha$};
\draw [decorate,decoration={brace,amplitude=5pt,raise=4pt},rotate around={-45:(0,0)}]
(5,2) -- (7,2) node [black,midway,above,xshift=12pt,yshift=6pt]
{$\alpha$};
\draw[fill=black,rotate around={-45:(0,0)}] (7,0) circle (2pt) (7,0.3) circle (2pt) (7,0.6) circle (2pt) (7,0.9) circle (2pt) (7,1.2) circle (2pt) (7,1.5) circle (2pt) (7,1.8) circle (2pt);
\draw[fill=black,rotate around={-45:(0,0)}] (0,0) circle (2pt) (0,0.3) circle (2pt) (0,0.6) circle (2pt) (0,0.9) circle (2pt) (0,1.2) circle (2pt) (0,1.5) circle (2pt) (0,1.8) circle (2pt);
\draw[rotate around={-45:(0,0)}] (0.5,1) node {$Q_\gamma$}; 
\draw[rotate around={-45:(0,0)}] (6.5,1) node {$P_\gamma$}; 
\end{tikzpicture}
\caption{The rectangle $\Omega$ and the sets $P_\gamma$ and $Q_\gamma$
in \cref{pr:phiexists,pr:narrowrectangles}. The shaded area shows the
rectangle $\Omega'$ in the proof of \cref{pr:narrowrectangles}.}
\label{fig:rectanglewithsquares}
\end{figure}
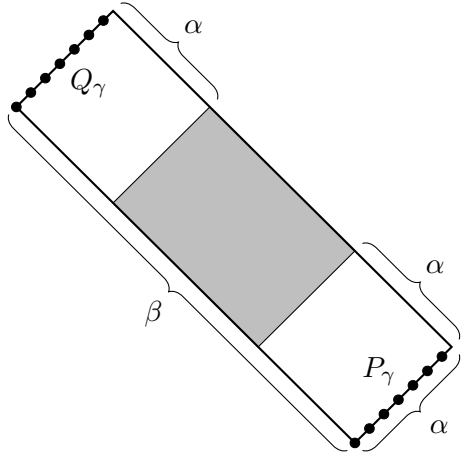
Let $M_\gamma$ be the size
of a maximal $\lfloor \alpha r\sqrt{\gamma}\rfloor$-decreasing subset
of $\sigma_\gamma\cap\Omega$ compatible with $P_\gamma\cup Q_\gamma$.
Then there is a constant $\Phi(r)$, independent of $\alpha$, $\beta$ and $\gamma$, such that
\[
M_\gamma/\alpha\beta\gamma\rightarrow \Phi(r)
\]
in $L^1$ as $\alpha\sqrt{\gamma}$ and $\beta\sqrt{\gamma}$
tend to infinity simultaneously in any manner.
\end{prop}
\begin{proof}
Clearly, $\Phi(0)$ exists and is $0$, so in the following we may assume that $r>0$.
First suppose $n_1:=\alpha r\sqrt{\gamma}$ and $n_2:=\beta r\sqrt{\gamma}$ are both integers. Put $s:=r^{-2}/2$ and define
$X_{\bfzero,(n_1,n_2)}$ as in \cref{pr:X}.
Then $X_{\bfzero,(n_1,n_2)}$ has the same distribution as
$M_\gamma$, so $M_\gamma/\alpha\beta\gamma$ converges in $L^1$ to
$\Phi(r):=r^2 c_s$ as $n_1$ and $n_2$ tends to infinity under the constraint
that they are integers.

Let $\alpha':=\lfloor \alpha r\sqrt{\gamma}\rfloor/r\sqrt{\gamma}$ and
$\beta':=\lfloor \beta r\sqrt{\gamma}\rfloor/r\sqrt{\gamma}$,
and let $\Omega'$, $L'_\gamma$, $P'_\gamma$, $Q'_\gamma$ and
$M'_\gamma$ be defined as $\Omega$,
$L_\gamma$, $P_\gamma$, $Q_\gamma$ and $M_\gamma$ but with $\alpha$ and $\beta$ replaced by $\alpha'$ and $\beta'$.
Then, $L'_\gamma=L_\gamma$ and $P'_\gamma=P_\gamma$.
Any $\lfloor \alpha'r\sqrt{\gamma}\rfloor$-decreasing subset
of $\sigma_\gamma\cap\Omega'$ compatible with $P'_\gamma\cup Q'_\gamma$
is also a $\lfloor \alpha r\sqrt{\gamma}\rfloor$-decreasing subset
of $\sigma_\gamma\cap\Omega$ compatible with $P_\gamma\cup Q_\gamma$,
so $M_\gamma\ge M'_\gamma$ and
$M_\gamma-M'_\gamma\le\#(\sigma_\gamma\cap(\Omega\setminus\Omega'))$.
Since $\alpha-\alpha'$ and $\beta-\beta'$ are both nonnegative and smaller than
$1/r\sqrt{\gamma}$ we have
\begin{multline*}
\mu(\Omega\setminus\Omega')=\alpha\beta-\alpha'\beta'=(\alpha-\alpha')\beta+(\beta-\beta')\alpha-(\alpha-\alpha')(\beta-\beta')\\
\le (\alpha-\alpha')\beta+(\beta-\beta')\alpha \le (\alpha+\beta)/r\sqrt{\gamma},
\end{multline*}
where, as always, $\mu$ denotes the Lebesgue measure on $\bbR^2$.

It follows that $\#(\sigma_\gamma\cap(\Omega\setminus\Omega'))$
has a Poisson distribution with mean smaller than
$(\alpha+\beta)\sqrt{\gamma}/r$, and thus
$\abs{M_\gamma-M'_\gamma}/\alpha\beta\gamma$ converges to
zero almost surely as $\alpha\sqrt{\gamma}$ and $\beta\sqrt{\gamma}$
tend to infinity.
We conclude
that $M_\gamma/\alpha\beta\gamma$ converges to $\Phi(r)$ in $L^1$.
\end{proof}

From the next proposition, \cref{th:narrowrectanglesresult} follows immediately.
\begin{prop}\label{pr:narrowrectangles}
With the same setup as in \cref{pr:phiexists},
define the random variable $\Lambda^{(\gamma)}$ as the size
of a maximal $\lfloor \alpha r\sqrt{\gamma}\rfloor$-decreasing subset
of $\sigma_\gamma\cap\Omega$. Then, as $\alpha\sqrt{\gamma}$
and $\beta/\alpha$ tends to infinity, we have
$\Lambda^{(\gamma)}/\alpha\beta\gamma\rightarrow \Phi(r)$ in $L^1$. Also, for any
fixed $\alpha$ and $\beta$, the inequality
\[
\abs{(\Lambda^{(\gamma)}/\alpha\beta\gamma)-\Phi(r)}<3\alpha/\beta
\]
holds \aas\ as $\gamma\rightarrow\infty$.
\end{prop}
\begin{proof}
Clearly, $\Lambda^{(\gamma)}\ge M_\gamma$.
Let $\Omega'$ be the (possibly empty) rectangle given by
$0<(x+y)/\sqrt2<\alpha$ and $\alpha<(y-x)/\sqrt2<\beta-\alpha$; see \cref{fig:rectanglewithsquares} for an illustration.
If $A$ is a
$\lfloor \alpha r\sqrt{\gamma}\rfloor$-decreasing subset
of $\sigma_\gamma\cap\Omega$, then
$A\cap\Omega'$ is a $\lfloor \alpha r\sqrt{\gamma}\rfloor$-decreasing subset
of $\sigma_\gamma\cap\Omega$ compatible with
$P_\gamma\cup Q_\gamma$. Thus,
$\Lambda^{(\gamma)}-M_\gamma\le\#(\sigma_\gamma\cap(\Omega\setminus\Omega'))$.
We have $\mu(\Omega\setminus\Omega')\le2\alpha^2$, so
$\#(\sigma_\gamma\cap(\Omega\setminus\Omega'))$ has a Poisson
distribution with mean at most $2\alpha^2\gamma$, and hence
$(\Lambda^{(\gamma)}-M_\gamma)/\alpha\beta\gamma$ converges almost surely
to zero as $\alpha\sqrt{\gamma}$ and $\beta/\alpha$ tends to infinity.
It follows from \cref{pr:phiexists} that
$\Lambda^{(\gamma)}/\alpha\beta\gamma\rightarrow \Phi(r)$ in $L^1$.
For fixed $\alpha$ and $\beta$, and for any $\delta>0$, it holds that
$(\Lambda^{(\gamma)}-M_\gamma)/\alpha\beta\gamma<2(1+\delta)\alpha/\beta$
\aas\ as $\gamma\rightarrow\infty$, and hence
$\abs{(\Lambda^{(\gamma)}/\alpha\beta\gamma)-\Phi(r)}<3\alpha/\beta$ \aas
\end{proof}

In \cref{sec:localparallelogram} we will need the following more flexible
version of \cref{pr:narrowrectangles} that incorporates the functional
$\functional_\rho$ from \cref{df:FrhoandL}.
\begin{prop}\label{pr:beta}
Let $\Omega$ be the open parallelogram
\[
\{(x,y)\in\bbR^2\ :\ \abs{ax+by}<1,\ \abs{ax-by}<\beta\}
\]
for some $a,b,\beta>0$,
and let $\rho$ be constant on $\Omega$. Let $u_{\rm linear}(x,y)=c(ax+by)$ for some $c\ge0$, and let $\{\sigma_\gamma\}_{\gamma>0}$ be Poisson point
processes on $\Omega$ with intensities $\gamma\rho$. Let $\Lambda^{(\gamma)}$ be
the size of a maximal $\lfloor 2c\sqrt{\gamma}\rfloor$-decreasing subset of $\sigma_\gamma$. Then
\[
\abs{(\Lambda^{(\gamma)}/\gamma) - \functional_\rho(u_{\rm linear})} \le 3\rho\mu(\Omega)/\beta
\]
\aas\ as $\gamma\rightarrow\infty$.
\end{prop}
\begin{proof}
If $\rho=0$, we have $\Lambda^{(\gamma)}=0$ and $\functional_\rho(u_{\rm linear})=0$
so the conclusion of the proposition is true.

In the following we assume that $\rho>0$.
Rescale the $x$- and $y$-axes and generate Poisson point processes with
intensities $\gamma'=2\gamma\rho/ab$ on the rectangle
$\abs{x+y}<1/\sqrt2$, $\abs{x-y}<\beta/\sqrt2$. In
\cref{pr:narrowrectangles}
this corresponds to $r:=c\sqrt{2ab/\rho}$, $\alpha:=1$, $\beta:=\beta$,
$\gamma:=\gamma'$, and the proposition yields that
$\abs{(\Lambda^{(\gamma)}/\beta\gamma')-\Phi(r)} < 3/\beta$ \aas\ as $\gamma\rightarrow\infty$.
It is straightforward to check that
$\mu(\Omega)=2\beta/ab$ and
$\functional_\rho(u_{\rm linear})=2\beta\rho \Phi(r)/ab$.
\end{proof}

Our next goal is to show some nice properties of $\Phi$, in particular
that it is increasing and concave. To this end, we need a couple
of lemmas which will be used again later on when we concern ourselves
with limit shapes of Young diagrams.

Let $\partial_-$ and $\partial_+$ denote the left and right one-sided derivative operators.
\begin{lemma}\label{lm:derivativelimit}
Let $F_1,F_2,\dotsc$ be random concave functions from
$\bbR_{\ge0}$ to $\bbR$, and suppose there is a deterministic function
$F:\bbR_{\ge0}\rightarrow\bbR$ such that
$F_n(x)\rightarrow F(x)$ in probability for any $x$. Then the following holds.
\begin{thmlist}
\item\label{lm:derivativelimitFisconcave}
$F$ is concave.
\item\label{lm:derivativelimitrightderivative}
$\partial_-F_n(x)\rightarrow F'(x)$ and $\partial_+F_n(x)\rightarrow F'(x)$ in probability for any point
$x>0$ where $F'(x)$ exists.
\end{thmlist}
\end{lemma}
\begin{proof}
\reflocal{lm:derivativelimitFisconcave}
Take any $0\le x<y$ and $0<t<1$. We must show that
$(1-t)F(x)+tF(y)\le F((1-t)x+ty)$.
For any $\eps>0$, \aas\ as $n$ tends to infinity we have
\begin{align}
F(x) &< F_n(x)+\eps, \nonumber \\
F(y) &< F_n(y)+\eps\ \text{and} \nonumber \\
F((1-t)x+ty) &> F_n((1-t)x+ty)-\eps. \label{eq:FgtFn}
\end{align}
The first two inequalities imply that
\[
(1-t)F(x)+tF(y)<(1-t)F_n(x)+tF_n(y)+\eps
\]
which is less than or equal to $F_n((1-t)x+ty)+\eps$ since $F_n$
is concave.
Combining this with \cref{eq:FgtFn}, we obtain
\[
(1-t)F(x)+tF(y)<F((1-t)x+ty)+2\eps.
\]
Since this holds with positive probability and $F$ is deterministic,
it holds deterministically, and since it holds for any $\eps>0$,
we must have $(1-t)F(x)+tF(y)\le F((1-t)x+ty)$.

\reflocal{lm:derivativelimitrightderivative}
Suppose $x>0$ is a point where $F'(x)$ exists, and take any $\eps>0$.
Since $F'(x)$ exists, there is a $\delta>0$ such that
\begin{equation}\label{eq:derivativeapprox}
\begin{aligned}
\frac{F(x)-F(x-\delta)}{\delta}-F'(x) &< \eps, \\
\frac{F(x+\delta)-F(x)}{\delta}-F'(x) &> -\eps.
\end{aligned}
\end{equation}
Since $F_n(x)\rightarrow F(x)$, $F_n(x-\delta)\rightarrow F(x-\delta)$ and
$F_n(x+\delta)\rightarrow F(x+\delta)$ in probability, the inequalities
\begin{equation}\label{eq:functionapprox}
\begin{aligned}
\frac{F_n(x)-F_n(x-\delta)}{\delta}
&\le\frac{F(x)-F(x-\delta)}{\delta}+\eps, \\
\frac{F_n(x+\delta)-F_n(x)}{\delta}
&\ge\frac{F(x+\delta)-F(x)}{\delta}-\eps
\end{aligned}
\end{equation}
hold \aas\ as $n\rightarrow\infty$. 
Since $F_n$ is concave, it has a left derivative and a right derivative
at $x$, and
\begin{equation}\label{eq:derivatesqueeze}
\frac{F_n(x+\delta)-F_n(x)}{\delta}
\le\partial_+F_n(x)
\le\partial_-F_n(x)
\le\frac{F_n(x)-F_n(x-\delta)}{\delta}.
\end{equation}
Combining \cref{eq:derivativeapprox,eq:functionapprox,eq:derivatesqueeze} yields
\[
F'(x)-2\eps
<\partial_+F_n(x)
\le\partial_-F_n(x)
<F'(x)+2\eps,
\]
and, since $\eps>0$ was chosen arbitrarily, we conclude that
$\partial_-F_n(x)\rightarrow F'(x)$ and $\partial_+F_n(x)\rightarrow F'(x)$ in probability.
\end{proof}

\begin{lemma}\label{lm:Lambdatolambda}
Let $\{\lambda^{(\gamma)}\}_{\gamma>0}$ be a family of
random Young diagrams and let
$\Lambda^{(\gamma)}_k:=\sum_{i=1}^k\lambda^{(\gamma)}_i$.
Let $a$ and $b$ be positive functions of $\gamma$ such that
$\lim_{\gamma\rightarrow\infty}a(\gamma)=\infty$.
Suppose there is a (deterministic) function
$G\,:\,\bbR_{\ge0}\rightarrow\bbR_{\ge0}$ such that
\[
b(\gamma)\Lambda^{(\gamma)}_{\lfloor a(\gamma)r\rfloor}\rightarrow G(r)
\]
in probability for any $r\ge0$. Then $G$ is increasing and concave, and
\[
a(\gamma)b(\gamma)\lambda^{(\gamma)}_{\lfloor a(\gamma)r\rfloor+1}\rightarrow G'(r)
\]
in probability for any $r>0$ where $G$ is differentiable.
\end{lemma}
\begin{proof}
For each $\gamma>0$, define the function
$F^{(\gamma)}\,:\,\bbR_{>0}\rightarrow\bbR_{\ge0}$ by
\[
F^{(\gamma)}(r):=\int_0^{a(\gamma)r} \lambda^{(\gamma)}_{\lfloor t\rfloor+1}\,dt.
\]
Since the integrand is a nonnegative decreasing function, $F^{(\gamma)}$ is increasing and concave, and since the integrand is piecewise constant,
\[
\abs{F^{(\gamma)}(r)-\Lambda^{(\gamma)}_{\lfloor a(\gamma)r\rfloor}}
=
\bigl(a(\gamma)r-\lfloor a(\gamma)r\rfloor\bigr)
\lambda^{(\gamma)}_{\lfloor a(\gamma)r\rfloor+1}
\le
\frac{a(\gamma)r-\lfloor a(\gamma)r\rfloor}{\lfloor a(\gamma)r\rfloor}
\sum_{i=1}^{\lfloor a(\gamma)r\rfloor}\!\!\!\lambda^{(\gamma)}_i
\]
which equals $\Lambda^{(\gamma)}_{\lfloor a(\gamma)r\rfloor}
$ times a factor that tends to zero as $\gamma\rightarrow\infty$.
It follows that
\[
b(\gamma)F^{(\gamma)}(r)\rightarrow G(r)
\]
in probability. Since all $F^{(\gamma)}$ are increasing, $G$ is increasing too. Furthermore, by \cref{lm:derivativelimit}, $G$ is concave and
$b(\gamma)\partial_+ F^{(\gamma)}(r)\rightarrow G'(r)$
for any $r$ where $G$ is differentiable. Since 
$\partial_+ F^{(\gamma)}(r)=a(\gamma)\lambda^{(\gamma)}_{\lfloor a(\gamma)r\rfloor+1}$ the lemma follows.
\end{proof}

\begin{prop}\label{pr:phiproperties}
$\Phi$ is increasing, concave, continuous and bounded by one.
Furthermore, $\Phi(0)=0$.
\end{prop}
\begin{proof}
Consider the setup of \cref{pr:narrowrectangles} with the specialization $\alpha=1$ and $\beta=\gamma$.
Then $\Lambda^{(\gamma)}/\gamma^2\rightarrow \Phi(r)$ in $L^1$
as $\gamma\rightarrow\infty$.
Let $\lambda^{(\gamma)}$ be
the random Young diagram corresponding to $\sigma_\gamma\cap\Omega$.
By \cref{pr:curtisgreene},
$\Lambda^{(\gamma)}=\sum_{i=1}^{\lfloor r\sqrt{\gamma}\rfloor}\lambda^{(\gamma)}_i$, and
\cref{lm:Lambdatolambda} with $a(\gamma)=\sqrt{\gamma}$,
$b(\gamma)=1/\gamma^2$ and $G=\Phi$ yields that $\Phi$ is increasing and
concave.

That $\Phi(0)=0$ and that $\Phi$ is bounded by one follows directly from its definition together with the law of large numbers.

Since $\Phi$ is concave,
it is automatically continuous on the open set $(0,\infty)$. It remains only to show that it is continuous at 0.

For any $\beta>0$, let $\Omega_\beta$ be the open
rectangle
\[
0<(x+y)/\sqrt2<1,\ \ 0<(y-x)/\sqrt2<\beta,
\]
and for any $\gamma>0$ and $\beta>0$, let $\sigma_{\gamma,\beta}$ be a Poisson point
process on $\Omega_\beta$ with homogeneous intensity $\gamma$.
Since $\Omega_\beta\subset
(-\tfrac\beta{\sqrt2},\tfrac1{\sqrt2})\times(0,\tfrac{1+\beta}{\sqrt2})$, by \cref{th:hammersley},
for any $\eps>0$, the size of the largest decreasing subset
of $\sigma_{\gamma,\beta}$ is smaller than $\frac\Gamma{\sqrt2}(1+\eps)(1+\beta)\sqrt{\gamma}$
\aas\ as $\gamma\rightarrow\infty$. (We know from the result of Vershik and Kerov \cite{VershikKerov} that $\Gamma=2$ but we will not need that now.)
Thus, for any $r>0$,
the maximum
size $\Lambda^{(\gamma)}$ of a $\lfloor r\sqrt{\gamma}\rfloor$-decreasing subset
of $\sigma_{\gamma,\beta}$ is smaller than
$\tfrac\Gamma{\sqrt2}r(1+\eps)(1+\beta)\gamma$ \aas\ as $\gamma\rightarrow\infty$.
By \cref{pr:narrowrectangles}, $\abs{\Lambda^{(\gamma)}/\beta\gamma-\Phi(r)}<3/\beta$ and thus
$\Lambda^{(\gamma)}/\gamma > \beta\Phi(r) - 3$ \aas\ as
$\gamma\rightarrow\infty$. So
$\tfrac\Gamma{\sqrt2}r(1+\eps)(1+\beta)>\Lambda^{(\gamma)}/\gamma > \beta\Phi(r) - 3$ \aas, and hence
$\Phi(r)<\bigl(\tfrac\Gamma{\sqrt2}r(1+\eps)(1+\beta)+3\bigr)/\beta$ for any $r,\eps,\beta>0$.
Letting $\beta\rightarrow\infty$ and $\eps\rightarrow0$ yields
$\Phi(r)\le \tfrac\Gamma{\sqrt2}r$, and it follows that $\Phi$ is continuous at 0.
\end{proof}

As a consequence of the nice properties of $\Phi$, the function $L$
from \cref{df:FrhoandL} is well behaved too.
\begin{lemma}\label{lm:Lcontinuous}
$L$ is continuous and increasing in both variables.
Furthermore, for any $\eta,\theta\ge0$ it holds that $L(\eta,0)=0$ and $L(\eta,\theta)\le\eta$.
\end{lemma}
\begin{proof}
By \cref{pr:phiproperties}, $L$ is continuous at all points
$(\eta,\theta)$ with $\eta>0$. At any point $(0,\theta)$, continuity of $L$ follows from the fact that $\Phi$ is bounded.

That $L$ is increasing in the second variable follows from the fact that
$\Phi$ is increasing (\cref{pr:phiproperties}). That $L(\eta,0)=0$ and $L(\eta,\theta)\le\eta$ for any $\eta,\theta\ge0$ follows from
the facts that $\Phi(0)=0$ and that $\Phi$ is bounded by one (\cref{pr:phiproperties}).

To show that $L$ is increasing in the first variable, take any $\theta\ge0$
and any $\eta'>\eta>0$.
Since $\Phi$ is concave (\cref{pr:phiproperties}),
\[
\Phi(\sqrt{2\theta/\eta'}) \ge \left(1-\sqrt{\eta/\eta'}\right)\Phi(0)
+\sqrt{\eta/\eta'}\,\Phi(\sqrt{2\theta/\eta})
=\sqrt{\eta/\eta'}\,\Phi(\sqrt{2\theta/\eta}),
\]
and it follows that
$
L(\eta',\theta)=\eta'\Phi(\sqrt{2\theta/\eta'})
\ge \sqrt{\eta'/\eta}\,L(\eta,\theta)\ge L(\eta,\theta)$.
\end{proof}

Finally, \cref{lm:Lcontinuous}, together with the following lemma,
shows that
the functional $\functional_\rho$ from \cref{df:FrhoandL} is well defined\label{sec:Fwelldefined}.
\begin{lemma}\label{lm:monotonic}
Let $u$ be a function from $\bbR^2$ to $\bbR$ that is increasing in both variables. Then the following holds.
\begin{thmlist}
\item\label{lm:monotonicmeasurable}
$u$ is measurable.
\item\label{lm:monotonicdifferentiable}
$u$ is differentiable almost everywhere.
\item\label{lm:monotonicpartialderivatives}
The partial derivatives of $u$ exist almost everywhere and are measurable.
\end{thmlist}
\end{lemma}
\begin{proof}
\reflocal{lm:monotonicmeasurable}
Let $a$ be a real number. We must show that $u^{-1}((-\infty,a])$
is measurable.
Define a function $g$ by letting $g(x):=\sup\{y\,:\,u(x,y) \le a\}$
whenever the supremum exists.
Then the domain of $g$ is an interval and $g$ is decreasing, so
it is measurable and its graph has measure zero. The inverse
image $u^{-1}((-\infty,a])$ is the
region below this graph (a measurable set) plus some subset of
the graph itself (a null set).

\reflocal{lm:monotonicdifferentiable}
This is proved in \cite[Sec.~6]{BurkillHaslamjones}.

\reflocal{lm:monotonicpartialderivatives}
This follows from \reflocal{lm:monotonicmeasurable}, \reflocal{lm:monotonicdifferentiable} and \cite[Lemma~2]{BurkillHaslamjones}.
\end{proof}

The set of points where $u$ is differentiable is denoted by
$\Diff(u)$.

\section{Maximizing \texorpdfstring{$\functional_\rho$}{the functional} by solving a PDE system}
\label{sec:generalpde}

In this section we show that the problem of maximizing the functional
$\functional_\rho$ can be reduced to a system of partial differential equations.
This will let us compute, in terms of $\Phi$, the maximum for a parallelogram density domain with constant density. Furthermore, if \cref{con:triangularlimitshape} holds, the PDE system simplifies significantly
and can be solved analytically for the uniform case as we will see in \cref{sec:uniform}.

Let us recall some standard facts from convex analysis.

Given a continuous convex function $f$ from an interval $I\subseteq\bbR$ to $\bbR$,
its \emph{Legendre transform} $\calL[f]$ is a function defined by
\[
\calL[f](s) = \sup_{r\in I} \bigl(rs-f(r)\bigr)
\]
for those $s$ for which the supremum exists.
It is well known that $\calL(f)$ is a continuous convex function and that its
domain is an interval. Furthermore, the \emph{Fenchel--Young inequality}
states that
\[
f(r)+\calL[f](s)\ge rs
\]
for any $r\in I$ and $s\in\dom \calL[f]$, with equality if and only if
$s\in\partial f(r)$, where the \emph{subdifferential set} $\partial f(r)$
is defined by
\[
\partial f(r)=\{s\ :\ f(z)-f(r)\ge (z-r)s\ \text{for any}\ z\in I\}.
\]
It is also known that $\partial f(r)$ is nonempty for any $r$.

\begin{defi}
Let $\Phi^\ast$ be the real function defined on $\bbR_{\ge0}$ by
\[
\Phi^\ast(s)=\inf_{r\ge0}\bigl(rs - \Phi(r)\bigr).
\]
\end{defi}
Note that the infimum exists since $\Phi$ is bounded by \cref{pr:phiproperties}.

\begin{lemma}\label{lm:phibound}
$\Phi^\ast$ is continuous and $-1\le\Phi^\ast(s)\le0$ for any $s\ge0$. Furthermore,
\begin{thmlist}
\item\label{lm:phiboundA}
$\Phi(r) + \Phi^\ast(s) \le rs$ for any $r, s\ge 0$, and
\item\label{lm:phiboundB}
for any $r$ there is an $s$ such that $\Phi(r) + \Phi^\ast(s) = rs$.
\end{thmlist}
\end{lemma}
\begin{proof}
Once we note that $\Phi^\ast(s) = -\calL[-\Phi](-s)$, it follows that
$\Phi^\ast$ is continuous, and part \reflocal{lm:phiboundA} follows from
the Fenchel--Young inequality while part \reflocal{lm:phiboundB} follows from the fact that all subdifferential sets are nonempty. To see that
$-1\le\Phi^\ast(s)\le0$ we observe that
\[
-1\le-\sup_{r\ge0}\Phi(r)\le\inf_{r\ge0}\bigl(rs-\Phi(r)\bigr)\le 0\cdot s-\Phi(0)
=0,
\]
where the first and last inequalities follow from \cref{pr:phiproperties}.
\end{proof}

Let us expand our terminology for doubly increasing functions
to include functions decreasing in $x$ and increasing in $y$.
Define a partial order $\le'$
on $\bbR^2$ by letting $(x_1,y_1)\le'(x_2,y_2)$ if
$x_1\ge x_2$ and $y_1\le y_2$.
For any subset $A$ of $\bbR^2$,
a function $v:A\rightarrow\bbR$ is \emph{decreasing in $x$ and increasing in $y$} if
$u(x_1,y_1)\le u(x_2,y_2)$ whenever $(x_1,y_1)\le'(x_2,y_2)$.
For $s\ge0$, we let $\calV_s(A)$ denote the set of functions
$v$ on $A$ that are decreasing in $x$ and increasing in $y$ and have
$\diam v(A)\le s$, and we let $\calV(A):=\bigcup_{s\ge0}\calV_s(A)$
denote the set of all bounded functions on $A$ decreasing in $x$ and increasing in $y$.
Let $\calV_{h,s}(A)$ denote the subset of $\calV_s(A)$ consisting of functions with values in $[h,h+s]$.

For any density domain $(\Omega,\rho)$, let
$\functional^\ast_\rho:\calV(\Omega)\rightarrow\bbR$
be a (nonlinear) functional given by
\[
\functional^\ast_\rho(v) = \int_{\Omega}\rho\,
\Phi^\ast\bigl(\lagomsqrt{-2v_xv_y/\rho}\bigr)\,d\mu,
\]
where the integrand is defined to be zero at points where $\rho=0$.
This functional is well defined by \cref{lm:monotonic} together with the fact
that $\Phi^\ast$ is continuous and bounded by \cref{lm:phibound}.

\begin{lemma}\label{lm:injectivity}
Let $\Omega$ be an open subset of $\bbR^2$,
let $u\in\calU(\Omega)$ and $v\in\calV(\Omega)$, and
let $A$ be the subset of $\Omega$ where
$u_xv_y-u_yv_x>0$. Then, the map $\varphi\,:\,A\rightarrow\bbR^2$
defined by $\varphi(x,y)=(u(x,y),v(x,y))$ is injective.
\end{lemma}
\begin{proof}
Suppose there are two
distinct points $p=(x,y)$ and $q=(x',y')$ in $A$
with $\varphi(p)=\varphi(q)$.
Without loss of generality we may assume that $x\le x'$.
Since $\Omega$ is open, for any sufficiently small
$\eps>0$ the point $p_\eps := (1-\eps)p + \eps q$ belongs to $\Omega$.
If $x=x'$ we have $u(p_\eps)=u(p)$ and $v(p_\eps)=v(p)$
and thus $u_y(p)=v_y(p)=0$.
If $y=y'$, we have $u(p_\eps)=u(p)$ and $v(p_\eps)=v(p)$
and thus $u_x(p)=v_x(p)=0$.
If $x<x'$ and $y<y'$ we have $u(p_\eps)=u(p)$ and thus
$u_x(p)=u_y(p)=0$.
If $x<x'$ and $y>y'$
we have $v(p_\eps)=v(p)$ and thus $v_x(p)=v_y(p)=0$.
In any of the four cases above, we conclude that
$u_xv_y-u_yv_x=0$ in $p$, and it follows that $p\not\in A$, a contradiction.
\end{proof}

We will need the following ``change of variables''
theorem that appears as Theorem~263D in~\cite{Fremlin}.
\begin{theo}\label{th:changeofvariables}
Let $D\subseteq\bbR^n$ be any measurable set,
and $\varphi:D\rightarrow\bbR^n$
a function differentiable relative to its domain\footnote{We say that $\varphi$ is \emph{differentiable relative to its domain} at a point $p\in D$ if there is a linear map $T(p)\,:\,\bbR^n\rightarrow\bbR^n$ (called a derivative of $\varphi$ relative to $D$ in $p$) such that for each $\eps>0$ there is a $\delta>0$ such that $\abs{\varphi(p)+T(p)(x-p)-\varphi(x)}\le\eps\abs{x-p}$ for any $x\in D$ with $\abs{x-p}<\delta$.} at each point of $D$.
For each $p\in D$, let $T(p)$ be a derivative of $\varphi$ relative to $D$
at $p$, and set $J(p):=\abs{\det T(p)}$. Then
$\mu(\varphi(D))\le \int_D J\,d\mu$
with equality if $\varphi$ is injective.
\end{theo}

Now we are ready for the main result of this section.
Recall that $\Diff(u)$ denotes the set of points where $u$ is differentiable.
\begin{theo}\label{th:pde}
Let $(\Omega,\rho)$ be a density domain.

Suppose, for some $r,s>0$, there are
$u\in\calU_r(\Omega)$ and $v\in\calV_s(\Omega)$ with the
following properties.
\begin{enumerate}
\item[(a)]\label{it:measurers}
The set $\{(u(x,y),v(x,y))\,:\,(x,y)\in\Diff(u)\cap\Diff(v)\}$ has measure
$rs$.
\item[(b)]\label{it:phiplusphistar}
The PDE system
\begin{gather}
\label{eq:first}
u_xv_y+u_yv_x = 0, \\
\label{eq:second}
\rho\left(\Phi(\lagomsqrt{2u_xu_y/\rho})+\Phi^\ast(\lagomsqrt{-2v_xv_y/\rho})\right)
= 2\sqrt{-u_xu_yv_xv_y}
\end{gather}
is satisfied almost everywhere in $\Omega$,
where the left-hand side of \cref{eq:second} is defined to be zero
at points where $\rho=0$.
\end{enumerate}
Then, $u$ is a maximizer of
$\functional_\rho$ in $\calU_r(\Omega)$ and $v$ is a maximizer of $\functional^\ast_\rho$
in $\calV_s(\Omega)$. Furthermore, $s = \Fmax'(r)$ if $\Fmax'(r)$ exists. 
\end{theo}

\begin{proof}
Let $u$ and $v$ be functions in $\calU_r(\Omega)$ and $\calV_s(\Omega)$,
respectively.
By \cref{lm:phibound},
\begin{align*}
\functional_\rho(u)+\functional^\ast_\rho(v)
&=\int_{\Omega} \rho\,\left(\Phi\bigl(\lagomsqrt{2u_xu_y/\rho}\bigr)
+\Phi^\ast\bigl(\lagomsqrt{-2v_xv_y/\rho}\bigr)\right)\,d\mu \\
&\le 2\int_{\Omega} \sqrt{-u_xu_yv_xv_y}\,d\mu
\end{align*}
with equality if and only if \cref{eq:second} holds
almost everywhere.

By the inequality of the geometric and arithmetic mean,
\begin{equation}\label{eq:geom_arit_mean}
2\int_{\Omega} \sqrt{-u_xu_yv_xv_y}\,d\mu
\le\int_{\Omega} (u_xv_y-u_yv_x)\,d\mu
\end{equation}
with equality if and only if \cref{eq:first} holds
almost everywhere.

Let $D:=\Diff(u)\cap\Diff(v)$.
Let $\varphi$ be the map from $D$ to $\bbR^2$
defined by $\varphi(x,y):=\bigl(u(x,y),v(x,y)\bigr)$, and
let $A$ be the subset of $D$ where $u_xv_y-u_yv_x$ is positive. 
It follows from \cref{lm:injectivity} that $\varphi$ is injective on $A$.
The right-hand side of~\cref{eq:geom_arit_mean} equals
\[
\int_A (u_xv_y-u_yv_x)\,d\mu,
\]
and since $\varphi$ is injective on $A$, by \cref{th:changeofvariables}
this equals $\mu(\varphi(A))$.
By the same theorem,
$\int_{\Omega} (u_xv_y-u_yv_x)\,d\mu \ge \mu(\varphi(D))$, so
$\mu(\varphi(A))\ge\mu(\varphi(D))$,
which implies that
$\mu(\varphi(A))=\mu(\varphi(D))$.
Thus, we obtain
\[
\mu(\varphi(A))
=\mu(\varphi(D))
\le rs
\]
with equality if and only if $u$ and $v$ have property \ref{it:measurers}.
(Note that $\mu(\varphi(D))\le rs$ holds trivially because $u\in\calU_r(\Omega)$
and $v\in\calV_r(\Omega)$.)

In conclusion, we have
\[
\functional_\rho(u)+\functional^\ast_\rho(v)\le rs
\]
with equality if and only if $u$ and $v$ have properties \ref{it:measurers} and \ref{it:phiplusphistar}. It follows
that such $u$ and $v$ are maximizers of $\functional_\rho$ and $\functional^\ast_\rho$
in $\calU_r(\Omega)$ and $\calV_s(\Omega)$, respectively.

It remains to show that $s = \Fmax'(r)$ if $\Fmax'(r)$ exists.
From above, it follows that
\begin{equation}\label{eq:FmaxplusFmaxast}
\Fmax(r) + \functional_\rho^\ast(v)\le rs
\end{equation}
for any $r,s>0$ and any $v\in\calV_s(\Omega)$, and that equality holds if there is a $u\in\calU_r(\Omega)$ such that
\ref{it:measurers} and \ref{it:phiplusphistar} hold.
If equality holds for
some particular $r,s,v$ and if $\Fmax'(r)$ exists, then the partial
derivative with respect to $r$ (while keeping $s$ and $v$ fixed) of
the left- and right-hand sides of \cref{eq:FmaxplusFmaxast} must coincide, so $\Fmax'(r)=s$.
\end{proof}

If \cref{con:triangularlimitshape} holds, the PDE system of
\cref{th:pde} can be written explicitly. We will exploit this fact
in \cref{sec:uniform} where we solve the system for the uniform case.
\begin{prop}\label{pr:pdeunderconjecture}
Suppose \cref{con:triangularlimitshape} holds. Then \cref{th:pde} holds if $\rho>0$ on $\Omega$ and \cref{eq:second} is replaced by
\begin{equation}\label{eq:seconduniform}
\min\{\lagomsqrt{u_xu_y/\rho},1\}+\min\{\lagomsqrt{-v_xv_y/\rho},1\}=1.
\end{equation}
\end{prop}
\begin{proof}
Suppose \cref{con:triangularlimitshape} holds. Then
\[
\Phi(r)=\begin{cases}
\sqrt2\,r-\frac{r^2}2 & \text{if $0\le r\le\sqrt2$,} \\
1 & \text{if $r>\sqrt2$,}
\end{cases}
\]
and $\Phi^\ast(s)=\Phi(s)-1$ for any $s\ge0$.

Let $p=\sqrt{u_xu_y/\rho}$ and $q=\sqrt{-v_xv_y/\rho}$.
If $p>1$, \cref{eq:seconduniform} implies that $q=0$ and hence
$\Phi(\sqrt2\,p)+\Phi^\ast(\sqrt2\,q)=1-1=0$ so
\cref{eq:second} is satisfied.
Analogously, if $q>1$, \cref{eq:seconduniform} implies that $p=0$ and hence
$\Phi(\sqrt2\,p)+\Phi^\ast(\sqrt2\,q)=0+0=0$ so
\cref{eq:second} is satisfied.
If $p,q\le1$, \cref{eq:seconduniform} can be written as
$p+q=1$ which implies $2(p+q)-(p+q)^2=1$ and hence
$2p-p^2+2q-q^2-1=2pq$.
The last equation can be written as
$\Phi(\sqrt2\,p)+\Phi^\ast(\sqrt2\,q)=2pq$ which is equivalent to
\cref{eq:second}.
\end{proof}

\Cref{th:pde} also lets us find a maximizer of $\functional_\rho$
on a parallelogram with constant density.
\begin{prop}\label{pr:parallelogrammaximizer}
Let $\Omega$ be the open parallelogram
\[
\{(x,y)\in\bbR^2\ :\ \abs{ax+by}<1,\ \abs{ax-by}<\beta\}
\]
for some $a,b,\beta>0$, and let $\rho$ be constant on $\Omega$.
Then, for any $c\ge0$, in $\calU_{2c}(\Omega)$
the functional $\functional_\rho$ is maximized by the function
$u(x,y)=c(ax+by)$, and the maximum value is $2\rho\beta\Phi(c\sqrt{2ab/\rho})/ab$ if $\rho>0$ and 0 if $\rho=0$.
\end{prop}
\begin{proof}
If $\rho=0$, the functional $\functional_\rho$ is identically zero,
so we may assume that $\rho>0$.
Note that $\sqrt{2u_xu_y/\rho}=c\sqrt{2ab/\rho}$ which is independent of
$x$ and $y$.
By \cref{lm:phiboundB}, there is a $d\ge0$ such that
$\Phi(c\sqrt{2ab/\rho})+\Phi^\ast(d\sqrt{2ab/\rho})=2abcd/\rho$.
Let $v(x,y):=d(by-ax)$.
Then, $u$ and $v$ satisfy the conditions in \cref{th:pde} with
$r=2c$ and $s=2\beta d$,
and hence $u$ is a maximizer of $\functional_\rho$ in $\calU_{2c}$.
The maximum value is $\mu(\Omega)\rho\Phi(\sqrt{2u_xu_y/\rho})=(2\beta/ab)\rho\Phi(c\sqrt{2ab/\rho})$.
\end{proof}

\section{Probabilistic behavior of the local parallelogram}
\label{sec:localparallelogram}
In \cref{def:kappa}, we defined the map $P\mapsto\kappa_P$
to convert $k$-decreasing sets
into doubly increasing functions. The next definition is a device
to go in the other direction.
\begin{defi}
For any $u$ in $\calU(\Omega)$, we let
\begin{multline*}
\decreasing(u):=\{(x,y)\in\Omega\ :\ u(x,y)\in\mathbb{Z}\ \text{and}\ u(x',y')<u(x,y)\\
\text{for any}\ (x',y')\in\Omega\setminus\{(x,y)\}
\ \text{such that}\ x'\le x\ \text{and}\ y'\le y\}.
\end{multline*}
\end{defi}
\begin{lemma}\label{lm:Iandkappa}
The following holds.
\begin{thmlist}
\item\label{lm:IandkappaIudecreasing}
For any $u$ in $\calU_r(\Omega)$,
$\decreasing(u)$ is $(\lfloor r\rfloor+1)$-decreasing, and if
$u\in\calU_{0,r}(\Omega)$,
$\decreasing(u)$ is $\lfloor r\rfloor$-decreasing.
\item\label{lm:IandkappaIkappaPisP}
If $P$ is a finite set of points in $\Omega\subset\bbR^2$
with distinct $x$-coordinates and distinct $y$-coordinates,
then $\decreasing(\kappa_P)=P$.
\end{thmlist}
\end{lemma}
\begin{proof}
\reflocal{lm:IandkappaIudecreasing}
Let $u\in\calU_r(\Omega)$.
Suppose $(x,y)$ and $(x',y')$ are distinct points in $\decreasing(u)$
such that $u(x,y)=u(x',y')$. By definition of $\decreasing(u)$, it follows that
$\{(x,y),(x',y')\}$ is a decreasing set. Hence, for each integer $k$, the fiber $u^{-1}(k)\cap \decreasing(u)$ is a decreasing set.
The image of $u$ contains at most $\lfloor r\rfloor+1$ integers,
so $\decreasing(u)$ is a union of $\lfloor r\rfloor+1$ decreasing
sets. If $u\in\calU_{0,r}(\Omega)$, the image of $u$ contains no
integer outside the set $\{0,1,\dotsc,\lfloor r\rfloor\}$. Since $\Omega$ is open and
$u\ge0$, for any point $(x,y)\in\Omega$ with $u(x,y)=0$ there is
an $x'<x$ such that $u(x',y)=0$, so the fiber
$u^{-1}(0)\cap \decreasing(u)$ is empty.

\smallskip
\reflocal{lm:IandkappaIkappaPisP}
Let us first show that $P\subseteq \decreasing(\kappa_P)$.
Take any point $(x,y)\in P$ and any point
$(x',y')\in\Omega\setminus\{(x,y)\}$ such that $x'\le x$ and $y'\le y$.
Let $Q$ be an increasing subset of
$P\cap\bigl((-\infty,x']\times(-\infty,y']\bigr)$ of cardinality
$\kappa_P(x',y')$. Since no two points in $P$ have the same $x$- or $y$-coordinates, $Q\cup\{(x,y)\}$ is an increasing subset of
$P\cap\bigl((-\infty,x]\times(-\infty,y]\bigr)$ of cardinality
$\kappa_P(x',y')+1$, and it follows that $\kappa_P(x',y')<\kappa_P(x,y)$.
This shows that $(x,y)\in \decreasing(\kappa_P)$, and since
$(x,y)$ was chosen arbitrarily in $P$, we conclude that
$P\subseteq \decreasing(\kappa_P)$.

To show that $\decreasing(\kappa_P)\subseteq P$, take any $(x,y)\in \decreasing(\kappa_P)$
and let $Q$ be an increasing subset of
$P\cap\bigl((-\infty,x]\times(-\infty,y]\bigr)$ of cardinality
$\kappa_P(x,y)$. Let $(x',y')$ be the point in $Q$ with maximal coordinates. Then, $\kappa_P(x',y')\ge\kappa_P(x,y)$ and,
since $(x,y)\in \decreasing(\kappa_P)$, it follows that $(x',y')=(x,y)$ and
hence that $(x,y)\in P$.
\end{proof}

The next lemma is essentially a reformulation of \cref{pr:beta} in
terms of the $D$ operator.
\begin{lemma}\label{lm:parallelogrambounds}
Let $\Omega$ be the open parallelogram
\[
\{(x,y)\in\bbR^2\ :\ \abs{ax+by}<1,\ \abs{ax-by}<\beta\}
\]
for some $a,b,\beta>0$,
and let $\rho$ be constant on $\Omega$. Let $u_{\rm linear}(x,y)=c(ax+by)$ for some $c\ge0$, and let $\{\sigma_\gamma\}_{\gamma>0}$ be Poisson point
processes on $\Omega$ with intensities $\gamma\rho$.
Then the following two statements hold \aas\ as $\gamma\rightarrow\infty$.
\begin{align*}
\forall w\in \calU_{2c}(\Omega),\,
\kors (\decreasing(w\sqrt\gamma)\cap\sigma_\gamma)/\gamma &\le \functional_\rho(u_{\rm linear})+3(\rho+1)\mu(\Omega)/\beta\\
\forall d\in\bbR\ \exists w\in \calU_{d,2c}(\Omega):\,
\kors (\decreasing(w\sqrt\gamma)\cap\sigma_\gamma)/\gamma &\ge \functional_\rho(u_{\rm linear})-3(\rho+1)\mu(\Omega)/\beta
\end{align*}
\end{lemma}
\begin{proof}
For any $w\in\calU_{2c}(\Omega)$, by \cref{lm:IandkappaIudecreasing},
$\decreasing(w\sqrt\gamma)\cap\sigma_\gamma$ is a $(\lfloor 2c\sqrt{\gamma}\rfloor+1)$-decreasing
subset of $\sigma_\gamma$.
Let $c'=c + \frac{1}{2\sqrt{\gamma}}$ and
let $u'_{\rm linear}(x,y)=c'(ax+by)$. Then, for any $w\in\calU_{2c}(\Omega)$,
$\decreasing(w\sqrt\gamma)\cap\sigma_\gamma$ is a $\lfloor 2c'\sqrt{\gamma}\rfloor$-decreasing
subset of $\sigma_\gamma$, and it follows from \cref{pr:beta} that the statement
\[
\forall w\in\calU_{2c}(\Omega),\ \ \kors(\decreasing(w\sqrt\gamma)\cap\sigma_\gamma)/\gamma \le \functional_\rho(u'_{\rm linear}) + 3\rho\mu(\Omega)/\beta
\]
holds \aas\ as $\gamma\rightarrow\infty$.

Note that $c'\rightarrow c$ as $\gamma\rightarrow\infty$.
By \cref{lm:Lcontinuous}, $L$ is continuous, so
\[
\functional_\rho(u_{\rm linear}) - \functional_\rho(u'_{\rm linear})
= \mu(\Omega) \bigl(L(\rho, c^2ab) - L(\rho, c'^2ab)\bigr)
\rightarrow0
\]
as $\gamma\rightarrow\infty$, and we conclude that the statement
\[
\forall w\in\calU_{2c}(\Omega),\ \kors(\decreasing(w\sqrt\gamma)\cap\sigma_\gamma)/\gamma \le \functional_\rho(u_{\rm linear}) + 3(\rho+1)\mu(\Omega)/\beta
\]
holds \aas\ as $\gamma\rightarrow\infty$.

Now for the second part.
This time, let $c'=c - \frac{1}{2\sqrt{\gamma}}$ and, as before,
let $u'_{\rm linear}(x,y)=c'(ax+by)$. We will soon let $\gamma$ tend
to infinity, so we may assume that $c'>0$.
Let $P_\gamma$ be any maximal $\lfloor 2c'\sqrt{\gamma}\rfloor$-decreasing subset of $\sigma_\gamma$.
It follows from \cref{pr:beta} that
\[
\kors P_\gamma/\gamma \ge \functional_\rho(u'_{\rm linear})-3\rho\mu(\Omega)/\beta
\]
\aas\ as $\gamma\rightarrow\infty$. Analogously to above,
since $c'\rightarrow c$
as $\gamma\rightarrow\infty$ and since $L$ is continuous, it follows
that
\[
\kors P_\gamma/\gamma \ge \functional_\rho(u_{\rm linear})-3(\rho+1)\mu(\Omega)/\beta
\]
\aas\ as $\gamma\rightarrow\infty$.
For any $d\in\bbR$, let $w=(\kappa_{P_\gamma}+\lceil d\sqrt\gamma\rceil)/\sqrt\gamma$.
By \cref{pr:decreasingincreasingrelation},
\[
0\le\kappa_{P_\gamma}\le\lfloor 2c'\sqrt{\gamma}\rfloor\le 2c'\sqrt{\gamma}= 2c\sqrt\gamma - 1,
\]
so $w$ belongs to $\calU_{d,2c}(\Omega)$.
Clearly, $D$ is invariant under translation by an integer, so
$\decreasing(w\sqrt\gamma)=\decreasing(\kappa_{P_\gamma})$ which equals $P_\gamma$
by \cref{lm:IandkappaIkappaPisP}. (Note that the points in $P_\gamma$
are in general position almost surely.)
We conclude that the statement
\[
\forall d\in\bbR\ \exists w\in \calU_{d,2c}(\Omega):\
\kors (\decreasing(w\sqrt\gamma)\cap\sigma_\gamma)/\gamma \ge \functional_\rho(u_{\rm linear})-3(\rho+1)\mu(\Omega)/\beta
\]
holds \aas\ as $\gamma\rightarrow\infty$.
\end{proof}

\section{Approximating \texorpdfstring{$\Omega$}{Omega} by a collection of parallelograms}
\label{sec:greatlemma}
Now when we have studied the behavior of a parallelogram density domain,
it is time to divide a general density domain into many small local parallelograms. It is vital that the number of such parallelograms is finite
since we want to infer an ``in probability'' result for the whole domain
from similar results for each local parallelogram. To this end we will rely heavily on the theory of Vitali coverings.

First a pair of technical lemmas.
\begin{lemma}\label{lm:disjointdistance}
Let $A$ and $B$ be disjoint closed subsets of $\bbR^2$,
and suppose $A$ is compact.
Then, the distance between $A$ and $B$ is positive.
\end{lemma}
\begin{proof}
Suppose not. Then there are sequences $a_i\in A$ and $b_i\in B$ such that
$\abs{a_i-b_i}\rightarrow0$. Since $A$ is compact, there is a subsequence
$a_{i_j}$ of $a_i$ that converges to some $a$ in $A$. By the triangle
inequality, $\abs{b_{i_j}-a}\le\abs{b_{i_j}-a_{i_j}}+\abs{a_{i_j}-a}$
which tends to zero. Since $B$ is closed, this implies that $a$ belongs
to $B$ which is a contradiction since $A$ and $B$ are disjoint.
\end{proof}

\begin{lemma}\label{lm:diamconvergence}
Let $\Omega$ be an open subset of $\bbR^2$ and let $C$ be a compact
subset of $\Omega$. Then there is a constant $K$ such that
\[
\diam w(C) \le \diam u(\Omega) + K\norm{w-u}_\Omega
\]
for any $u,w\in\calU(\Omega)$.
\end{lemma}
\begin{proof}
By \cref{lm:disjointdistance}, the distance between
$C$ and $\bbR^2\setminus\Omega$ is positive, so there exists
a $d>0$ such that for each $(x,y)\in C$ we have
$[x-d,x+d]\times[y-d,y+d]\subset\Omega$.
It follows that, for any $u,w\in\calU(\Omega)$,
\[
\norm{w-u}_\Omega
\ge
\sup_{(x,y)\in C}\norm{w-u}_{[x,x+d]\times[y,y+d]}
\ge \bigl(\sup w(C)-\sup u(\Omega)\bigr)d^2
\]
and
\[
\norm{w-u}_{\Omega}
\ge
\sup_{(x,y)\in C}\norm{w-u}_{[x-d,x]\times[y-d,y]}
\ge \bigl(\inf u(\Omega)-\inf w(C)\bigr)d^2.
\]
Thus we can choose $K:=2/d^2$.
\end{proof}

Next, we make the idea of a local parallelogram precise.
\begin{defi}
Let $u\in\calU(\Omega)$ and let $\iota>0$.
A \emph{$(u,\iota)$-parallelogram} is a closed parallelogram of the form
\begin{align*}
P=\{(x,y)\in\bbR^2\,:\,
&\abs{\tilde{u}^P_x(x-x_P)+\tilde{u}^P_y(y-y_P)} \le \iota c_P,\\
&\abs{\tilde{u}^P_x(x-x_P)-\tilde{u}^P_y(y-y_P)} \le c_P\},
\end{align*}
where $(x_P,y_P)$ is a point in $\Diff(u)$, $c_P>0$,
$\tilde{u}^P_x:=\max\{\iota^3, u_x(x_P,y_P)\}$
and
$\tilde{u}^P_y:=\max\{\iota^3, u_y(x_P,y_P)\}$.

Also, for notational convenience, define
$u^P_x=u_x(x_P,y_P)$ and $u^P_y=u_y(x_P,y_P)$.

We say that $P$ is \emph{well behaved} if
$\tilde{u}^P_x=u^P_x$ and $\tilde{u}^P_y=u^P_y$.
\end{defi}

\begin{lemma}\label{lm:uiparallelogrambound}
For any point $(x,y)$ in a $(u,\iota)$-parallelogram $P$, we have
\[
\abs{x-x_P} + \abs{y-y_P} \le c_P(1+\iota)\iota^{-3}
\]
and
\[
\abs{u^P_x(x-x_P)+u^P_y(y-y_P)} \le c_P(1+\iota).
\]
\end{lemma}
\begin{proof}
We have
\begin{equation}\label{eq:utildex}
\begin{split}
& \quad 2\tilde{u}^P_x\abs{x-x_P}\\
&\le\abs{\tilde{u}^P_x(x-x_P) + \tilde{u}^P_y(y-y_P)}
+\abs{\tilde{u}^P_x(x-x_P) - \tilde{u}^P_y(y-y_P)}\\
&\le c_P (1+\iota),
\end{split}
\end{equation}
where the first inequality is the triangle inequality and the second
inequality follows from the definition of a $(u,\iota)$-parallelogram.
For the same reason,
\begin{equation}\label{eq:utildey}
2\tilde{u}^P_y\abs{y-y_P} \le c_P (1+\iota),
\end{equation}
so
\[
\abs{x-x_P}+\abs{y-y_P} \le c_P (1+\iota)
\left(\frac1{2\tilde{u}^P_x}+\frac1{2\tilde{u}^P_y}\right)
\le c_P (1+\iota)\iota^{-3}.
\]

Finally, by \cref{eq:utildex,eq:utildey},
\[
\abs{u^P_x(x-x_P)+u^P_y(y-y_P)}
\le\tilde{u}^P_x\abs{x-x_P}+\tilde{u}^P_y\abs{y-y_P}
\le c_P(1+\iota).
\]
\end{proof}

Let us recall the definition of a regular Vitali covering.
As always, we let $\mu$ denote the Lebesgue measure on $\bbR^2$.
\begin{defi}
Let $A\subseteq\bbR^2$ and let $\calC$ be a collection
of closed subsets of $\bbR^2$.
\begin{itemize}
\item
$\calC$ is a \emph{Vitali covering} of $A$ if,
for any $p\in A$ and any $\delta>0$, there is a $C\in\calC$ such that
$p\in C$ and $0<\diam C<\delta$.
\item
$\calC$ is \emph{regular} if
there is a constant $K$ such that $(\diam C)^2\le K \mu(C)$ for any $C\in\calC$.
\end{itemize}
\end{defi}

\begin{lemma}\label{lm:vitali}
Let $u\in \calU(\Omega)$ and let $T$ be a subset
of $\Diff(u)$ where $u_x$ and $u_y$ are both bounded.
Then, for any $\iota>0$, the family of all $(u,\iota)$-parallelograms is a regular Vitali covering of $T$.
\end{lemma}
\begin{proof}
For any $p\in T$, the diameter of
a $(u,\iota)$-parallelogram $P$ centered at $(x_P,y_P)=p$ is bounded by
\[
2\sqrt{(x-x_P)^2+(y-y_P)^2}\le 2\bigl(\abs{x-x_P}+\abs{y-y_P}\bigr),
\]
which is at most $2c_P(1+\iota)\iota^{-3}$ by \cref{lm:uiparallelogrambound}. By choosing $c_P$ small enough we
can make the diameter arbitrarily small, so the family of
$(u,\iota)$-parallelograms is a Vitali covering of $A$.
To see that it is regular, note that
$
\mu(P)=2\iota c_P^2/\tilde{u}^P_x\tilde{u}^P_y,
$
so the quotient
\[
\frac{\diam P}{\sqrt{\mu(P)}}
\le \frac{2c_P(1+\iota)\iota^{-3}}{\sqrt{\mu(P)}}
\]
is bounded since $u_x$ and $u_y$ are bounded on $T$.
\end{proof}

The following lemma divides a density domain with a doubly increasing
function $u$ into a finite number of local parallelograms such that,
within each parallelogram, the situation is close to the condition
of \cref{lm:parallelogrambounds}, that is, the density is nearly constant and $u$ is nearly a linear function aligned with the
parallelogram.

We will use the ordo notation $o_\iota(1)$ to represents a
function of $\iota$ that tends to zero as $\iota$ tends to zero.

\begin{lemma}\label{lm:great}
Let $(\Omega,\rho)$ be a density domain and let $u\in \calU(\Omega)$ and $\eps>0$.
Then, for any $0<\iota<1$,
there is a measurable set $S_\iota\subseteq \Diff(u)$
and a finite disjoint collection $\calP_\iota$ of
$(u,\iota)$-parallelograms such that the following holds.
\begin{thmlist}
\item\label{lm:greatmisc}
For each $P\in\calP_\iota$ it holds that $P\subset\Omega$, $(x_P,y_P)\in S_\iota$ and $c_P<1$. Also, $S_\iota\subseteq\bigcup\calP_\iota$.
\item\label{lm:greatbounded}
$\rho$, $u_x$ and $u_y$ are bounded on $\bigcup_{0<\iota<1}S_\iota$.
\item\label{lm:greatparallelogramscovereverything}
$\norm{\rho}_{\Omega\setminus S_\iota}<\eps + o_\iota(1)$.
\intuition{$S_\iota$ nearly covers the density domain.}
\item\label{lm:greatScoversparallelograms}
For each $P\in\calP_\iota$ it holds that
$\mu(P\cap S_\iota)/\mu(P)>1-\iota$.
\intuition{$S_\iota$ nearly covers each parallelogram.}
\item\label{lm:greatrhoconstant}
For each $P\in\calP_\iota$
it holds that $\abs{\rho(x,y)-\rho(x_P,y_P)}\le\iota$
for any $(x,y)\in P\cap S_\iota$.
\intuition{$\rho$ is nearly constant on each parallelogram.}
\item\label{lm:greatulinear}
For each $P\in\calP_\iota$ it holds that
\[
\abs{u(x,y)-\bigl(u(x_P,y_P)+u^P_x(x-x_P)
+u^P_y(y-y_P)\bigr)} \le \iota^5(\abs{x-x_P}+\abs{y-y_P})
\]
for any $(x,y)\in P$.
\intuition{$u$ is nearly linear on each parallelogram.}
\item\label{lm:greatLconstant}
For each $P\in\calP_\iota$ it holds that
\[
\abs{\norm{L(\rho,u_xu_y)}_P/\mu(P)-L(\rho(x_P,y_P),u^P_xu^P_y)}<\iota.
\]
\intuition{$L(\rho,u_xu_y)$ is nearly constant on each parallelogram.}
\item\label{lm:greatuvarieslittle}
For each $P\in\calP_\iota$ it holds that
\[
\sup_{(x,y)\in P}\abs{u(x,y)-u(x_P,y_P)}<
\begin{cases}
c_P\iota(1+5\iota) & \text{if
$P$ is well behaved,} \\
c_P(1+7\iota) & \text{otherwise.}
\end{cases}
\]
\intuition{$u$ does not vary too much inside each parallelogram.}
\item\label{lm:greatwvarieslittle}
There is a function $d:(0,1)\rightarrow\bbR_{>0}$ such that
\[
\sup_{w\in\calU(\Omega):\ \norm{w-u}_\Omega<d(\iota)}\ \sup_{P\in\calP_\iota}
\left(\left(\frac{\diam w(P)}{2\iota c_P}\right)^2\tilde{u}^P_x\tilde{u}^P_y-u^P_x u^P_y\right) < o_\iota(1).
\]
\end{thmlist}
\end{lemma}
\begin{proof}
By \cref{lm:monotonic}, $\mu(\Omega\setminus \Diff(u))=0$.

Since $\norm{\rho}_{\Omega}<\infty$, we can choose a subset
$T$ of $\Diff(u)$ with finite measure such that $\rho$,
$u_x$ and $u_y$ are all smaller than some positive constant
$C$ there,
and such that
\begin{equation}\label{eq:Tlarge}
\norm{\rho}_{\Omega\setminus T} < \eps.
\end{equation}
For each $0<\iota<1$,
let $T_\iota$ be the set of points in $T$ that are
Lebesgue points of $L(\rho,u_xu_y)$
with respect to the family of $(u,\iota)$-parallelograms, that is, $T_\iota$
is the set of points $(x_0,y_0)\in T$ such that for each $\eps'>0$
we have
$\abs{\norm{L(\rho,u_xu_y)}_P/\mu(P)
- L(\rho,u_xu_y)|_{(x_0,y_0)}}<\eps'$
for all sufficiently small $(u,\iota)$-parallelograms $P$ centered at $(x_0,y_0)$.
By \cref{lm:vitali}, the family of $(u,\iota)$-parallelograms is a regular Vitali
covering of $T$, so,
by Lebesgue's differentiation theorem (see e.g.~\cite{Folland}), $\mu(T_\iota)=\mu(T)$ for any $\iota$.

For any $0<\iota<1$ and any positive integer $j$, let
\[
\tilde{S}_\iota^j=\{(x,y)\in T_\iota\,:\,(j-1)\iota\le\rho(x,y)<j\iota\}
\]
and let $S_\iota^j$ be the set of points in $\tilde{S}_\iota^j$
at which the density of $\tilde{S}_\iota^j$ is $1$. By Lebesgue's density
theorem, $\mu(\tilde{S}_\iota^j\setminus S_\iota^j)=0$.
For any $\iota$, $T_\iota$ is the union of a finite number of sets of the form
$\tilde{S}_\iota^j$, so the union $\hat{S}_\iota:=\bigcup_j S_\iota^j$ of all $S_\iota^j$ for
a fixed $\iota$ has the same measure as $T_\iota$ and hence as $T$.

For any $0<\iota<1$ and any positive integer $j$, let
$\calA_\iota^j$ be the family of $(u,\iota)$-parallelograms $P$ with center
in $S_\iota^j$ and $c_P<1$ such that
\begin{enumerate}[label=\upshape(\Roman*),ref=(\Roman*)]
\item\label{it:a}
$\mu(P\cap S_\iota^j)/\mu(P)>1-\iota$,
\item\label{it:b}
$\abs{\norm{L(\rho,u_xu_y)}_P/\mu(P) - L(\rho(x_P,y_P),u^P_xu^P_y)}<\iota$, and
\item\label{it:c}
the $(u,\iota)$-parallelogram $P'$ concentric with $P$
but with $c_{P'}=(1+\iota)c_P$ is contained in $\Omega$, and
\[
\abs{u(x,y)
-\bigl(u(x_P,y_P)+u^P_x(x-x_P)+u^P_y(y-y_P)\bigr)}
\le\iota^5\bigl(\abs{x-x_P}+\abs{y-y_P}\bigr)
\]
for any point $(x,y)\in P'$.
\end{enumerate}
Let $\calA_\iota=\bigcup_j\calA_\iota^j$.

We claim that $\calA_\iota^j$ is a Vitali covering
of $S_\iota^j$. To see this, take any $(x,y)\in S_\iota^j$ and note
the following:
\begin{itemize}
\item
Since the density of $\tilde{S}_\iota^j$ is 1 at $(x,y)$ (by the choice
of $S_\iota^j$) and $\mu(\tilde{S}_\iota^j\setminus S_\iota^j)=0$,
\ref{it:a} holds for any sufficiently small $(u,\iota)$-parallelogram
centered at $(x,y)$.
\item
Since $(x,y)$ is a Lebesgue point of $T$ (by the choice of $T_\iota$),
\ref{it:b} holds for any sufficiently small $(u,\iota)$-parallelogram
centered at $(x,y)$.
\item
Since $\Omega$ is open and $u$ is differentiable at $(x,y)$,
\ref{it:c} holds for any sufficiently small $(u,\iota)$-parallelogram
centered at $(x,y)$.
\end{itemize}
By \cref{lm:vitali}, the family of all $(u,\iota)$-parallelograms
is regular, so it follows that
$\calA_\iota^j$ is a regular Vitali covering of $S_\iota^j$,
and hence $\calA_\iota$ is a regular Vitali covering
of $\hat{S}_\iota$.
By Vitali's cover theorem (see e.g.~\cite{EvansGariepyBook}), there is a finite disjoint subfamily
$\calP_\iota$ of $\calA_\iota$ such that 
$\mu(\hat{S}_\iota\cap\cup\calP_\iota)\ge(1-\iota)\mu(\hat{S}_\iota)$ and hence
\begin{equation}\label{eq:muSP}
\mu(T\cap\cup\calP_\iota)\ge(1-\iota)\mu(T).
\end{equation}
Finally, define
\[
S_\iota:=
\bigcup_{j=1}^\infty\bigcup_{P\in\calP_\iota\cap\calA_\iota^j}P\cap S_\iota^j.
\]

Let us check that $S_\iota$ and $\calP_\iota$ have the properties claimed in the lemma.

\smallskip
\reflocal{lm:greatmisc} and \reflocal{lm:greatbounded} follow directly from the definitions.

\smallskip
\reflocal{lm:greatparallelogramscovereverything}
%\item\label{it:ZERO}
We have
\begin{multline*}
\mu(S_\iota)
=\sum_{j=1}^\infty\sum_{P\in\calP_\iota\cap\calA_\iota^j}\mu(P\cap S_\iota^j)
\ge \{\text{by \ref{it:a}}\}
\ge\sum_{j=1}^\infty\sum_{P\in\calP_\iota\cap\calA_\iota^j}(1-\iota)\mu(P)\\
=(1-\iota)\mu(\cup\calP_\iota)
\ge \{\text{by \cref{eq:muSP}}\} \ge(1-\iota)^2\mu(T),
\end{multline*}
so \reflocal{lm:greatparallelogramscovereverything} follows from \cref{eq:Tlarge}
and the fact that $\rho$ is bounded in $T$.

\smallskip
\reflocal{lm:greatScoversparallelograms} follows from \ref{it:a}.

\smallskip
\reflocal{lm:greatrhoconstant} follows from the definition of the sets $\tilde{S}_\iota^j$.

\smallskip
\reflocal{lm:greatulinear} follows from \ref{it:c}.

\smallskip
\reflocal{lm:greatLconstant} follows from \ref{it:b}.

\smallskip
\reflocal{lm:greatuvarieslittle} requires some reasoning.
From \ref{it:c} we know that, for each $P\in\calP_\iota$,
\begin{equation}\label{eq:uvarieslittle}
\abs{u(x,y)-\bigl(u(x_P,y_P)+u^P_x(x-x_P)
+u^P_y(y-y_P)\bigr)} \le \iota^5(\abs{x-x_P}+\abs{y-y_P})
\end{equation}
for any $(x,y)\in P'$, where $P'\subset\Omega$ is defined as in \ref{it:c}.

Consider a parallelogram $P$ in $\calP_\iota$.
By \cref{lm:uiparallelogrambound}, inside $P'$,
\begin{equation}\label{eq:absxplusabsy}
\abs{x-x_P}+\abs{y-y_P}\le c_{P'} (1+\iota)\iota^{-3}.
\end{equation}
Also, inside $P'$,
\begin{equation}\label{eq:udiam}
\begin{split}
&\quad\abs{u(x,y) - u(x_P,y_P)} \le\{\text{triangle ineq.}\}\le\\
&\le \abs{u^P_x(x-x_P)+u^P_y(y-y_P)}\\
&\quad + \abs{u(x,y) - (u(x_P,y_P)+u^P_x(x-x_P)+u^P_y(y-y_P))}\\
&\le \abs{u^P_x(x-x_P)+u^P_y(y-y_P)} + \iota^5(\abs{x-x_P}+\abs{y-y_P})\\
&\le\abs{u^P_x(x-x_P)+u^P_y(y-y_P)} + \iota^2(1+\iota)c_{P'},
\end{split}
\end{equation}
where the last inequality follows from \cref{eq:absxplusabsy}.

If $P$ is well behaved, inside $P'$ we have
\[
\abs{u^P_x(x-x_P)+u^P_y(y-y_P)} = \abs{\tilde{u}^P_x(x-x_P)+\tilde{u}^P_y(y-y_P)}
\le \iota c_{P'}
\]
by the definition of a $(u,\iota)$-parallelogram. Hence, by \cref{eq:udiam}, inside $P'$,
\[
\abs{u(x,y) - u(x_P,y_P)} \le c_{P'}\iota(1+\iota+\iota^2)
= c_P(1+\iota)\iota(1+\iota+\iota^2)< c_P\iota(1+5\iota),
\]
where the last inequality follows from the
fact that $\iota < 1$.

If $P$ is not well behaved, by \cref{lm:uiparallelogrambound},
at least we have
\[
\abs{u^P_x(x-x_P)+u^P_y(y-y_P)}\le c_{P'}(1+\iota).
\]
Hence, by \cref{eq:udiam}, inside $P'$ we have
\[
\abs{u(x,y) - u(x_P,y_P)} \le c_{P'}(1+\iota)(1+\iota^2)
=c_P(1+\iota)^2(1+\iota^2)<c_P(1+7\iota),
\]
where we have used again than $\iota<1$.

\smallskip
\reflocal{lm:greatwvarieslittle} requires some reasoning as well.
Consider a parallelogram $P$ in $\calP_\iota$. Let $P'$ be the larger
concentric $(u,\iota)$-parallelogram as defined in \ref{it:c}.
By \cref{lm:diamconvergence} applied
to the open set $\interior{P'}$ and the compact set $P$
there is a $K_P>0$ such that
$\diam w(P) \le \diam u(P')+K_P\norm{w-u}_\Omega$ for any $w\in\calU(\Omega)$.

It follows from \reflocal{lm:greatuvarieslittle} that,
if $P$ is well behaved,
\[
\diam w(P) \le 2\iota c_P(1+5\iota) + K_P\norm{w-u}_\Omega
\le 2\iota c_P(1+6\iota)
\]
for any $w\in\calU(\Omega)$ with $\norm{w-u}_\Omega\le 2\iota^2 c_P/K_P$,
and if $P$ is not well behaved,
\[
\diam w(P) \le 2c_P(1+7\iota) + K_P\norm{w-u}_\Omega
\le 2c_P(1+7\iota+\iota^2)<18c_P
\]
for any $w\in\calU(\Omega)$ with $\norm{w-u}_\Omega\le 2\iota^2 c_P/K_P$.

Thus, if $P$ is well behaved, for any
$w\in\calU(\Omega)$ with $\norm{w-u}_\Omega\le 2\iota^2 c_P/K_P$ we have
\begin{multline*}
\left(\frac{\diam w(P)}{2\iota c_P}\right)^2
\tilde{u}^P_x\tilde{u}^P_y-u^P_xu^P_y
=\left[\left(\frac{\diam w(P)}{2\iota c_P}\right)^2-1\right]
u^P_x u^P_y\\
\le [(1+6\iota)^2-1]
u^P_x u^P_y
= o_\iota(1)
\end{multline*}
since $u_x$ and $u_y$ are bounded in $\bigcup_\iota S_\iota$.

If $P$ is not well behaved, at least one of $u^P_x$ and $u^P_y$ is smaller
than $\iota^3$, and at least one of $\tilde{u}^P_x$ and $\tilde{u}^P_y$ equals $\iota^3$. It follows that, for any
$w\in\calU(\Omega)$ with $\norm{w-u}_\Omega\le 2\iota^2 c_P/K_P$ we have
\[
\left(\frac{\diam w(P)}{2\iota c_P}\right)^2
\tilde{u}^P_x\tilde{u}^P_y-u^P_xu^P_y
\le(9/\iota)^2
\tilde{u}^P_x\tilde{u}^P_y-u^P_xu^P_y=o_\iota(1).
\]

We conclude that
\[
\left(\frac{\diam w(P)}{2\iota c_P}\right)^2
\tilde{u}^P_x\tilde{u}^P_y-u^P_xu^P_y < o_\iota(1)
\]
for any
$w\in\calU(\Omega)$ with $\norm{w-u}_\Omega\le 2\iota^2 c_P/K_P$
whether $P$ is well behaved or not. Thus, we can choose
$d(\iota):=2\iota^2\min_{P\in\calP_\iota}(c_P/K_P)$.
\end{proof}

\section{Semicontinuity of \texorpdfstring{$\functional_\rho$}{the functional} and a related probabilistic result}\label{sec:semicontinuity}
Our proof of \cref{th:mainF} will rely heavily on the following
result.
\begin{prop}\label{pr:semicontinuity}
$\functional_\rho$ is upper semicontinuous in the $L^1(\Omega)$-norm.
\end{prop}
\begin{proof}
Let $u\in \calU(\Omega)$.
We must show that, for any $\eps'>0$ there is a $\delta>0$ such that
$\functional_\rho(w)<\functional_\rho(u)+\eps'$ for any $w\in\calU(\Omega)$ with
$\norm{w-u}_\Omega<\delta$.

Choose any $\eps>0$ smaller than $\eps'$ and apply \cref{lm:great}.
Let $d:(0,1)\rightarrow\bbR_{>0}$ be the function defined in
\cref{lm:greatwvarieslittle}
and consider any family $\{w^{(\iota)}\in\calU(\Omega)\}_{0<\iota<1}$ such that $\norm{w^{(\iota)}-u}_\Omega < d(\iota)$.
Let
\[
q^{(\iota)}:=
\sup_{w\in\calU(\Omega):\ \norm{w-u}_\Omega<d(\iota)}\ \sup_{P\in\calP_\iota}
\left(\left(\frac{\diam w(P)}{2\iota c_P}\right)^2
\tilde{u}^P_x\tilde{u}^P_y-u^P_xu^P_y\right).
\]
Consider any $P\in \calP_\iota$ and let $\rho_P:=\sup_{S_\iota\cap P}\rho$.
By \cref{lm:Lcontinuous}, $L$ is increasing in the first variable,
so
\[
\norm{L(\rho,w^{(\iota)}_xw^{(\iota)}_y)}_{P\cap S_\iota} \le 
\norm{L(\rho_P,w^{(\iota)}_xw^{(\iota)}_y)}_P.
\]
Let $r_P:=\diam w^{(\iota)}(P)$.
By \cref{pr:parallelogrammaximizer},
on $\calU_{r_P}(\interior P)$, $\functional_{\rho_P}$ is maximized by the
function
\[
u^P_{\rm linear}(x,y)
=\frac{r_P}{2\iota c_P}\bigl(\tilde{u}^P_x(x-x_P)+\tilde{u}^P_y(y-y_P)\bigr),
\]
so
\begin{multline*}
\norm{L(\rho_P,w^{(\iota)}_xw^{(\iota)}_y)}_P
\le \functional_{\rho_P}(u^P_{\rm linear})
=\mu(P)L(\rho_P,(\tfrac{r_P}{2\iota c_P})^2\tilde{u}^P_x\tilde{u}^P_y)\\
\le\mu(P) L(\rho_P,q^{(\iota)}+u^P_xu^P_y)=:\text{RHS},
\end{multline*}
where the last inequality uses the fact that $L$ is increasing in the
second variable (\cref{lm:Lcontinuous}).
By \cref{lm:greatrhoconstant},
$\abs{\rho_P-\rho(x_P,y_P)}\le\iota$,
and by \cref{lm:greatwvarieslittle},
$q^{(\iota)}<o_\iota(1)$.
By \cref{lm:Lcontinuous}, $L$ is continuous and hence
uniformly continuous on the set
$\{(\rho(x,y),u_x(x,y)u_y(x,y))\,:\,(x,y)\in \bigcup_{\iota>0}S_\iota\}$
which is bounded by \cref{lm:greatbounded}.
Again by \cref{lm:Lcontinuous}, $L$ is increasing in the second variable, so
$\text{RHS} < \mu(P)\bigl(L(\rho(x_P,y_P),u^P_xu^P_y)+o_\iota(1)\bigr)$,
where $o_\iota(1)$ is independent of $P$.
By \cref{lm:greatLconstant},
\[
L(\rho(x_P,y_P),u^P_xu^P_y) <
\tfrac{1}{\mu(P)}\norm{L(\rho,u_xu_y)}_P+\iota,
\]
so
\[
\norm{L(\rho,w^{(\iota)}_xw^{(\iota)}_y)}_{P\cap S_\iota}
< \norm{L(\rho,u_xu_y)}_P+o_\iota(1)\mu(P).
\]
By \cref{lm:greatmisc}, $S_\iota\subseteq\bigcup\calP_\iota$, so summing over all $P$ in $\calP_\iota$ yields
\[
\norm{L(\rho,w^{(\iota)}_xw^{(\iota)}_y)}_{S_\iota}-\norm{L(\rho,u_xu_y)}_{\cup\calP_\iota}
<o_\iota(1)\mu(\cup\calP_\iota)=o_\iota(1).
\]
By \cref{lm:Lcontinuous}, $L(\rho,w^{(\iota)}_xw^{(\iota)}_y)\le\rho$,
so by \cref{lm:greatparallelogramscovereverything} it now follows that
\[
\norm{L(\rho,w^{(\iota)}_xw^{(\iota)}_y)}_\Omega\le
\norm{L(\rho,u_xu_y)}_\Omega+\eps + o_\iota(1).
\]
Since $\eps<\eps'$, the lemma follows.
\end{proof}

In the proof of \cref{th:main}, we will need the following probabilistic analogue to \cref{pr:semicontinuity}.
\begin{lemma}\label{lm:notbetterthanfunctional}
Let $(\Omega,\rho)$ be a density domain, and let $\{\sigma_\gamma\}_{\gamma>0}$ be Poisson point
processes on $\Omega$ with intensity functions $\gamma\rho$.
Then, for any $u\in \calU(\Omega)$ and any $\eps'>0$ there is
a $\delta>0$ such that
\[
\sup_{w\in \calU(\Omega),\ \norm{w-u}_\Omega<\delta}\kors (\decreasing(w\sqrt\gamma)\cap\sigma_\gamma)/\gamma < \functional_\rho(u)+\eps'
\]
holds \aas\ as $\gamma\rightarrow\infty$.
\end{lemma}
\begin{proof}
Choose any $\eps>0$ smaller than $\eps'/2$ and apply \cref{lm:great}.
Let $d:(0,1)\rightarrow\bbR_{>0}$ be the function defined in
\cref{lm:greatwvarieslittle}.
Consider any $P\in\calP_\iota$ and any $w\in\calU(\Omega)$ such that
$\norm{w-u}_\Omega<d(\iota)$, and let $r_P:=\diam w(P)$ and
$\rho_P:=\sup_{S_\iota\cap P}\rho$.
Define $u^P_{\rm linear}\in\calU_{r_P}(\interior P)$ by
\[
u^P_{\rm linear}(x,y)
=\frac{r_P}{2\iota c_P}\bigl(\tilde{u}^P_x(x-x_P)+\tilde{u}^P_y(y-y_P)\bigr).
\]
Since $\sigma_\gamma\cap P\cap S_\iota$ is a subset of a Poisson point
process on $P$ with homogeneous intensity $\rho_P$,
\cref{lm:parallelogrambounds} yields that, \aas\ as $\gamma\rightarrow\infty$,
\[
\tfrac{1}{\gamma}\kors (\decreasing(w\sqrt\gamma)\cap\sigma_\gamma\cap P\cap S_\iota)
\le \functional_{\rho_P}(u^P_{\rm linear})+3\iota\mu(P)(\rho_P+1)=:\text{RHS}.
\]
As in the proof of the upper semicontinuity, we obtain
\[
\text{RHS} <
\norm{L(\rho,u_xu_y)}_P + \mu(P)o_\iota(1),
\]
where $o_\iota(1)$ is independent of $P$.
Summing over all $P$ in $\calP_\iota$ yields that, \aas\ as $\gamma\rightarrow\infty$,
for any $w\in\calU(\Omega)$ such that $\norm{w-u}_\Omega<d(\iota)$,
\[
\tfrac{1}{\gamma}\kors (\decreasing(w\sqrt\gamma)\cap\sigma_\gamma\cap S_\iota)
-\norm{L(\rho,u_xu_y)}_{\cup\calP_\iota}
<o_\iota(1)\mu(\cup\calP_\iota)=o_\iota(1).
\]
By the law of large numbers, \aas\ as $\gamma\rightarrow\infty$ we have
$\frac{1}{\gamma}\kors(\sigma_\gamma\setminus S_\iota)<
\norm{\rho}_{\Omega\setminus S_\iota}+\eps$
which is bounded by $2\eps+o_\iota(1)$ by \cref{lm:greatparallelogramscovereverything}, so
\aas
\[
\sup_{w\in\calU(\Omega):\ \norm{w-u}_\Omega<d(\iota)}\tfrac{1}{\gamma}\kors (\decreasing(w\sqrt\gamma)\cap\sigma_\gamma)
-\norm{L(\rho,u_xu_y)}_{\cup\calP_\iota}
<2\eps+o_\iota(1).
\]
Since $\eps<\eps'/2$, the lemma now follows.
\end{proof}

\section{Compactness of the set of doubly increasing functions}
\label{sec:compactness}
In this section we put probabilistic matters aside and concern ourselves with the topology of the set of doubly increasing functions.

First, we show that doubly increasing functions can be extended to larger domains in a natural way.
\begin{lemma}\label{lm:increasingextension}
Let $A\subseteq B\subseteq\bbR^2$.
For any $u\in\calU(A)$ there is a $w\in\calU(B)$ such that
the restriction of $w$ to $A$ is $u$ and the images $u(A)$ and $w(B)$
have the same closure.
\end{lemma}
\begin{proof}
For each $(x,y)\in B$, let $P(x,y):=\{(x',y')\in A\,:\,x'\le x,\,y'\le y\}$.
Define $w$ by letting $w(x,y):=\sup_{P(x,y)} u$ if $P(x,y)$ is nonempty, and
$w(x,y):=\inf_A u$ if $P(x,y)$ is empty. It is straightforward to verify
that $w$ is doubly increasing.
\end{proof}

Our proof of the existence of maximizers of $\functional_\rho$
(\cref{th:mainFmaximizer}) will need two key ingredients.
The first one is the semicontinuity of
$\functional_\rho$ (\cref{pr:semicontinuity}) and the other one is the
following result about the topology of the set of doubly increasing functions, which will be essential also in the proof of \cref{th:main}.
\begin{prop}\label{pr:compactness}
If $\Omega\subseteq\bbR^2$ is open, then $\calU_{h,r}(\Omega)$ is a compact subset of $L^1(\Omega)$.
\end{prop}
\begin{proof}
Let $\{u_n\}_{n=1}^\infty$ be a sequence of elements in
$\calU_{h,r}(\Omega)$. We need to show that there is a convergent
subsequence.

Let $Q=\{q_1,q_2,\dotsc\}$ be a countable dense subset of $\Omega$.
We define a sequence $S_1,S_2,\dotsc$ of subsequences of
$\{u_n\}_{n=1}^\infty$ recursively as follows. First, let $S_1$ be the
original sequence $\{u_n\}_{n=1}^\infty$. Then, for $n=1,2,\dotsc$,
let $S_{n+1}$ be a subsequence of $S_n$ that converges at the point
$q_n$. This is possible since $[h,h+r]$ is a compact set.
Finally, construct another subsequence $\{w_n\}_{n=1}^\infty$
of $\{u_n\}_{n=1}^\infty$ by letting $w_n$ be the $n$th element
of $S_n$. Clearly, $\{w_n\}$ converges at each point of $Q$.
By \cref{lm:increasingextension}, we can choose a
$w\in\calU_{h,r}(\Omega)$ such that
$\lim_{n\rightarrow\infty}w_n(q)=w(q)$ for any $q\in Q$.
We claim that $\lim_{n\rightarrow\infty}w_n(p)=w(p)$ for any
continuity point $p$ of $w$. Take any $\eps>0$. Since $w$ is continuous
at $p$, we can pick a $\delta>0$ such that
$\abs{w(p')-w(p)}<\eps/2$ for any $p'$ with $\abs{p'-p}<\delta$.
Let $A^-:=\{p'\in\Omega\,:\,\abs{p'-p}<\delta,\ p'\ \text{strictly south-west of}\ p\}$
and $A^+:=\{p'\in\Omega\,:\,\abs{p'-p}<\delta,\ p'\ \text{strictly north-east of}\ p\}$.
Since $\Omega$ is open, $A^-$ and $A^+$ are both open and nonempty,
so there are $q^-\in Q\cap A^-$ and $q^+\in Q\cap A^+$.
For all sufficiently large $n$, we have $\abs{w_n(q^-)-w(q^-)}<\eps/2$
and $\abs{w_n(q^+)-w(q^+)}<\eps/2$, and hence
\begin{multline*}
w(p)-\eps\le w(q^-)-\eps/2\le w_n(q^-)\le w_n(p)\\
\le w_n(q^+)\le w(q^+)+\eps/2 \le w(p)+\eps.
\end{multline*}
Since $\eps$ was chosen arbitrarily, we conclude that
$\{w_n\}_{n=1}^\infty$ converges to $w$ at any continuity point of $w$.

By \cref{lm:monotonic}, $w$ is continuous almost everywhere,
so by the theorem of bounded convergence, $w_n$ converges to $w$ in
the $L^1(\Omega)$-norm.
\end{proof}

\section{A gluing lemma and our first main theorem}
\label{sec:mainproof}
To prove our first main theorem, \cref{th:main}, we need one final
lemma. We know that we can find large unions of decreasing subsets within
each local parallelogram. The following lemma glues those unions together to form a global union of decreasing subsets.
\begin{lemma}\label{lm:notworsethanfunctional}
Let $\{\sigma_\gamma\}_{\gamma>0}$ be Poisson point
processes on $\Omega$ with intensity functions $\gamma\rho$,
and let $r\ge0$.
For any $u\in \calU_{0,r}(\Omega)$ and any $\eps'>0$,
there is a family $\{w_\gamma\in\calU_{0,r}\}_{\gamma>0}$
(dependent on $\{\sigma_\gamma\}_{\gamma>0}$)
such that $w_\gamma\rightarrow u$ and, \aas\ as $\gamma\rightarrow\infty$,
\[
\kors (\decreasing(w_\gamma\sqrt\gamma)\cap\sigma_\gamma)/\gamma > \functional_\rho(u)-\eps'.
\]
\end{lemma}
\begin{proof}
Choose any $\eps>0$ smaller than $\eps'$ and apply \cref{lm:great}.
We will only consider $\iota$ in the interval $(0,1/2)$.

For each $P\in \calP_\iota$, let $P'$ denote an open parallelogram
with the same center as $P$ but a factor $1-2\iota$ as wide and high.
Let $\rho_P:=\inf_{S_\iota\cap P'}\rho$ and let
$\bar{u}_P\in \calU(P)$ be defined by
\[
\bar{u}_P(x,y):=u(x_P,y_P)+u^P_x(x-x_P)+u^P_y(y-y_P)
\]
if $P$ is well behaved and $\bar{u}_P(x,y):=u(x,y)$ otherwise.
Let
$\tau_P^{(\gamma)}$ be a Poisson point process on $P'\setminus S_\iota$
with homogenous intensity $\gamma\rho_P$. Then
$(\sigma_\gamma\cap P')\cup\tau_P^{(\gamma)}$ is a
superset of a Poisson point process 
$\tilde{\sigma}_\gamma$ on $P'$ with homogeneous intensity $\gamma\rho_P$.

Let $r_P:=\diam\bar{u}_P(P')$.
By \cref{lm:parallelogrambounds},
for any well-behaved $P\in\calP_\iota$ it holds \aas\ as $\gamma$ tends
to infinity that there is a $w_P^{(\gamma)}\in \calU_{u(x_P,y_P)-\frac{r_P}{2},r_P}(P')$ such that
\begin{equation}\label{eq:tildesigma}
\tfrac{1}{\gamma}\kors (\decreasing(w_P^{(\gamma)}\sqrt\gamma)\cap\tilde{\sigma}_\gamma)/\mu(P')
\ge L(\rho_P,u^P_xu^P_y)-3(\rho_P+1)\iota.
\end{equation}
Since $\tilde{\sigma}_\gamma\subseteq(\sigma_\gamma\cap P')\cup\tau^{(\gamma)}_P$, we have
\begin{multline}\label{eq:tildesigmatau}
\kors (\decreasing(w_P^{(\gamma)}\sqrt\gamma)\cap\tilde{\sigma}_\gamma)
\le\kors(\decreasing(w_P^{(\gamma)}\sqrt\gamma)\cap((\sigma_\gamma\cap P')\cup\tau_P^{(\gamma)}))\\
\le\kors(\decreasing(w_P^{(\gamma)}\sqrt\gamma)\cap\sigma_\gamma\cap P')+\kors\tau_P^{(\gamma)},
\end{multline}
and by the law of large numbers,
\begin{equation}\label{eq:tau}
\kors\tau_P^{(\gamma)} \le \gamma(\rho_P+1)\mu(P'\setminus S_\iota)
\end{equation}
\aas\ as $\gamma\rightarrow\infty$.
Combining \cref{eq:tildesigma,eq:tildesigmatau,eq:tau}, we obtain
\begin{multline*}
\bigl(\tfrac{1}{\gamma}\kors (\decreasing(w_P^{(\gamma)}\sqrt\gamma)\cap\sigma_\gamma\cap P')
+(\rho_P+1)\mu(P'\setminus S_\iota)\bigr)/\mu(P')\\
\ge
L(\rho_P,u^P_xu^P_y)-3(\rho_P+1)\iota.
\end{multline*}

By \cref{lm:greatrhoconstant} and the fact
that $L$ is uniformly continuous (\cref{lm:Lcontinuous}) on the bounded
set $\{(\rho(x,y),u_x(x,y)u_y(x,y))\,:\,(x,y)\in \bigcup_\iota S_\iota\}$,
we obtain
\begin{multline}\label{eq:Pprime}
\bigl(\tfrac{1}{\gamma}\kors (\decreasing(w_P^{(\gamma)}\sqrt\gamma)\cap\sigma_\gamma\cap P')
+(\rho_P+1)\mu(P'\setminus S_\iota)\bigr)/\mu(P')\\
> L(\rho(x_P,y_P),u^P_xu^P_y)-o_\iota(1),
\end{multline}
where $o_\iota(1)$ is independent of $P$.
But the inequality in \cref{eq:Pprime} holds also (even deterministically for any $w_P^{(\gamma)}$)
for any $P\in\calP_\iota$ that is not well behaved by
the fact that $\rho$, $u_x$ and $u_y$ are bounded on $\bigcup_\iota S_\iota$ and $u^P_x$ or $u^P_y$
is at most $\iota^3$ if $P$ is not well behaved together with the
fact that $L$ is continuous and $L(\rho,0)=0$ by \cref{lm:Lcontinuous}.

Since $w^{(\gamma)}_P\in\calU_{u(x_P,y_P)-\frac{r_P}{2},r_P}(P')$,
for any well-behaved $P\in\calP_\iota$ we have
\begin{equation}\label{eq:wclosetouxpyp}
\sup_{(x,y)\in P'}\,\abs{w_P^{(\gamma)}(x,y)-u(x_P,y_P)}\le r_P/2 = (1-2\iota)\iota c_P.
\end{equation}

Let $\tilde{\calP}_\iota\subseteq\calP_\iota$ be the set of well-behaved parallelograms
in $\calP_\iota$, and
let $w_\gamma\in \calU(\bigcup_{P\in\tilde{\calP}_\iota}P'\cup(\Omega\setminus\interior\cup\tilde{\calP}_\iota))$
be defined by $w_\gamma(x,y):=w_P^{(\gamma)}(x,y)$ if $(x,y)\in P'$ where $P\in\tilde{\calP}_\iota$ and
$w_\gamma(x,y):=u(x,y)$ if $(x,y)\in\Omega\setminus\interior\cup\tilde{\calP}_\iota$.
We claim that $w_\gamma$ is doubly increasing. To show this, it suffices to check
the inequality condition for any pair of points
$(x_1,y_1)$ and $(x_2,y_2)$ where
\begin{itemize}
\item
$(x_1,y_1)\le(x_2,y_2)$ or $(x_1,y_1)\ge(x_2,y_2)$, and
\item
$(x_1,y_1)$ belongs to $P'$ for some
$P\in\tilde{\calP}_\iota$ and $(x_2,y_2)$ lies
on the boundary of $P$.
\end{itemize}
\Cref{fig:wellbehavedparallelograms} shows an example.
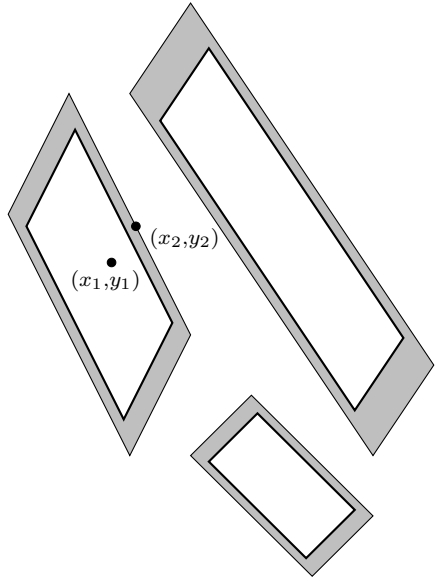
\begin{figure}
\begin{tikzpicture}[scale=0.8]
\draw[fill=lightgray] (0,0) -- (-2,4) -- (-1,6) -- (1,2) -- cycle;
\draw[thick,fill=white] (-0.1,0.6) -- (-1.7,3.8) -- (-0.9,5.4) -- (0.7,2.2) -- cycle;
\draw[fill=lightgray] (4,0) -- (0,6) -- (1,7.5) -- (5,1.5) -- cycle;
\draw[thick,fill=white] (3.7,0.75) -- (0.5,5.55) -- (1.3,6.75) -- (4.5,1.95) -- cycle;
\draw[fill=lightgray] (3,-2) -- (1,0) -- (2,1) -- (4,-1) -- cycle;
\draw[thick,fill=white] (2.9,-1.7) -- (1.3,-0.1) -- (2.1,0.7) -- (3.7,-0.9) -- cycle;
%\draw (-1,3.2) node [black,midway,right,xshift=0pt,yshift=0pt] {$(x_1,x_2)$};
\draw[fill] (-0.3,3.2) circle(2pt) (0.1,3.8) circle(2pt);
\node at (-0.4,2.9) {$\scriptstyle (x_1,y_1)$};
\node at (0.9,3.6) {$\scriptstyle (x_2,y_2)$};
\end{tikzpicture}
\caption{The situation in the proof of \cref{lm:notworsethanfunctional}. We see three disjoint well-behaved parallelograms and their slightly smaller primed counterparts. The shaded area
is initially excluded from the domain of $w_\gamma$.}
\label{fig:wellbehavedparallelograms}
\end{figure}
We have
\begin{equation}\label{eq:ubar}
s\cdot\bigl(\bar{u}_P(x_2,y_2)-u(x_P,y_P)\bigr)=\iota c_P,
\end{equation}
where $s$ is $+1$ if $(x_1,y_1)\le(x_2,y_2)$ and $-1$ if
$(x_1,y_1)\ge(x_2,y_2)$.
Also,
\begin{equation}\label{eq:wubar}
\abs{w_\gamma(x_1,y_1)-u(x_P,y_P)}\le(1-2\iota)\iota c_P
\end{equation}
by \cref{eq:wclosetouxpyp},
and
\begin{multline}\label{wubartwo}
\abs{w_\gamma(x_2,y_2)-\bar{u}_P(x_2,y_2)}=\abs{u(x_2,y_2)-\bar{u}_P(x_2,y_2)}\\
\le(\abs{x_2-x_P}+\abs{y_2-y_P})\iota^5
\le 2\iota^2 c_P,
\end{multline}
where the first inequality follows from \cref{lm:greatulinear}
and the last inequality follows from \cref{lm:uiparallelogrambound}
together with the fact that $\iota<1$.
Combining \cref{eq:ubar,eq:wubar,wubartwo} and the triangle inequality, we obtain
\begin{multline*}
s\cdot\bigl(w_\gamma(x_2,y_2)-w_\gamma(x_1,y_1)\bigr)
\ge\\
s\cdot\bigl(\bar{u}_P(x_2,y_2)-u(x_P,y_P)\bigr)
- \abs{w_\gamma(x_1,y_1)-u(x_P,y_P)}
- \abs{w_\gamma(x_2,y_2)-\bar{u}_P(x_2,y_2)}\\
\ge \iota c_P - (1-2\iota)\iota c_P - 2\iota^2 c_P = 0,
\end{multline*}
and we have showed that $w_\gamma$ is increasing.

For any $P\in\tilde\calP_\iota$, inside $P'$ we have
\begin{multline*}
\abs{w_P^{(\gamma)}(x,y)-u(x,y)}
\le\abs{w_P^{(\gamma)}(x,y)-u(x_P,y_P)} + \abs{u(x,y)-u(x_P,y_P)}\\
\le (1-2\iota)\iota c_P + (1+5\iota)\iota c_P,
\end{multline*}
where the first inequality is the triangle inequality and the second
inequality follows from \cref{eq:wclosetouxpyp,lm:greatuvarieslittle}. Since $c_P<1$ by \cref{lm:greatmisc}, it
follows that
\begin{equation}\label{eq:wclosetou}
\sup_{P\in\tilde{\calP}_\iota}\sup_{(x,y)\in P'}\abs{w_P^{(\gamma)}(x,y)-u(x,y)}\rightarrow0
\end{equation}
as $\iota\rightarrow0$. 

By \cref{lm:increasingextension}, $w_\gamma$ can be extended to a doubly increasing function on all of $\Omega$. Let us choose such an
extension and denote it by the same symbol $w_\gamma$.
Since $w_\gamma$
coincides with $u$ outside the interiors of the well-behaved
parallelograms, $w_\gamma\in\calU_{0,r}(\Omega)$.

We have
\begin{multline*}
\norm{w_\gamma-u}_\Omega = \sum_{P\in\tilde{\calP}_\iota}\norm{w_\gamma-u}_{P'}
+\sum_{P\in\tilde{\calP}_\iota}\norm{w_\gamma-u}_{P\setminus P'}\\
\le
\sum_{P\in\tilde{\calP}_\iota}\mu(P')\sup_{(x,y)\in P'}\abs{w_P^{(\gamma)}(x,y)-u(x,y)} + r\sum_{P\in\tilde{\calP}_\iota}\mu(P\setminus P')
=o_\iota(1)
\end{multline*}
by \cref{eq:wclosetou}.
Also, by \cref{eq:Pprime},
\begin{multline*}
\tfrac{1}{\gamma}\kors (\decreasing(w_\gamma\sqrt\gamma) \cap \sigma_\gamma) \ge
\sum_{P\in\calP_\iota}\tfrac{1}{\gamma}\kors(\decreasing(w_\gamma\sqrt\gamma)\cap\sigma_\gamma\cap P')
\\
\ge-(1+\sup_{S_\iota}\rho)\mu(\cup\calP_\iota\setminus S_\iota)
+\sum_{P\in\calP_\iota}\bigl(L(\rho(x_P,y_P),u^P_xu^P_y)-o_\iota(1)\bigr)\mu(P').
\end{multline*}
By \cref{lm:great}, part
\reflocal{lm:greatScoversparallelograms} and \reflocal{lm:greatLconstant}, this is greater than
\begin{multline*}
-o_\iota(1)+\sum_{P\in\calP_\iota}\left(\norm{L(\rho,u_xu_y)}_P\frac{\mu(P')}{\mu(P)} - o_\iota(1)\mu(P')\right)\\
\ge(1-2\iota)\norm{L(\rho,u_xu_y)}_{\cup\calP_\iota}-o_\iota(1)
\end{multline*}
which is greater than
$\norm{L(\rho,u_xu_y)}_\Omega - \eps - o_\iota(1)$, by \cref{lm:greatparallelogramscovereverything} and the fact that
$L(\rho,u_xu_y)\le\rho$ by \cref{lm:Lcontinuous}. Since
$\eps<\eps'$, we are done.
\end{proof}

Finally, we are ready to prove our first main theorem.
\main
\begin{proof}
Take any $\eps>0$ and $r\ge0$. Let $M_\eps$ denote the set of
$u\in\calU_{0,r}(\Omega)$ such that $\norm{u-u_{\rm max}}_\Omega<\eps$ for
some $u_{\rm max}\in\calU_{0,r}(\Omega)$ with
$\functional_\rho(u_{\rm max})=\Fmax(r)$. Let $M_\eps^c=\calU_{0,r}(\Omega)\setminus M_\eps$ denote the complementary set.
\begin{claim}\label{cl:supissmall}
$\sup_{u\in M_\eps^c} \functional_\rho(u)<\Fmax(r)$.
\end{claim}
Suppose there is a sequence $u_1,u_2,\dotsc \in M_\eps^c$ with
$\functional_\rho(u_i)\rightarrow \Fmax(r)$. By \cref{pr:compactness},
$M_\eps^c$ is compact, so
there is a convergent subsequence
$u_{i_1},u_{i_2},\dotsc\rightarrow u\in M_\eps^c$.
Since $\functional_\rho$ is upper semicontinuous (\cref{pr:semicontinuity}),
$\functional_\rho(u)\ge \limsup \functional_\rho(u_{i_j})=\Fmax(r)$, which is
a contradiction. This proves \cref{cl:supissmall}.

By the claim above we can choose an
$\eps'>0$ such that
\begin{equation}\label{eq:epsprim}
\sup_{u\in M_\eps^c}\functional_\rho(u) < \Fmax(r) - 3\eps'.
\end{equation}
Without loss of generality, we assume that $\eps'<\eps$.
Let $\{\sigma_\gamma\}_{\gamma>0}$ be Poisson point processes
on $\Omega$ with intensity functions $\gamma\rho$ such that $\tau_\gamma$
approaches $\sigma_\gamma$ as $\gamma\rightarrow\infty$.
\begin{claim}\label{cl:anywbad}
\Aas\ as $\gamma\rightarrow\infty$,
any $w\in M_\eps^c$ satisfies the inequality
\[
\kors(\decreasing(w\sqrt\gamma)\cap\sigma_\gamma)/\gamma < \Fmax(r)-2\eps'.
\]
\end{claim}
\begin{claim}\label{cl:anyIbad}
\Aas\ as $\gamma\rightarrow\infty$, any
$\lfloor r\sqrt{\gamma}\rfloor$-decreasing subset $P$ of $\tau_\gamma$ such that $\kappa_P/\sqrt{\gamma}\in M_\eps^c$ satisfies the inequality
$\kors P/\gamma < \Fmax(r)-\eps'$.
\end{claim}
\begin{claim}\label{cl:goodI}
\Aas\ as $\gamma\rightarrow\infty$,
there is a $\lfloor r\sqrt{\gamma}\rfloor$-decreasing subset $P$ of $\tau_\gamma$ such that $\kors P/\gamma > \Fmax(r)-\eps'$.
\end{claim}

For each $u\in\calU_{0,r}(\Omega)$, choose $\delta_u$
according to \cref{lm:notbetterthanfunctional}. The open balls
$\{B_{\delta_u}(u)\,:\,u\in M_\eps^c\}$ cover $M_\eps^c$, which is
compact by \cref{pr:compactness},
so there
is a finite subcover $\{B_{\delta_u}(u)\,:\,u\in A\subset M_\eps^c\}$,
$A$ finite. Now \cref{cl:anywbad} follows from \cref{lm:notbetterthanfunctional} (and the choice of
$\delta_u$) together with \cref{eq:epsprim}.

To show \cref{cl:anyIbad}, consider a
$\lfloor r\sqrt{\gamma}\rfloor$-decreasing subset $P$ of $\tau_\gamma$.
We have
$\kors P/\gamma\le
\kors(P\cap\sigma_\gamma)/\gamma + \kors(\tau_\gamma\setminus\sigma_\gamma)/\gamma$,
which is smaller than
$\frac1\gamma\kors(P\cap\sigma_\gamma) + \eps'$
\aas\ as $\gamma\rightarrow\infty$.
By \cref{lm:IandkappaIkappaPisP}, $\decreasing(\kappa_P)=P$,
so we obtain the inequality
$\kors P/\gamma<\frac1\gamma\kors(\decreasing(\kappa_P)\cap\sigma_\gamma)+\eps'$.
Now \cref{cl:anyIbad} follows from \cref{cl:anywbad} by
putting $w=\kappa_P/\sqrt{\gamma}$.

Take $u\in\calU_{0,r}(\Omega)$
such that $\functional_\rho(u)>\Fmax(r)-\tfrac13\eps'$.
By \cref{lm:notworsethanfunctional}, \aas\ as $\gamma\rightarrow\infty$,
there is a $w_\gamma\in\calU_{0,r}(\Omega)$ such that
$\kors(\decreasing(w_\gamma\sqrt\gamma)\cap\sigma_\gamma)/\gamma > \functional_\rho(u)-\frac13\eps'>\Fmax(r)-\tfrac23\eps'$
and hence
$\kors(\decreasing(w_\gamma\sqrt\gamma)\cap\tau_\gamma)/\gamma > \Fmax(r)-\eps'$
since 
$\kors(\sigma_\gamma\setminus\tau_\gamma)/\gamma$ tends to zero
in probability as $\gamma\rightarrow\infty$.
By \cref{lm:IandkappaIudecreasing}, $\decreasing(w_\gamma\sqrt\gamma)\cap\tau_\gamma$
is a $\lfloor r\sqrt{\gamma}\rfloor$-decreasing subset of $\tau_\gamma$,
and \cref{cl:goodI} follows.

Now \reflocal{th:mainlimitsurface} follows from \cref{cl:anyIbad,cl:goodI}.

\begin{claim}\label{cl:anywbadagain}
\Aas\ as $\gamma\rightarrow\infty$,
any $w\in\calU_{0,r}(\Omega)$ satisfies the inequality
$\kors(\decreasing(w\sqrt\gamma)\cap\sigma_\gamma)/\gamma < \Fmax(r)+\eps'$.
\end{claim}
\begin{claim}\label{cl:anyIbadagain}
\Aas\ as $\gamma\rightarrow\infty$, any
$\lfloor r\sqrt{\gamma}\rfloor$-decreasing subset $P$ of $\tau_\gamma$ satisfies the inequality
$\kors P/\gamma < \Fmax(r)+2\eps'$.
\end{claim}

There is a finite subcover
$\{B_{\delta_u}(u)\,:\,u\in A'\subset \calU_{0,r}(\Omega)\}$ of $\calU_{0,r}(\Omega)$, so
\cref{cl:anywbadagain} follows from \cref{lm:notbetterthanfunctional}.

By \cref{pr:decreasingincreasingrelation}, for any 
$\lfloor r\sqrt{\gamma}\rfloor$-decreasing set $P$, we have
$\kappa_P/\sqrt{\gamma}\in\calU_{0,r}(\Omega)$. With this in mind,
\cref{cl:anyIbadagain} follows from \cref{cl:anywbadagain}
the same way as \cref{cl:anyIbad} followed from \cref{cl:anywbad}.

Finally, \reflocal{th:mainlimittoFmax} follows from \cref{cl:goodI,cl:anyIbadagain}
since $\eps>0$ was chosen arbitrarily and $\eps'<\eps$. 
\end{proof}

\limitshape
\begin{proof}
By \cref{pr:curtisgreene}, the maximal size $\Lambda^{(\gamma)}$ of a
$\lfloor r\sqrt\gamma\rfloor$-decreasing subset of $\tau_\gamma$ is
$\Lambda^{(\gamma)}=\sum_{i=1}^{\lfloor r\sqrt\gamma\rfloor}\lambda^{(\gamma)}_i$,
and by \cref{th:mainlimittoFmax}, $\Lambda^{(\gamma)}/\gamma\rightarrow \Fmax(r)$ in
probability as $\gamma\rightarrow\infty$. Now the corollary follows from \cref{lm:Lambdatolambda}
with $a(\lambda)=\sqrt\gamma$, $b(\gamma)=1/\gamma$ and $G=\Fmax$.
\end{proof}

\rhombuslimitshapeexists
\begin{proof}
By \cref{pr:parallelogrammaximizer}, for any $r\ge0$, the maximum value of $\functional_1$ on $\calU_r(\Omega)$ is $\Fmax(r)=\Phi(r)$, so
$\Fmax'(r)=\Phi'(r)$ whenever the derivative exists.
Now the theorem follows from \cref{cor:limitshape}.
\end{proof}

\section{Concavity and existence of maximizers}\label{sec:mainFproof}
Before we are ready to prove \cref{th:mainF}, we need another property of
$\Phi$ and a technical lemma.

The following proposition is stated in terms of the constant $\Gamma$ from
\cref{th:hammersley}. Of course we know that $\Gamma=2$ from the result of Vershik and Kerov \cite{VershikKerov}, but we do not need that in the proof of \cref{th:mainF}. As discussed in \cref{sec:future}, a
proof of \cref{con:triangularlimitshape} would yield a conceptually new
proof of the Logan--Shepp--Vershik--Kerov limit shape, so to avoid a
circular dependence we take care not to rely on that result.
\begin{prop}\label{pr:phioneisone}
$\Phi(r)=1$ if $r\ge\Gamma/\sqrt2$, where $\Gamma$ is the constant from
\cref{th:hammersley}.
\end{prop}
\begin{proof}
By \cref{pr:phiproperties}, $\Phi$ is continuous at
$r=\Gamma/\sqrt2$, so it suffices to show that $\Phi(r)=1$ for any $r>\Gamma/\sqrt2$.
For any $\beta>0$, consider the density domain $(\Omega_\beta, 1)$
where $\Omega_\beta$ is the open
rectangle $\abs{x+y}<1$, $\abs{x-y}<\beta$, and
for any $\gamma>0$ and $\beta>0$, let $\sigma_{\gamma,\beta}$ be a Poisson point
process on $\Omega_\beta$ with homogeneous intensity $\gamma$.
Since $\Omega_\beta\subset(-\tfrac{1+\beta}2,\tfrac{1+\beta}2)^2$, by \cref{th:hammersley}, for any $\eps>0$, the size of the largest increasing subset
of $\sigma_{\gamma,\beta}$ is smaller than $(\Gamma+\eps)(1+\beta)\sqrt{\gamma}$
\aas\ as $\gamma\rightarrow\infty$. Thus, by \cref{pr:decreasingincreasingrelation}, the maximum
size of a $\lfloor(\Gamma+\eps)(1+\beta)\sqrt{\gamma}\rfloor$-decreasing subset
of $\sigma_{\gamma,\beta}$ is $\kors\sigma_{\gamma,\beta}$ \aas\ By the law of large numbers, we have
$\kors\sigma_{\gamma,\beta}/\gamma\rightarrow2\beta$ in probability as $\gamma\rightarrow\infty$, so, by \cref{th:mainlimittoFmax},
$\Fmax((\Gamma+\eps)(1+\beta))= 2\beta$ for any $\eps>0$ and $\beta>0$. On the other hand, by \cref{pr:parallelogrammaximizer},
$\Fmax((\Gamma+\eps)(1+\beta))=2\beta\Phi((\Gamma+\eps)(1+\beta)/\sqrt2)$. It follows that $\Phi((\Gamma+\eps)(1+\beta)/\sqrt2)=1$.
\end{proof}

\begin{lemma}\label{lm:concavitycomposition}
Let $A$ be a convex subset of a vector space and let $B\subseteq\bbR$.
Let $\psi:A\rightarrow B$ be a concave function and let
$\phi:B\rightarrow\bbR$ be an increasing concave function.
Then $\phi\circ\psi:A\rightarrow\bbR$ is concave.
\end{lemma}
\begin{proof}
Since $\psi$ is concave, for any $t\in[0,1]$ we have
\[
\psi((1-t)x+ty)\ge (1-t)\psi(x)+t\psi(y)
\]
and thus, since $\phi$ is increasing,
\[
\phi(\psi((1-t)x+ty))\ge\phi((1-t)\psi(x)+t\psi(y)),
\]
which, since $\phi$ is concave, is at least
\[
(1-t)\phi(\psi(x))+t\phi(\psi(y)).
\]
Thus, $\psi\circ\phi$ is concave.
\end{proof}

Finally, we have all we need to prove our second main theorem.
\mainF
\begin{proof}
\reflocal{th:mainFmaximizer}
By \cref{pr:semicontinuity},
$\functional_\rho$ is upper semicontinuous in $L^1(\Omega)$,
and by \cref{pr:compactness}, $\calU_{0,r}(\Omega)$ is compact, so $\functional_\rho$ attains its maximum there.

\smallskip
\reflocal{th:mainFrhoconcave}
For any point $(x,y)\in\Omega$, let $\calU_{(x,y)}$ be the
set of $u$ in $\calU(\Omega)$ such that $u_x(x,y)$ and $u_y(x,y)$ exist,
and define the function
$\chi_{(x,y)}:\calU_{(x,y)}\rightarrow\bbR^2$ by letting
$\chi_{(x,y)}(u)=(u_x(x,y),u_y(x,y))$. Clearly, $\chi_{(x,y)}$ is
a linear function. 

Let $\supp\rho=\{(x,y)\in\Omega\,:\,\rho(x,y)>0\}$.
For any $(x,y)\in\supp\rho$, define $\psi\,:\,\bbR_{\ge0}^2\rightarrow\bbR_{\ge0}$ by $\psi_{(x,y)}(r,s)=\sqrt{2rs/\rho(x,y)}$.
We can easily check that $\psi_{(x,y)}$ is concave; its Hessian determinant
is negative semidefinite on $\bbR_{>0}^2$. Also, $\Phi$ is concave and
increasing by \cref{pr:phiproperties}, and hence $\Phi\circ\psi_{(x,y)}$ is concave by \cref{lm:concavitycomposition}.

Since $\chi_{(x,y)}$ is linear, $\Phi\circ\psi_{(x,y)}\circ\chi_{(x,y)}$ is
concave. By \cref{lm:monotonicdifferentiable}, for any $u\in\calU(\Omega)$,
$\chi_{(x,y)}$ is defined at $u$ for almost every $(x,y)\in\Omega$,
so the integral
\[
\int_{\supp\rho}\rho(x,y)(\Phi\circ\psi_{(x,y)}\circ\chi_{(x,y)})(u)\,d\mu
\]
is well defined and concave as a function of $u$.

\smallskip
\reflocal{th:mainFmaxcontinuousincreasingconcave}
That $\Fmax$ is increasing follows from the fact that $\calU_{r_1}(\Omega)\subseteq\calU_{r_2}(\Omega)$ whenever $r_1\le r_2$.
To show that $\Fmax$ is concave, we must check that
$\Fmax((1-t)r_1+tr_2)\ge (1-t)\Fmax(r_1)+t\Fmax(r_2)$
for any $r_1,r_2\ge0$ and any $t\in(0,1)$. By \reflocal{th:mainFmaximizer}, there are $u^{(1)}\in\calU_{0,r_1}$ and
$u^{(2)}\in\calU_{0,r_2}$ such that $\functional_\rho(u^{(1)})=\Fmax(r_1)$ and
$\functional_\rho(u^{(2)})=\Fmax(r_2)$. Let $u:=(1-t)u^{(1)}+tu^{(2)}$.
Clearly, $u\in\calU_{0,(1-t)r_1+tr_2}$.
By \reflocal{th:mainFrhoconcave}, $\functional_\rho$ is concave, so
\begin{multline*}
\Fmax((1-t)r_1+tr_2)\ge\functional_\rho(u)=\functional_\rho((1-t)u^{(1)}+tu^{(2)})
\\
\ge(1-t)\functional_\rho(u^{(1)}) + t\functional_\rho(u^{(2)})
=(1-t)\Fmax(r_1) + t\Fmax(r_2).
\end{multline*}
This shows that $\Fmax$ is concave.

A concave function is automatically continuous on any open interval, so to
show that $\Fmax$ is continuous we need only to show that it is continuous at zero.
Choose any $\eps>0$.
Since $\int_\Omega\rho\,d\mu$ is finite,
there is a measurable subset $S$ of $\Omega$ and positive constants $\rho_0,\rho_1$
such that $\int_{\Omega\setminus S} \rho\,d\mu<\eps$ and $\rho_0<\rho<\rho_1$ on $S$.
There is also a $c>0$ such that
$\int_{\Omega\setminus[-c,c]^2}\rho\,d\mu<\eps$.
Choose any $r>0$ and any $u\in\calU_r(\Omega)$. By \cref{lm:increasingextension},
there is a $w\in\calU_r(\bbR^2)$ that coincides with $u$ on $\Omega$,
and by Tonelli's theorem,
\[
\int_{\Omega\cap[-c,c]^2}u_x\,d\mu
\le\int_{[-c,c]^2}w_x\,d\mu
=\int_{-c}^c\left(\int_{-c}^c w_x(x,y)\,dx\right)dy.
\]
Now, by Lebesgue's theorem for increasing functions in one dimension,
\[
\int_{-c}^c w_x(x,y)\,dx\le w(c,y)-w(-c,y),
\]
which is at most $r$ since $w\in\calU_r(\bbR^2)$.
Thus, $\int_{\Omega\cap[-c,c]^2}u_x\,d\mu\le \int_{-c}^cr\,dy=2cr$
and, analogously, $\int_{\Omega\cap[-c,c]^2}u_y\,d\mu\le 2cr$.
Choose any $\delta>0$ and
let $T$ be the subset of $S\cap[-c,c]^2$ where $\sqrt{2u_xu_y/\rho}\ge\delta$. Since $\rho>\rho_0$ on $S$, on $T$ we have
$u_xu_y\ge \delta^2\rho/2>\delta^2\rho_0/2$, so $T\subseteq T_1\cup T_2$,
where $T_1$ is the subset of $S\cap[-c,c]^2$ where $u_x\ge \delta\sqrt{\rho_0/2}$ and $T_2$ is the subset of $S\cap[-c,c]^2$ where $u_y\ge \delta\sqrt{\rho_0/2}$. By Markov's inequality,
\[
\mu(T_1)\le\frac1{\delta\sqrt{\rho_0/2}}\int_{S\cap[-c,c]^2} u_x\,d\mu
\le \frac{2cr}{\delta\sqrt{\rho_0/2}},
\]
and analogously for $T_2$, so $\mu(T)\le\mu(T_1)+\mu(T_2)\le 4cr/\delta\sqrt{\rho_0/2}$.
Combining all above, bearing in mind that $L(\rho,u_xu_y)\le\rho$
by \cref{lm:Lcontinuous} and that $\Phi$ is increasing by \cref{pr:phiproperties}, we obtain
\begin{align*}
\functional_\rho(u)&= \int_\Omega L(\rho,u_xu_y)\,d\mu \\
&\le \int\limits_{\Omega\setminus S} \rho\,d\mu\ +
\!\!\!\int\limits_{\Omega\setminus [-c,c]^2}\!\!\!\rho\,d\mu\ +
\int\limits_T \rho\,d\mu
\ +\!\!\!\!\!\int\limits_{(S\cap[-c,c]^2)\setminus T}\!\!\!\!\!\rho\,\Phi(\lagomsqrt{2u_xu_y/\rho})\,d\mu \\
&<2\eps + \mu(T)\rho_1 +
\mu((S\cap[-c,c]^2)\setminus T)\rho_1\Phi(\delta) \\
&\le2\eps + \frac{4cr\rho_1}{\delta\sqrt{\rho_0/2}} +
\mu(\Omega)\rho_1\Phi(\delta).
\end{align*}
Since $\eps$, $r$ and $\delta$ were chosen freely above,
we have shown that
\[
\Fmax(r)\le2\eps + \frac{4cr\rho_1}{\delta\sqrt{\rho_0/2}} +
\mu(\Omega)\rho_1\Phi(\delta)
\]
for any $\eps,\delta,r>0$. Here, $c$, $\rho_0$ and $\rho_1$ depend on $\eps$ but not on $r$ or $\delta$.
By \cref{pr:phiproperties}, $\Phi$ is continuous, so for
any $\eps>0$ we can choose $\delta>0$ such that
$\mu(\Omega)\rho_1\Phi(\delta)<\eps$. After that, we can choose $r>0$
such that $4cr\rho_1/(\delta\sqrt{\rho_0/2})<\eps$. It follows that
$\Fmax(r)<4\eps$ and since $\Fmax$ is an increasing function, we conclude that
$\lim_{r\rightarrow0+}\Fmax(r)=0$. Hence, $\Fmax$ is continuous.

To show that $\Fmax(r)\rightarrow\norm{\rho}_\Omega$ as
$r\rightarrow\infty$, let us again choose any $\eps>0$ and take
$S$, $\rho_0$, $\rho_1$ and $c$ as above. Choose any $r\ge0$ and
let $u\in\calU_r(S\cap[-c,c]^2)$ be defined by $u(x,y):=r(x+y)\,/\,4c$.
By \cref{lm:increasingextension}, $u$ can be extended to $\Omega$ without
increasing the diameter of its image, so from
now on we consider $u$ to be an element of $\calU_r(\Omega)$.
Since $u_x=u_y=r\,/\,4c$ almost everywhere on $S\cap[-c,c]^2$,
again bearing in mind that $L(\rho,u_xu_y)\le\rho$ and $\Phi$ is increasing we obtain
\begin{align*}
0\le\norm{\rho}_\Omega - \functional_\rho(u)&=
\int_\Omega \bigl(\rho - L(\rho,u_xu_y)\bigr)\,d\mu \\
&\le \int\limits_{\Omega\setminus S} \rho\,d\mu\ +
\!\!\!\int\limits_{\Omega\setminus [-c,c]^2}\!\!\!\rho\,d\mu
\ +\!\!\!\!\!\int\limits_{S\cap[-c,c]^2}\!\!\!\!\!\rho\,\bigl(1-\Phi(\lagomsqrt{2u_xu_y/\rho})\bigr)\,d\mu \\
&<2\eps +
\mu(S\cap[-c,c]^2)\rho_1\bigl(1-\Phi(\sqrt{2/\rho_1}\cdot r/4c)\bigr) \\
&\le2\eps + 4\mu(\Omega)\rho_1\bigl(1-\Phi(\sqrt{2/\rho_1}\cdot r/4c)\bigr).
\end{align*}
Since $\eps$ and $r$ were chosen freely above, we have shown that,
for any $\eps>0$ and $r\ge0$, there is a $u\in\calU_r(\Omega)$ such that
\[
0\le\norm{\rho}_\Omega - \functional_\rho(u) \le 2\eps + 4\mu(\Omega)\rho_1\bigl(1-\Phi(\sqrt{2/\rho_1}\cdot r/4c)\bigr).
\]
Here, $c$ and $\rho_1$ depend on $\eps$ but not on $r$.
By \cref{pr:phioneisone}, the last term tends to zero
as $r$ tends to infinity, and we conclude that
$\norm{\rho}_\Omega - \Fmax(r)\rightarrow 0$ as $r\rightarrow\infty$.
\end{proof}

\section{Essentially unique maximizers}\label{sec:uniquemaximizers}
In this section we show that, under reasonable assumptions, the maximizer
of the functional $\functional_\rho$ is essentially unique.

\begin{defi}
Let $\psi$ be a real-valued function on a convex subset $A$ of
a vector space and let
$C$ be a subset of $A$. Then, $\psi$ is said to be
\emph{strictly concave from $C$} if, for any $x\in C$
and $y\in A$ with
$x\ne y$, and any $t\in(0,1)$, we have
$\psi((1-t)x+ty)>(1-t)\psi(x)+t\psi(y)$.
\end{defi}

\begin{defi}
Let $\phi$ be a function from some $B\subseteq\bbR$ to $\bbR$,
and let $C$ be a subset of $B$.
We say that $\phi$ is \emph{strictly increasing from $C$}
if $\phi(x)<\phi(y)$ for any $C\ni x<y\in B$.
\end{defi}

\begin{lemma}\label{lm:advancedconcavitycomposition}
Let $A$ be a convex subset of a vector space and let $B\subseteq\bbR$.
Let $\psi:A\rightarrow B$ be a concave function and let
$\phi:B\rightarrow\bbR$ be an increasing function.
Suppose $\psi$ is strictly concave on the line through any
two distinct points $x\ne y$ with $\psi(x)=\psi(y)$.
Suppose also that $\phi$ is
strictly concave from some subset $C$ of $B$ and
strictly increasing from $C$.
Then $\phi\circ\psi:A\rightarrow\bbR$
is strictly concave from $\psi^{-1}(C)$.
\end{lemma}
\begin{proof}
Take any $x\in\psi^{-1}(C)$ and $y\in A$ with $x\ne y$, and take
any $t\in(0,1)$.
We must show that
\[
(\phi\circ\psi)((1-t)x+ty)>(1-t)(\phi\circ\psi)(x)+t(\phi\circ\psi)(y).
\]

Since $\psi$ is concave, we have
\begin{equation}\label{eq:psiconcave}
\psi((1-t)x+ty)\ge(1-t)\psi(x)+t\psi(y)
\end{equation}
and thus, since $\phi$ is increasing,
\begin{equation}\label{eq:phiincreasing}
\phi(\psi((1-t)x+ty))\ge\phi((1-t)\psi(x)+t\psi(y)).
\end{equation}
Since $\phi$ is strictly concave from $C$, we have
\[
\phi((1-t)\psi(x)+t\psi(y))\ge(1-t)\phi(\psi(x))+t\phi(\psi(y)).
\]
with equality only if $\psi(x)=\psi(y)$.
If $\psi(x)=\psi(y)$, by the assumptions in the lemma,
$\psi$ is strictly concave on the line through $x$ and $y$
and hence the inequality in \cref{eq:psiconcave} holds strictly.
Since $\phi$ is strictly increasing from $C$, this implies that
the inequality in \cref{eq:phiincreasing} holds strictly too.

Thus, $\psi\circ\phi$ is strictly concave from $\psi^{-1}(C)$.
\end{proof}

\begin{prop}\label{pr:uniquemaximizer}
Suppose $\Phi$ is strictly concave on $[0,\sqrt2]$.
Then, if $u^{(1)}$ and $u^{(2)}$ are two maximizers of the operator $\functional_\rho$
in $\calU_r(\Omega)$, the two sets
\begin{align*}
\{(x,y)\ :\ 0<u^{(1)}_x(x,y)u^{(1)}_y(x,y)<\rho(x,y)\} & \text{\ and}\\
\{(x,y)\ :\ 0<u^{(2)}_x(x,y)u^{(2)}_y(x,y)<\rho(x,y)\} & \mbox{}
\end{align*}
are almost equal,
and the partial derivatives of
$u^{(1)}$ and $u^{(2)}$ agree almost everywhere on that set.
\end{prop}
\begin{proof}
For $i=1,2$, let $D_i\subseteq\Omega$ be the set of points $(x,y)$ where
$0<u^{(i)}_x u^{(i)}_y<\rho$. Let $A:=\bbR_{>0}^2$, $B:=\bbR_{\ge0}$ and
$C:=[0,\sqrt2)$. For any $(x,y)\in\supp\rho$, let $\chi_{(x,y)}$ and $\psi_{(x,y)}$ be defined
as in the proof of \cref{th:mainFrhoconcave}. It follows from \cref{pr:phiproperties} together with our assumption that $\Phi$ is strictly concave on $[0,\sqrt2]$
that the assumptions of \cref{lm:advancedconcavitycomposition} are
satisfied for $\psi_{(x,y)}$, $\Phi$, $A$, $B$ and $C$, so
$\Phi\circ\psi_{(x,y)}$ is strictly concave from
$\psi_{(x,y)}^{-1}(C)=\{(r,s)\in\bbR_{>0}^2\,:\,rs<\rho(x,y)\}$.

Let $D$ be the set of
points in $D_1\cup D_2$ where $(u^{(1)}_x,u^{(1)}_y)\ne (u^{(2)}_x,u^{(2)}_y)$.
We claim that $\mu(D)=0$.

Let $w:=(u^{(1)}+u^{(2)})/2$. If $(x,y)\in D$, the points
$p:=\chi_{(x,y)}(u^{(1)})$ and $q:=\chi_{(x,y)}(u^{(2)})$ are distinct
and at least one of them belongs to $\psi_{(x,y)}^{-1}(C)$.
Hence, with $L_{(x,y)}$ as a shorthand for $\rho(x,y)\cdot(\Phi\circ\psi_{(x,y)}\circ\chi_{(x,y)})$, we obtain
\begin{multline*}
L_{(x,y)}(w)
=\rho(x,y)(\Phi\circ\psi_{(x,y)})((p+q)/2)\\
>\rho(x,y)[(\Phi\circ\psi_{(x,y)})(p)+(\Phi\circ\psi_{(x,y)})(q)]/2
=[L_{(x,y)}(u^{(1)})
+L_{(x,y)}(u^{(2)})]/2.
\end{multline*}
Suppose $\mu(D)>0$. Then,
\[
\int_{D}L_{(x,y)}(w)\,d\mu
>\frac{1}{2}\left(\int_{D}L_{(x,y)}(u^{(1)})\,d\mu
+\int_{D}L_{(x,y)}(u^{(2)})\,d\mu
\right).
\]
Also, by \cref{th:mainFrhoconcave},
\[
\int\limits_{(\supp\rho)\setminus D}\!\!\!\!\!L_{(x,y)}(w)\,d\mu
\ge\frac{1}{2}\left(\,\int\limits_{(\supp\rho)\setminus D}\!\!\!\!\!L_{(x,y)}(u^{(1)})\,d\mu
\ +\!\!\!\!\!\int\limits_{(\supp\rho)\setminus D}\!\!\!\!\!L_{(x,y)}(u^{(2)})\,d\mu
\right),
\]
so it follows that
\[
\int_{\supp\rho}L_{(x,y)}(w)\,d\mu
>\frac{1}{2}\left(\int_{\supp\rho}L_{(x,y)}(u^{(1)})\,d\mu
+\int_{\supp\rho}L_{(x,y)}(u^{(2)})\,d\mu
\right).
\]
This means that $\functional_\rho(w)>[\functional_\rho(u^{(1)})+\functional_\rho(u^{(2)})]/2$,
which is impossible since $u^{(1)}$ and $u^{(2)}$ both are maximizers
of $\functional_\rho$. We conclude that $\mu(D)=0$.
\end{proof}

Recall the definition of $\calV$ from \cref{sec:generalpde}.
\begin{lemma}\label{lm:uniquefromderivativesgeneral}
Let $\Omega$ be an open subset of $\bbR^2$, and let $u\in\calU_{-c,2c}(\Omega)$ and $v\in\calV_{-d,2d}(\Omega)$ for some $c,d>0$.
Suppose that $u$ and $v$ are everywhere differentiable with nonzero partial derivatives and that the image of $\Omega$ under the map
$\varphi_u\,:\,(x,y)\mapsto(u(x,y),v(x,y))$ is $(-c,c)\times(-d,d)$.
Then, for any
$w\in\calU_{-c,2c}(\Omega)$ whose partial derivatives coincide with those of
$u$ almost everywhere, it holds that $w=u$ everywhere.
\end{lemma}
\begin{proof}
Take any $w\in\calU_{-c,2c}(\Omega)$ whose partial derivatives coincide with those of $u$ almost everywhere.
Let $(x_0,y_0)$ be any point in $\Omega$ and let $u_0:=u(x_0,y_0)$,
$v_0:=v(x_0,y_0)$ and $w_0=w(x_0,y_0)$. We need to show that $w_0=u_0$, but
by symmetry%
\footnote{To be precise: Suppose $\Omega$, $u$, $v$ and $w$ satisfy the assumptions in
the lemma, and define $\Omega'$, $u'$, $v'$ and $w'$ by
$\Omega':=-\Omega=\{(-x,-y)\,:\,(x,y)\in\Omega\}$,
$u'(x,y)=-u(-x,-y)$, $v'(x,y)=-v(-x,-y)$ and $w'(x,y)=-w(-x,-y)$.
Then, $\Omega'$, $u'$, $v'$ and $w'$ satisfy the assumptions in the lemma
too, and $w(x,y)\le u(x,y)$ if and only if $w'(x,y)\ge u'(x,y)$.}
it is enough to show that $w_0\ge u_0$, so this will be our goal.

Since $u$ has positive partial derivatives,
for any sufficiently small $\delta>0$ there is an $x_1<x_0$ and an $y_1<y_0$
such that $u(x_1,y_0)=u(x_0,y_1)=u_0-\delta$. Since $v_x<0$ and $v_y>0$ everywhere, we have $v_-:=v(x_0,y_1)<v(x_0,y_0)<v(x_1,y_0)=:v_+$.
Let $S:=(-c,u_0-\delta)\times(v_-,v_+)$ and $T:=\varphi_u^{-1}(S)$.
We claim that every point in $T$ is south-west
of $(x_0,y_0)$. To see this,
first note that $v\le v_-$ for any point south-east
of $(x_0,y_1)$ and $v\ge v_+$ for any point north-west of $(x_1,y_0)$. Then
note that $u\ge u_0-\delta$ for any point north-east of $(x_0,y_1)$ or
$(x_1,y_0)$. Thus, every point in $T$ is south of $(x_1,y_0)$ and west of $(x_0,y_1)$, and we have proved the claim.
\Cref{fig:uandv} illustrates the situation.
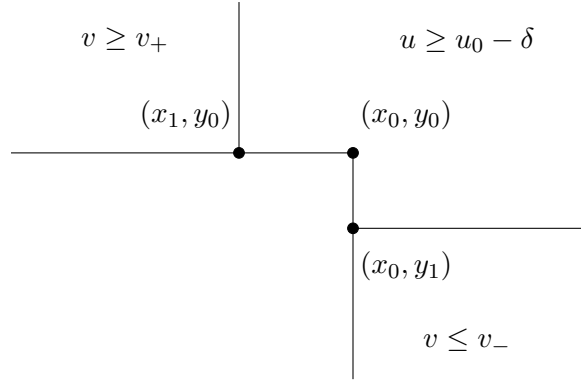
\begin{figure}
\begin{tikzpicture}
\draw[fill] (0,1) circle(2pt) (1.5,0) circle(2pt) (1.5,1) circle(2pt);
\draw (-0.7,1.5) node {$(x_1,y_0)$} (2.2,-0.5) node {$(x_0,y_1)$} (2.2,1.5) node {$(x_0,y_0)$};
\draw (-3,1) -- (0,1) -- (0,3);
\draw (1.5,-2) -- (1.5,0) -- (4.5,0);
\draw (0,1) -- (1.5,1) -- (1.5,0);
\draw (-1.5,2.5) node {$v\ge v_+$} (3,-1.5) node {$v\le v_-$} (3,2.5) node {$u\ge u_0-\delta$};
\end{tikzpicture}
\caption{The situation in the proof of \cref{lm:uniquefromderivativesgeneral}}
\label{fig:uandv}
\end{figure}

Let $J:=u_xv_y-u_yv_x$
be the Jacobian determinant of $\varphi_u$. By \cref{th:changeofvariables},
\[
\int_T J d\mu\ge \mu(\varphi_u(T))=\mu(S),
\]
where the last equality follows from the surjectivity of $\varphi_u$ onto
the set $(-c,c)\times(-d,d)$.
Let $T'$ be the subset of $T$ where
$w$ is differentiable and
the partial derivatives of $w$ and $u$ coincide, and
define the function $\varphi_w\,:\,T'\rightarrow[-c,c]\times[-d,d]$ by $\varphi_w(x,y)=(w(x,y),v(x,y))$. Then $\varphi_w$ is injective by \cref{lm:injectivity}, and by \cref{th:changeofvariables},
\[
\int_{T'} J d\mu = \mu(\varphi_w(T')).
\]
Since $\mu(T\setminus T')=0$, the two integrals above are equal, so
\begin{equation}\label{eq:mumu}
\mu(\varphi_w(T'))\ge\mu(S).
\end{equation}
Recall that any point in $T$ is south-west of $(x_0,y_0)$, so
$w(x,y)\le w_0$ for any point $(x,y)\in T$. It follows that
\[
\varphi_w(T')\subseteq[-c,w_0]\times(v_-,v_+),
\]
so $\mu(\varphi_w(T'))\le (w_0+c)(v_+-v_-)$. On the other hand, $\mu(S)=(u_0-\delta+c)(v_+-v_-)$, so \cref{eq:mumu} yields that $w_0\ge u_0-\delta$.
Since this holds for arbitrarily small positive $\delta$, we conclude that
$w_0\ge u_0$.
\end{proof}

\begin{lemma}\label{lm:uniquefromderivatives}
Let $\Omega$ be the open rhombus $a\abs{x}+b\abs{y}<1$
and let $c$ be a positive number. Let $u\in\calU_{-c,2c}$ and suppose that
$u_x=ca$ and $u_y=cb$ almost everywhere in $\Omega$. Then $u(x,y)=c(ax+by)$
everywhere on $\Omega$.
\end{lemma}
\begin{proof}
This follows from \cref{lm:uniquefromderivativesgeneral} with
$d=1$ and $v(x,y)=by-ax$.
\end{proof}

Our final result in this section relies on the fact due to Vershik and Kerov \cite{VershikKerov} that the constant $\Gamma$ in \cref{th:hammersley} equals 2.
\begin{prop}
Suppose $\Phi$ is strictly concave on $[0,\sqrt2]$. Let
$\Omega$ be the open rhombus $a\abs{x}+b\abs{y}<1$ and $\rho>0$
be constant, and let
$c$ be any positive number smaller than or equal to $\sqrt{\rho/ab}$.
Then, in $\calU_{-c,2c}(\Omega)$
the functional $\functional_\rho$ is uniquely maximized by the function
$u_{\rm linear}(x,y)=c(ax+by)$.
\end{prop}
\begin{proof}
By \cref{pr:parallelogrammaximizer}, $u_{\rm linear}$ is a maximizer of
$\functional_\rho$. If $c<\sqrt{\rho/ab}$, the uniqueness follows from
\cref{pr:uniquemaximizer} together with \cref{lm:uniquefromderivatives}.
Suppose $c=\sqrt{\rho/ab}$
and suppose there is another maximizer $u\in\calU_{-c,2c}$.
By scale invariance we can assume without loss of generality that
$a=b=\rho=1$.
We have $\functional_1(u_{\rm linear})=\mu(\Omega)\Phi(\sqrt2)$ which equals $\mu(\Omega)$
by \cref{pr:phioneisone} together with the result from \cite{VershikKerov} that $\Gamma=2$, so we must have $\Phi(\sqrt{2u_xu_y})=1$ and hence
$u_xu_y\ge1$ almost everywhere. It follows that
\begin{equation}\label{eq:integrallarge}
\int_\Omega \sqrt{u_xu_y}\,d\mu\ge 2.
\end{equation}
On the other hand, by the inequality of the geometric and
arithmetic mean,
\begin{equation}\label{eq:integralsmallone}
\int_\Omega \sqrt{u_xu_y}\,d\mu\le\frac12\int_\Omega (u_x+u_y)d\mu
\end{equation}
with equality if and only if $u_x=u_y$ almost everywhere in $\Omega$.
Define a map $\zeta\,:\,\bbR^2\rightarrow\bbR^2$ by $\zeta(\alpha,\beta)=(\frac{\alpha+\beta}{2},\frac{\alpha-\beta}{2})$.
Then, $\Omega=\zeta((-1,1)^2)$ and
\begin{multline}\label{eq:integralsmalltwo}
\frac12\int_\Omega (u_x+u_y)\,d\mu
=\frac14\int_{(-1,1)^2} \bigl((u_x+u_y)\circ\zeta\bigr)\,d\mu
=\frac12\int_{(-1,1)^2} \tfrac{\partial}{\partial\alpha}(u\circ\zeta)\,d\mu
\\
=\{\text{Tonelli's theorem}\}
=\frac12\int_{-1}^{1}\left(\int_{-1}^{1} \tfrac{\partial}{\partial\alpha}(u\circ\zeta)\,d\alpha\right)\,d\beta
\\
\le\{\text{since $u\circ\zeta$ is increasing in the first variable}\}
\le
\frac12\int_{-1}^{1}2c\,d\beta=2c=2.
\end{multline}
%By the change of variables $x=\alpha+\beta$, $y=\alpha-\beta$,
%we get
%\begin{multline}\label{eq:integralsmalltwo}
%\frac12\int_\Omega (w_x+w_y)d\mu=\iint_\Omega w_\alpha\,d\alpha\,d\beta
%\\=\int_{-1/2}^{1/2}\left(\int_{-1/2}^{1/2} w_\alpha\,d\alpha\right)\,d\beta
%\le
%\int_{-1/2}^{1/2}2c\,d\beta=2c=2.
%\end{multline}
Combining \cref{eq:integrallarge,eq:integralsmallone,eq:integralsmalltwo},
we see that $u_x=u_y=1$ must hold almost everywhere in $\Omega$.
By \cref{lm:uniquefromderivatives}, $u=u_{\rm linear}$ almost everywhere.
\end{proof}

\section{The uniform case}\label{sec:uniform}
In this section we suppose \cref{con:triangularlimitshape}
holds true and explore the consequences for the case of uniformly random permutations. The limit surfaces turn out to have a surprisingly
simple parameterization in terms of trigonometric functions,
and we are able to recover the Logan--Shepp--$\!$Vershik--Kerov
limit-shape result mentioned in \cref{sec:background} and depicted in
\cref{fig:famouslimitshape}. Level plots of some limit surfaces are shown
in \cref{fig:limitsurfaces}.
\begin{figure}%
\centerfloat
\subfloat[$\alpha=\pi/3$]{%
\quad%
\begin{tikzpicture}[scale=4.8]
  \draw[->] (-0.6, 0) -- (0.6, 0) node[right] {$x$};
  \draw[->] (0, -0.6) -- (0, 0.6) node[above] {$y$};
  \tikzmath{\myalpha = pi/3;}
  \foreach \psiplusphi in {
  -1.04646,
  -0.62934,
  -0.43881,
  -0.28133,
  -0.13784,
   0.00000,
   0.13784,
   0.28133,
   0.43881,
   0.62934,
   1.04646
   }
{
  \draw[domain=abs(\psiplusphi)-pi:pi-abs(\psiplusphi), smooth, variable=\psiminusphi]%
  (1/2,-1/2) -- plot ({1/(2*pi)*((\psiplusphi-\psiminusphi)+sin(deg(\psiplusphi-\psiminusphi)))},%
  {1/(2*pi)*((\psiplusphi+\psiminusphi)+sin(deg(\psiplusphi+\psiminusphi)))})%
  -- (-1/2,1/2);
}
  \foreach \psiplusphi in {\myalpha-pi,pi-\myalpha}
{
  \draw[thick,domain=abs(\psiplusphi)-pi:pi-abs(\psiplusphi), smooth, variable=\psiminusphi,densely dashed]%
  plot ({-1/(2*pi)*((\psiplusphi-\psiminusphi)+sin(deg(\psiplusphi-\psiminusphi)))},%
  {1/(2*pi)*((\psiplusphi+\psiminusphi)+sin(deg(\psiplusphi+\psiminusphi)))});
}
\end{tikzpicture}%
%\quad%
\label{subfig:onethird}%
}%
\subfloat[$\alpha=\pi/2$]{%
\quad%
\begin{tikzpicture}[scale=4.8]
  \draw[->] (-0.6, 0) -- (0.6, 0) node[right] {$x$};
  \draw[->] (0, -0.6) -- (0, 0.6) node[above] {$y$};
  \foreach \pitimesu in {-1,-0.8,...,1}
{
  \tikzmath{\psiplusphi = rad(asin(\pitimesu));}
  \draw[domain=abs(\psiplusphi)-pi:pi-abs(\psiplusphi), smooth, variable=\psiminusphi]%
  (1/2,-1/2) -- plot ({1/(2*pi)*((\psiplusphi-\psiminusphi)+sin(deg(\psiplusphi-\psiminusphi)))},%
  {1/(2*pi)*((\psiplusphi+\psiminusphi)+sin(deg(\psiplusphi+\psiminusphi)))})%
  -- (-1/2,1/2);
}
  \foreach \psiplusphi in {-pi/2,pi/2}
{
  \draw[thick,domain=abs(\psiplusphi)-pi:pi-abs(\psiplusphi), smooth, variable=\psiminusphi,densely dashed]%
  plot ({-1/(2*pi)*((\psiplusphi-\psiminusphi)+sin(deg(\psiplusphi-\psiminusphi)))},%
  {1/(2*pi)*((\psiplusphi+\psiminusphi)+sin(deg(\psiplusphi+\psiminusphi)))});
}
\end{tikzpicture}%
%\quad%
\label{subfig:onehalf}%
}\\
\subfloat[$\alpha=2\pi/3$]{%
\quad%
\begin{tikzpicture}[scale=4.8]
  \draw[->] (-0.6, 0) -- (0.6, 0) node[right] {$x$};
  \draw[->] (0, -0.6) -- (0, 0.6) node[above] {$y$};
  \tikzmath{\myalpha = 2*pi/3;}
  \foreach \psiplusphi in {
  -2.09360,
  -1.19831,
  -0.82575,
  -0.52615,
  -0.25698,
   0.00000,
   0.25698,
   0.52615,
   0.82575,
   1.19831,
   2.09360
   }
{
  \draw[domain=abs(\psiplusphi)-pi:pi-abs(\psiplusphi), smooth, variable=\psiminusphi]%
  (1/2,-1/2) -- plot ({1/(2*pi)*((\psiplusphi-\psiminusphi)+sin(deg(\psiplusphi-\psiminusphi)))},%
  {1/(2*pi)*((\psiplusphi+\psiminusphi)+sin(deg(\psiplusphi+\psiminusphi)))})%
  -- (-1/2,1/2);
}
  \foreach \psiplusphi in {\myalpha-pi,pi-\myalpha}
{
  \draw[thick,domain=abs(\psiplusphi)-pi:pi-abs(\psiplusphi), smooth, variable=\psiminusphi,densely dashed]%
  plot ({-1/(2*pi)*((\psiplusphi-\psiminusphi)+sin(deg(\psiplusphi-\psiminusphi)))},%
  {1/(2*pi)*((\psiplusphi+\psiminusphi)+sin(deg(\psiplusphi+\psiminusphi)))});
}
\end{tikzpicture}%
%\quad%
\label{subfig:twothirds}%
}%
\subfloat[the family of level curves for $u$ and $v$]{%
\quad%
\begin{tikzpicture}[scale=4.8]
  \draw[->] (-0.6, 0) -- (0.6, 0) node[right] {$x$};
  \draw[->] (0, -0.6) -- (0, 0.6) node[above] {$y$};
  \foreach \pitimesu in {-1,-0.9,...,1}
{
  \tikzmath{\psiplusphi = 3*pi/4*\pitimesu;}
  \draw[domain=abs(\psiplusphi)-pi:pi-abs(\psiplusphi), smooth, variable=\psiminusphi]%
  (1/2,-1/2) -- plot ({1/(2*pi)*((\psiplusphi-\psiminusphi)+sin(deg(\psiplusphi-\psiminusphi)))},%
  {1/(2*pi)*((\psiplusphi+\psiminusphi)+sin(deg(\psiplusphi+\psiminusphi)))})%
  -- (-1/2,1/2);
  \tikzmath{\psiplusphi = 3*pi/4*\pitimesu;}
  \draw[domain=abs(\psiplusphi)-pi:pi-abs(\psiplusphi), smooth, variable=\psiminusphi,densely dashed]%
  (-1/2,-1/2) -- plot ({-1/(2*pi)*((\psiplusphi-\psiminusphi)+sin(deg(\psiplusphi-\psiminusphi)))},%
  {1/(2*pi)*((\psiplusphi+\psiminusphi)+sin(deg(\psiplusphi+\psiminusphi)))})%
  -- (1/2,1/2);
}
\end{tikzpicture}%
%\quad%
\label{subfig:uvlevelcurves}
}%
\caption[]{\subref{subfig:onethird}--\subref{subfig:twothirds} show
11 evenly distributed level curves (solid) of $u$ in \cref{th:uniform} for three
different values of $\alpha$ (corresponding to the marked points
in \cref{fig:famouslimitshape}). The lowest and highest level curves
correspond to $\psi+\phi=\pm\alpha$, and outside of these curves $u$ has a constant value of $\pm r/2$.
In the region between the dashed curves $\psi-\phi=\pm(\pi-\alpha)$, any maximizer of $\functional_\rho$ coincides with $u$.
Finally, \subref{subfig:uvlevelcurves} shows curves of the form
$\psi+\phi=\text{const}$ (solid) and $\psi-\phi=\text{const}$ (dashed). Those
are the possible level curves of $u$ and $v$, respectively, for any $0<\alpha<\pi$.}%
\label{fig:limitsurfaces}%
\end{figure}
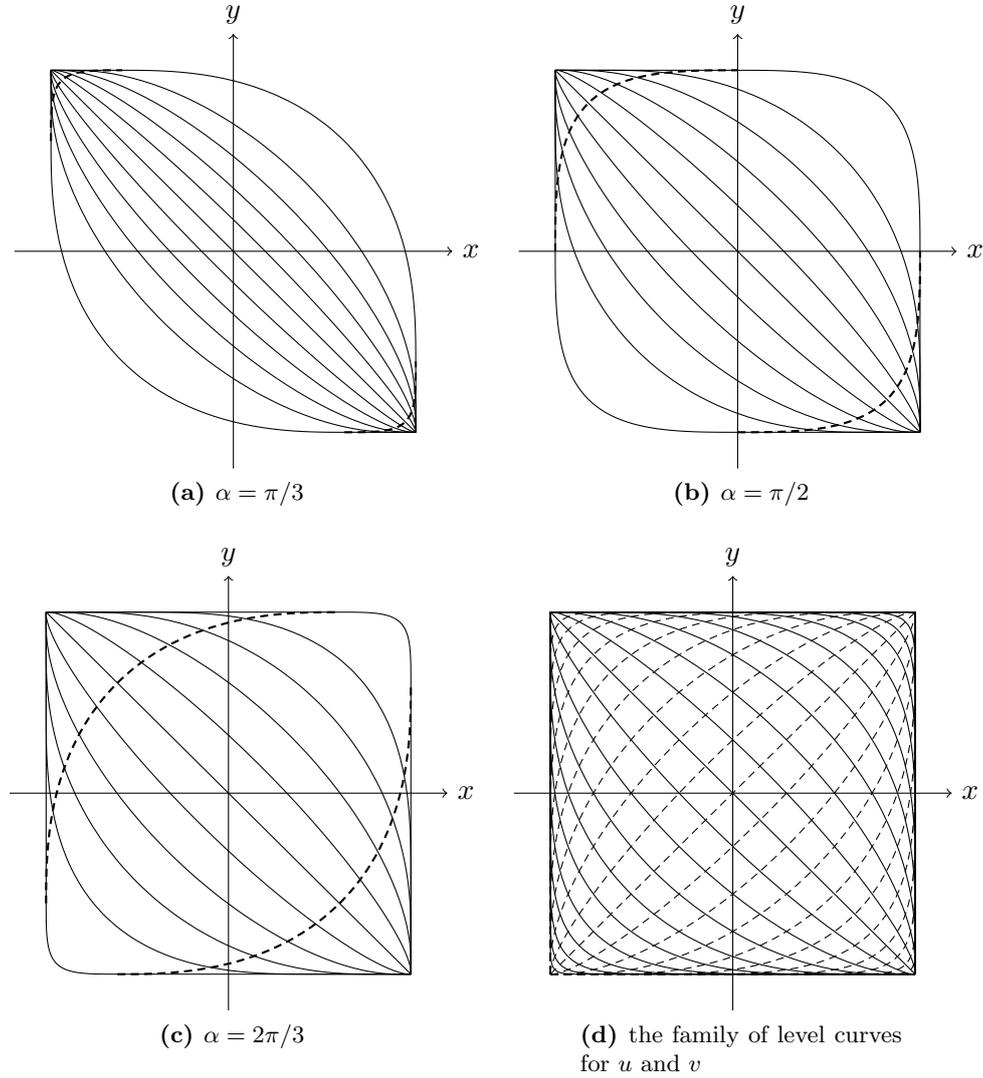
\begin{theo}\label{th:uniform}
Suppose \cref{con:triangularlimitshape} holds.
Let $\Omega$ be the open square $-\frac12<x,y<\frac12$, let $\rho=1$ and
take any $0<r<2$. Then, there is a maximizer $u$ to
$\functional_\rho$ in $\calU_{-r/2,r}(\Omega)$ given by the following.

Let $0<\alpha<\pi$ be defined by
$r = \frac{2}{\pi}(\sin\alpha-\alpha\cos\alpha)$, and define
$\phi$ and $\psi$ by
\begin{align*}
x &= \tfrac{1}{\pi}(\phi+\tfrac12\sin2\phi),\ \ \ -\tfrac\pi2<\phi<\tfrac\pi2, \\
y &= \tfrac{1}{\pi}(\psi+\tfrac12\sin2\psi),\ \ \ -\tfrac\pi2<\psi<\tfrac\pi2.
\end{align*}
Then
\[
u = \begin{cases}
\frac1\pi\bigl(\sin(\psi+\phi)-(\psi+\phi)\cos\alpha\bigr)
& \text{if $\abs{\psi+\phi} \le \alpha$,}\\
-r/2 & \text{if $\psi+\phi<-\alpha$,}\\
r/2 & \text{if $\psi+\phi>\alpha$.}\\
\end{cases}
\]
Furthermore, every maximizer to $\functional_\rho$ in $\calU_{-r/2,r}(\Omega)$
coincides with $u$ in the region where
$\abs{\psi-\phi} \le \pi-\alpha$.
\end{theo}
\begin{proof}
Recall the definition of $\calV$ from \cref{sec:generalpde}.
Let $s=r+2\cos\alpha$ and define $v\in\calV(\Omega)$ by
\[
v := \begin{cases}
\frac1\pi\bigl(\sin(\psi-\phi)+(\psi-\phi)\cos\alpha\bigr)
& \text{if $\abs{\psi-\phi}\le \pi-\alpha$,}\\
-s/2 & \text{if $\psi-\phi<-(\pi-\alpha)$,}\\
s/2 & \text{if $\psi-\phi>\pi-\alpha$.}\\
\end{cases}
\]
We claim that $u$ and $v$ satisfy the assumption of
\cref{th:pde}.
When $\psi+\phi$ spans over $(-\alpha,\alpha)$,
$u$ spans over $(-r/2,r/2)$. Analogously, when $\psi-\phi$ spans
over $(-(\pi-\alpha),\pi-\alpha)$, $v$ spans over $(-s/2,s/2)$.
From this, property \ref{it:measurers} in \cref{th:pde} follows.

In order to verify property \ref{it:phiplusphistar},
by \cref{pr:pdeunderconjecture} it is enough to check that
the equations
\begin{align}
u_xv_y+u_yv_x&=0, \label{eq:uxvyplusuyvx} \\
\min\{\lagomsqrt{u_xu_y/\rho},1\}+\min\{\lagomsqrt{-v_xv_y/\rho},1\}&=1
\label{eq:minplusmin}
\end{align}
hold almost everywhere in $\Omega$.
To this end, first we make the following simple calculations.
\begin{alignat*}{2}
u_x&=\frac{\partial u/\partial\phi}{dx/d\phi}
= \frac{\cos(\psi+\phi)-\cos\alpha}{2\cos^2\phi}
\ \ &&\text{if $\abs{\psi+\phi}\le\alpha$},\\
u_y&=\frac{\partial u/\partial\psi}{dy/d\psi}
= \frac{\cos(\psi+\phi)-\cos\alpha}{2\cos^2\psi}
\ \ &&\text{if $\abs{\psi+\phi}\le\alpha$},\\
v_x&=\frac{\partial v/\partial\phi}{dx/d\phi}
= -\frac{\cos(\psi-\phi)+\cos\alpha}{2\cos^2\phi}
\ \ &&\text{if $\abs{\psi-\phi}\le\pi-\alpha$},\\
v_y&=\frac{\partial v/\partial\psi}{dy/d\psi}
= \frac{\cos(\psi-\phi)+\cos\alpha}{2\cos^2\psi}
\ \ &&\text{if $\abs{\psi-\phi}\le\pi-\alpha$}.
\end{alignat*}
Note that
\begin{equation}\label{eq:psiphi}
\abs{\psi+\phi}+\abs{\psi-\phi}<\pi
\end{equation}
in $\Omega$.

Consider a point in $\Omega$ where $\abs{\psi+\phi}\ge\alpha$. Then $u_x=u_y=0$, so \cref{eq:uxvyplusuyvx} holds there.
Also, from \cref{eq:psiphi} it follows that
$\abs{\psi-\phi}<\pi-\alpha$
and hence $\sqrt{-v_xv_y}=\frac{\cos(\psi-\phi)+\cos\alpha}{2\cos\phi\cos\psi}=\frac{\cos(\psi-\phi)+\cos\alpha}{\cos(\psi-\phi)+\cos(\psi+\phi)}\ge1$, so \cref{eq:minplusmin} holds too.

Now consider a point where $\abs{\psi-\phi}\ge\pi-\alpha$. Then
$v_x=v_y=0$, so \cref{eq:uxvyplusuyvx} holds there.
Also, from \cref{eq:psiphi} it follows that
$\abs{\psi+\phi}<\alpha$
and hence $\sqrt{u_xu_y}=\frac{\cos(\psi+\phi)-\cos\alpha}{2\cos\phi\cos\psi}=\frac{\cos(\psi+\phi)-\cos\alpha}{\cos(\psi-\phi)+\cos(\psi+\phi)}\ge1$, so \cref{eq:minplusmin} holds too.

Finally, consider a point where $\abs{\psi+\phi}<\alpha$ and
$\abs{\psi-\phi}<\pi-\alpha$. Then \cref{eq:uxvyplusuyvx} and
\cref{eq:minplusmin} both follow
from our expressions for the partial derivatives. 

We have shown that $u$ is a maximizer of $\functional_\rho$. It remains
to be shown that every maximizer of $\functional_\rho$ in $\calU_{-r/2,r}(\Omega)$ coincides with $u$ within the region where $\abs{\psi-\phi}\le\pi-\alpha$. Let $w\in\calU_{-r/2,r}(\Omega)$ be a maximizer of $\functional_\rho$.
Let $R$ be the subregion of $\Omega$ where $\abs{\psi+\phi}<\alpha$ and
$\abs{\psi-\phi}<\pi-\alpha$.
Since $\sqrt{u_xu_y}$ and $\sqrt{-v_xv_y}$ are both positive in $R$, it follows from
\cref{eq:minplusmin} that they are both smaller than one there.
Then, by \cref{pr:uniquemaximizer}, the partial derivatives of
$u$ and $w$ coincide almost everywhere in $R$, and by \cref{lm:uniquefromderivativesgeneral}, $u$ and $w$ coincide everywhere in $R$. For any point $p$ in $\Omega$
where $\psi+\phi\le-\alpha$,
north-east of $p$ there are points in $R$ with $w$-values arbitrarily close to $-r/2$. Analogously, for any point $p$ in $\Omega$
where $\psi+\phi\ge\alpha$,
south-west of $p$ there are points in $R$ with $w$-values arbitrarily close to $r/2$. It follows that $w$ coincides with $u$ also in the
region $\abs{\psi+\phi}\ge\alpha$.
Finally, for any point $p$ in $\Omega$ where $\abs{\psi-\phi}=\pi-\alpha$, both south-west and north-east of $p$ there are
points in $R$ with $w$-values arbitrarily close to $u(p)$, so $w$ coincides with $u$ at $p$ as well.
\end{proof}

\subsection{Recovering the limit shape of Logan, Shepp and Vershik, Kerov}
\begin{figure}
\begin{tikzpicture}[scale=2.5]
  \draw[thick, domain=0:pi, smooth, variable=\myalpha]%
  plot ({2/pi*(sin(deg(\myalpha))-\myalpha*cos(deg(\myalpha)))},%
  {2/pi*(sin(deg(\myalpha))-\myalpha*cos(deg(\myalpha)))+2*cos(deg(\myalpha))});
  \foreach \myalpha in {pi/3, pi/2, 2*pi/3}
  {
    \tikzmath{\myr = 2/pi*(sin(deg(\myalpha))-\myalpha*cos(deg(\myalpha)));
    \mys = \myr + 2*cos(deg(\myalpha));}
	  \draw[fill] (\myr,\mys) circle(1pt);
  }
  \node at (0.55,1.3) {$\alpha=\pi/3$};
  \node at (1,0.7) {$\alpha=\pi/2$};
  \node at (1.6,0.3) {$\alpha=2\pi/3$};
  \draw[->] (-0.2, 0) -- (2.2, 0) node[right] {$r$};
  \draw[->] (0, -0.2) -- (0, 2.2) node[above] {$s$};
  \draw (2,0.05) -- (2,-0.05) node[below] {2};
  \draw (0.05,2) -- (-0.05,2) node[left] {2};
\end{tikzpicture}
\caption{The Logan--Shepp--$\!$Vershik--Kerov limit shape of a Young diagram
drawn from the Plancherel distribution. We have marked points
with three specific $\alpha$-values in the parameterization given by
\cref{eq:limitshapeparameterization}. Their corresponding limit surfaces
are depicted in \cref{fig:limitsurfaces}.}
\label{fig:famouslimitshape}
\end{figure}
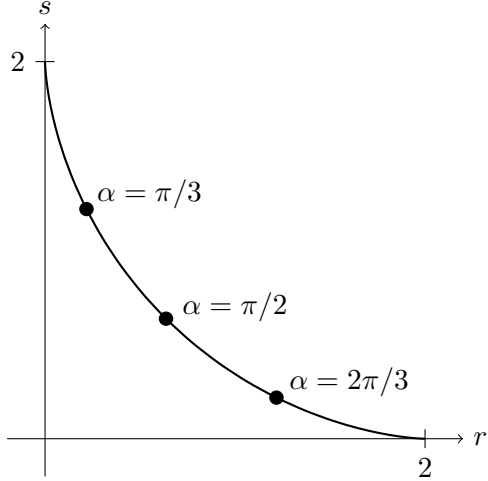

As a consequence of \cref{th:uniform}, under the assumption that
\cref{con:triangularlimitshape} holds we are able to recover the celebrated result of Logan, Shepp~\cite{LoganShepp} and Vershik, Kerov~\cite{VershikKerov} on the limit shape of the Young diagram associated with a uniformly random permutation under the Robinson--Schensted correspondence. In the proof of \cref{th:uniform} we showed that $u$ and $v$ satisfy the assumption of \cref{th:pde}. One consequence of that theorem is that
$s=\Fmax'(r)$ whenever $\Fmax'(r)$ exists. By \cref{cor:limitshape},
it follows that the limit shape in the $r$-$s$ plane is parameterized by
\begin{equation}\label{eq:limitshapeparameterization}
\begin{aligned}
r &= \tfrac2\pi(\sin\alpha-\alpha\cos\alpha), \\
s &= \tfrac2\pi(\sin\alpha-\alpha\cos\alpha) + 2\cos\alpha,
\end{aligned}
\end{equation}
where $0<\alpha<\pi$; see \cref{fig:famouslimitshape} for an illustration.
 By the substitution $\alpha=\theta+\frac\pi2$, this is equivalent to the parameterization of the Logan--Shepp--$\!$Vershik--Kerov limit shape
given in \cite[Section~1.20]{RomikBook}, namely
\begin{align*}
r &= \left(\tfrac{2\theta}{\pi}+1\right)\sin\theta+\tfrac2\pi\cos\theta, \\
s &= \left(\tfrac{2\theta}{\pi}-1\right)\sin\theta+\tfrac2\pi\cos\theta,
\end{align*}
where $-\pi/2<\theta<\pi/2$.

\subsection{Comparison with the limiting surface of a random square tableau}
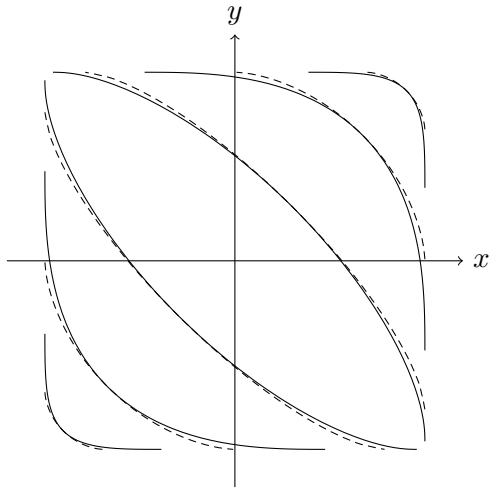
\begin{figure}
\begin{tikzpicture}[scale=5,declare function={%
		xfromphi(\phi) = 1/pi * (\phi + 1/2 * sin(deg(2*\phi)));
		s(\xminusy,\L) = sqrt(4*\L*(1-\L)-(\xminusy)^2);
		arctantwo(\a,\b) = atan(\a/\b);
		xplusy(\xminusy,\L) = 2/pi*\xminusy*rad(arctantwo((1-2*\L)*\xminusy,s(\xminusy,\L))) + 2/pi*rad(arctantwo(s(\xminusy,\L),1-2*\L)) - 1;
	}]%

  \draw[->] (-0.6, 0) -- (0.6, 0) node[right] {$x$};
  \draw[->] (0, -0.6) -- (0, 0.6) node[above] {$y$};
  \foreach \pitimesu in {-0.8, -0.5, -0.2}
{
  \tikzmath{\psiplusphi = 3*pi/4*\pitimesu;}
  \draw[domain=abs(\psiplusphi)-pi:pi-abs(\psiplusphi), smooth, variable=\psiminusphi, ]%
  plot ({xfromphi((\psiplusphi-\psiminusphi)/2)},
  	{xfromphi((\psiplusphi+\psiminusphi)/2)})
  plot ({-xfromphi((\psiplusphi-\psiminusphi)/2)},
  	{-xfromphi((\psiplusphi+\psiminusphi)/2)});
  \tikzmath{
  	\phimid = \psiplusphi / 2;
  	\xmid = xfromphi(\phimid);
  	\xpittelmid = \xmid + 1/2;
  	\L = (1 - cos(180*\xpittelmid)) / 2;
  	\xpittelstart = sqrt(1-(1-2*\L)^2);
  	\xstart = \xpittelstart - 1/2;
  }
  \draw[domain=-\xstart-1/2+0.001:\xstart+1/2-0.001, smooth, variable=\xminusy,densely dashed
  ]%
  plot ({(xplusy(\xminusy,\L)+\xminusy)/2},
	  {(xplusy(\xminusy,\L)-\xminusy)/2})
  plot ({-(xplusy(\xminusy,\L)+\xminusy)/2},
	  {-(xplusy(\xminusy,\L)-\xminusy)/2});

}
\end{tikzpicture}
\caption{Six curves of the form $\psi+\phi=\mathrm{const}$ (solid) and six level curves of the limit surface of a random square tableau (dashed) chosen so that they coincide with the solid curves at the symmetry axis $x=y$.}\label{fig:comparison}
\end{figure}

Pittel and Romik \cite{PittelRomik07}, and more recently Ma\'slanka and \'Sniady \cite{MaslankaSniady19}, studied uniform random standard Young tableaux
on an $n\times n$ square shape. After rescaling the shape to a unit square in the $x$-$y$ plane and interpreting the entries as $z$-coordinates, the tableau gives rise to a two-dimensional surface. Pittel and Romik found the limit of this surface as $n$ tends to infinity, and its level curves, shown in \cite[Fig. 1(d)]{PittelRomik07}
and \cite[Fig. 13]{MaslankaSniady19}, look very similar to the level curves of $u$ in
\cref{th:uniform}, shown in \cref{fig:limitsurfaces}. As can be seen in \cref{fig:comparison}, however, these families of curves are not the same.

\section{Future research}\label{sec:future}
The present work has generated plenty of open questions. The most significant one is \cref{con:triangularlimitshape}, of course,
and since the triangular shape is arguably the simplest
of all possible limit shapes, intuitively the phenomenon should have
a simple explanation, even though it has evaded the author so far.
As we have seen, a proof of \cref{con:triangularlimitshape} would
yield, as a by-product, a new proof of the Logan--Shepp--$\!$Vershik--Kerov limit shape in \cref{fig:famouslimitshape}, a proof very different from the known proofs.
Logan, Shepp \cite{LoganShepp} and Vershik, Kerov \cite{VershikKerov}
found the limit shape independently of
each other, but both proofs were based on the \emph{hook-length formula}, an almost magically simple formula
for computing the number of standard Young tableaux of a specific shape
(see e.g.~\cite{StanleyEnum2}).

Another category of open questions concerns the regularity of the maximizers of
the functional $\functional_\rho$. Is there always a continuous maximizer,
or even a differentiable one? Under what conditions can we find a
maximizer that satisfies the PDE system of \cref{th:pde}?
In the uniform case, provided \cref{con:triangularlimitshape} holds true, the maximizers $u$ of the form given by \cref{th:uniform}
have the following property: If $u_1$ and $u_2$ are
such maximizers associated to $r=r_1$ and $r=r_2$, respectively,
where $r_1\le r_2$, then every level curve of $u_1$ is also a level curve
of $u_2$. Is this true in general?

Our definitions of increasing and decreasing sets and density domains can be generalized
to higher dimension. For instance, we might say that a finite set of points in $\bbR^3$ is \emph{increasing} if the equivalences
$x<x'\Leftrightarrow y<y'\Leftrightarrow z<z'$ hold for any pair of points
$(x,y,z)$ and $(x',y',z')$ in the set. Can our results be generalized
to this setting?

Finally, we offer another conjecture, based on evidence from
computer-generated limit shapes for various density domains.
\begin{conj}\label{con:limitshapeconvex}
The limit shape $\Fmax'$ that appears in \cref{cor:limitshape} is
always convex.
\end{conj}
Somewhat surprisingly, we have the following implication, relying
on the result of Vershik and Kerov \cite{VershikKerov} that the
constant $\Gamma$ in \cref{th:hammersley} equals 2.
\begin{prop}
\Cref{con:limitshapeconvex} implies \cref{con:triangularlimitshape}.
\end{prop}
\begin{proof}
At the end of the proof of \cref{pr:phiproperties}, we showed that
$\Phi(r)\le\frac{\Gamma}{\sqrt2}r$ for any $r\ge0$, where $\Gamma=2$ by \cite{VershikKerov}.
Since $\Phi$ is concave, it follows that $\Phi'(r)\le\sqrt2$ for any $r\ge0$ where the derivative exists. On the other hand, by \cref{pr:phioneisone}, $\Phi'(r)=0$ for any $r>\sqrt2$.

Suppose \cref{con:limitshapeconvex} holds.
Since $\Fmax=\Phi$ for the diamond region $\abs{x}+\abs{y}<1/\sqrt2$
with $\rho=1$, it follows that $\Phi'$ is convex. Above, we saw that
$\Phi'(0)\le\sqrt2$ and $\Phi'(\sqrt2+\eps)=0$ for any small $\eps>0$,
and, together with the convexity, this yields that $\Phi'(r)\le\sqrt2-r$ for $0\le r\le \sqrt2$. But
\[
\int_0^{\sqrt2}\Phi'(r)\,dr=\Phi(\sqrt2)=1=\int_0^{\sqrt2}(\sqrt2-r)\,dr
\]
by \cref{pr:phioneisone}, so $\Phi'(r)$ must be equal to $\sqrt2-r$ in the interval $0\le r\le \sqrt2$.
\end{proof}

\section{Acknowledgement}
This work was supported by the Swedish Research Council (reg.no.~2020-04157).
The author is grateful to Prof.~Svante Jansson for
feedback on an early draft and to two anonymous referees for
insightful suggestions that helped improving the paper.

\bibliographystyle{abbrv}
\bibliography{bib-file/increasing}

\end{document}